\documentclass[final]{siamltex}

\usepackage[percent]{overpic}
\usepackage{graphicx,graphics,epsfig}
\usepackage{amsmath,amssymb,dsfont}
\usepackage{mathrsfs}
\usepackage{arydshln}
\usepackage{url,cite}
\usepackage{enumerate}
\usepackage[usenames,dvipsnames]{color}
\usepackage[normalem]{ulem}

\def\dot{\boldsymbol{\cdot}}

\def\etal{et al.~\/}

\def\fa{\hbox{ \;for all\; }}

\def\eref#1{{\rm (\ref{#1})}}
\def\pn{\par\noindent}
\def\pf{\pn{\textsc{Proof}}\ignorespaces}%
\def\qed{~\relax\ifmmode\hskip2em \Box
 \else\unskip\nobreak\hskip1em \hfill$\Box$
 \fi \newline}

\def\N{\mathbb{N}}

\def\R{\mathbb{R}}

\def\calb{\mathcal{B}}
\def\cale{\mathcal{E}}

\def\calk{\mathcal{K}}
\def\call{\mathcal{L}}
\def\caln{\mathcal{N}}

\def\calt{\mathcal{T}}
\def\calu{\mathcal{U}}

\def\Matlab{MATLAB}

\newtheorem{assumption}[theorem]{Assumption}
\newtheorem{example}[theorem]{Example}
\DeclareMathOperator*{\arginf}{arg\,inf}
\DeclareMathOperator*{\peaks}{peaks}
\DeclareMathOperator*{\franke}{franke}
\DeclareMathOperator*{\Span}{span}

\title{$H^2$--Convergence of least-squares kernel collocation methods}

\author{
Ka~Chun Cheung\thanks{Department of Mathematics, Hong Kong Baptist University, Kowloon Tong, Hong Kong.}
\and Leevan Ling\footnotemark[1]$~^,$\thanks{Correspondence to L. Ling (E-mail: {\tt lling@hkbu.edu.hk})}
\and Robert Schaback\thanks{Institut f\"ur Numerische und Angewandte Mathematik, Georg-August-Universit\"at
    G\"ottingen, Germany. }
}

\begin{document}

\maketitle

\begin{abstract}
The strong-form asymmetric kernel-based collocation method, commonly referred to
as the Kansa method, is easy to implement and hence is widely used for solving
engineering problems and partial differential equations despite the lack of
theoretical support. The simple least-squares (LS) formulation, on the other hand,
makes the  {study}
of its solvability and convergence
 rather nontrivial.
In this paper, we focus on
general second order linear elliptic differential equations in $\Omega \subset
\R^d$
{under} Dirichlet boundary conditions.
With kernels that reproduce $H^m(\Omega)$ and some smoothness assumptions
 {on the solution}, we provide  denseness
{conditions}
for a constrained least-squares  {method} and
a  {class of}  weighted least-squares
 {algorithms} to be convergent.
Theoretically,
we identify  some $H^2(\Omega)$ convergent  LS formulations
that {have}  {an} optimal error
{behavior like } $h^{m-2}$.
We also demonstrate the effects of various collocation settings
 {on the respective}
convergence rates, as well as how these
formulations perform with high order
kernels and when coupled with the stable
evaluation technique for the Gaussian kernel.

\end{abstract}



\pagestyle{myheadings}
\thispagestyle{plain}
\markboth{K.C. Cheung, L. Ling, R. Schaback}{Convergence for LS kernel collocation}

\section{Introduction}
Mathematical models or differential equations are meaningful only if
they can somehow mirror the overly complicated real world.
Similarly, numerical methods are useful only if they can produce
approximations guaranteed to converge to the outcome that the
mathematical model predicts. It could take tens of years for some
good numerical strategies to mature and become a well-established
class of numerical  {methods}
with  {a}  {complete and rigid}
theoretical framework. Take the
finite element method as an  {example. It}
waited for a quarter of a
century to get its rigorous mathematical foundation. This paper
aims to continue our theoretical contributions to the unsymmetric
radial basis  {function}
collocation  {method}, which is also known as
the Kansa method in the community and we shall use this name
throughout  {this paper} for brevity.

To quickly overview the  {development}
of the Kansa method and its connection to the radial basis
 {function (RBF)}
scattered data interpolation problem, let us
look at some of its cornerstones
 \cite{Fasshauer-Meshapprmethwith:07,%
Fasshauer+McCourt-KernApprMethusin:15,Wendland-ScatDataAppr:05}.
An RBF is a smooth scalar function $\phi:\R^+\to \R$,
which usually is induced from a kernel function
$\Phi:\R^d\times\R^d\to \R$ in today's applications,
such that the interpolant of an interpolation problem is given as a
linear combination
\begin{equation}\label{RBFexpansion}
  u = \sum_{j=1}^{n_Z} \lambda_j \phi(\| \dot - z_j \|_2)
  = \sum_{j=1}^{n_Z} \lambda_j \Phi( \dot, z_j),
\end{equation}
of shifted   RBFs
in which the set $Z=\{z_1,\ldots,z_{n_Z}\}$
 {contains} \emph{trial centers}
 {that specify the shifts of the
kernel function} in the expansion.  {Dealing with scaling}
has been
another huge topic in Kansa methods
\cite{Golbabai+MohebianfarETAL-varishapparastra:15,%
Kansa+Carlson-Impraccumultinte:92,Tsai+KolibalETAL-goldsectsearalgo:10}
for a decade, but we will ignore this point for the sake of
 {simplicity.}

Impressed by the meshfree nature, simplicity to program, dimension
independence, and  {arbitarily} high convergence rates
interpolations, {E.J.} Kansa
\cite{Kansa-Multscatdataappr:90,Kansa-Multscatdataappr:90a}
proposed  {to modify} the {RBF interpolation}
method to solve partial differential
equations (PDEs) in the early 90s. Using the same RBF expansion
\eref{RBFexpansion}, Kansa imposed strong-form
{collocation conditions}
instead of interpolation conditions for identifying
the unknown coefficients. Consider a PDE given by $\call u=f$ in
$\Omega$ and $\calb u =g$ on $\Gamma=\partial \Omega$. The
Kansa method collocates the PDE at the trial centers $Z$ to
yield exactly ${n_Z}$ conditions:
\begin{equation}\label{KansaStandard}
\begin{array}{rcl}
\call u(z_i) &=&  \sum  \lambda_j \call\phi(\| z_i- z_j\|), \quad \mbox{for $z_i\in Z\cap \Omega$},
    \\
    \calb u(z_i) &=&  \sum  \lambda_j \calb\phi( \|z_i- z_j\|), \quad \mbox{for $z_i\in Z\cap \Gamma$},
\end{array}
\end{equation}
for identifying the unknown $\lambda_j$ or equivalently,
a numerical approximation to $u$
from the \emph{trial space}
\begin{equation}\label{KansaTrialZ}
\calu_Z = \calu_{Z,\Omega,\Phi} := \Span\{ \Phi(\dot - z_j):\, z_j\in Z\}.
\end{equation}
This approach requires no re-formulation of the
PDE and no triangularization. As long as one knows how
to program for an interpolation problem,  it only takes
minutes to understand and code up something for the Kansa method.
Since invented, the Kansa method has been widely used in
vast numbers of applications in physics and engineering \cite{Chen+FanETAL-METHAPPRPARTSOLU:11,Kansa+Geiser-Numesolutimeinvi:13,Li+LiETAL-Compsuppradibasi:15,Pang+ChenETAL-SpacadveequaKans:15}.

Since the differential and boundary operators of a
PDE are independently applied to yield different rows
{of the final linear system of equations}, it is
easy to see why any Kansa {system}
matrix is unsymmetric. While this
has some implications for the choice of linear solvers, the unsymmetric
matrix {places}
the Kansa method far away from the approximation theories
from which RBFs interpolation theories were built. Though the technique
introduced by Kansa is very successful in a large variety
of  {applications in Engineerings and Science},
there were no proven results {about it} for over 10  years.
After many unsuccessful attempts to establish such a foundation,
Hon and Schaback \cite{Hon+Schaback-unsycollradibasi:01} showed in 2001
that there are extremely rare cases where the original
approach can fail because the underlying linear system can be singular.
{This put an end to all attempts to prove stability of the Kansa method
in general.}
One {workaround}
is to apply {\em symmetric {collocation}}
\cite{Fasshauer-Solvdiffequawith:99,Franke+Schaback-Solvpartdiffequa:98}
that  {mimics}
scattered Hermite  {interpolation.}
 {While the Kansa trial space basis in \eref{RBFexpansion}
is independent of the collocation, the symmetric method takes a basis that
is itself dependent on the collocation.}
This  approach
 {yields positive definite}
symmetric {system matrices}
at the expense of
{higher smoothness requirements
and less stability.} {On the positive side,
symmetric collocation can be proven \cite{schaback:2015-3} to be
  error-optimal, because it is a
pointwise optimal recovery of the solution from discrete input data.}

The situation for the Kansa method remained the same until 2006,
when we provided the first solvability results for an extended
Kansa method. In order to ensure solvability, {{\em overtesting}
is applied. Keeping the trial space \eref{KansaTrialZ}
based on a set $Z$ of {\em trial centers}, the standard Kansa system
\eref{KansaStandard} is modified by taking} {another, but usually
  larger} discrete
set {$X$} of collocation points 
that is sufficiently fine relative
to the {set $Z$ of} trial centers. 
Readers are referred to the
original articles \cite{Ling+OpferETAL-Resumeshcolltech:06}  and an
extension \cite{Schaback-ConvUnsyKernMesh:07} to the corresponding
weak problems for details.
In 2008, we had a partial answer to the convergence of an
overdetermined Kansa formulation \cite{Ling+Schaback-Stabconvunsymesh:08}.
Our analysis  {was} carried out based on the continuous and discrete
maximum norms. We  {showed}
that the $\ell^\infty$-minimizer of a residual
functional converges to the exact solution at the optimal speed,
 {i.e. with}
the same convergence rate as the interpolant converges to the exact
solution. From then on, we  {attempted}
to extend the theories to the
least-squares (LS) minimizer \cite{Kwok+Ling-convleasKansmeth:09}
and numerically  {verified} in  {extended}
precision  {arithmetic}
that  {the} LS-minimizer
also converges at {the}
optimal rate \cite{Lee+LingETAL-convnumealgounsy:09}.
Recently, in \cite{Schaback-WellProbhaveUnif:16},  we {gave} 
an $L_\infty$ convergence rate
of ${m-2-d/2}={m-3}$ for  an overdetermined Kansa method in $H^m$ for $m>3$.
In this study, we continue to work on the overdetermined Kansa method and
concentrate on the popular LS {solution}.
In Section~\ref{sec:main}, we will provide all the necessary
assumptions and  {prove} error estimates for a
{ {constrained}
least-squares} (CLS) and a class of {\em weighted least-squares}
(WLS) formulations. 
The convergence for the CLS formulation will then be given in
Section~\ref{sec:CLS}. In Section~\ref{sec:WLS} and \ref{sec:UZ},
the  {theory}
for WLS formulations in two trial spaces will be given.
Lastly, we will numerically verify the accuracy and convergence
 {rates}
of some proven convergent formulations in Section~\ref{sec:num}.

\section{Notations, assumptions and main theorems}\label{sec:main}

Throughout the paper, the notation $C_{\dot}$ will be reserved for
generic constants whose  {subscripts indicate}
the  {dependencies} of the constant.

We consider a general second order elliptic differential
 {equation} in some bounded domain $\Omega\subset \R^d$ subject to the
Dirichlet boundary condition on $\Gamma=\partial\Omega$: 
\begin{equation}\label{eq:PDE}
    \begin{array}{r@{\, = \,}ll}
      \call u & f &\mbox{ in } \Omega,\\
      u     &  g&\mbox{ on } \Gamma,
    \end{array}
\end{equation}
where
\begin{equation}\label{eq:L}
\begin{array}{rcl}
\call u
&:=&
\displaystyle{ \sum_{i,j=1}^d \frac\partial{\partial x^j}
\left( a^{ij}(x) \frac\partial{\partial x^i} u(x)\right)
+ \sum_{j=1}^d \frac\partial{\partial x^j}
\left( b^j(x) u(x) \right)}\\
&&
\quad + \displaystyle{ \sum_{i=1}^d c^i(x)
\frac\partial{\partial x^i} u(x) + d(x) u(x) = f(x)}.
\end{array}
\end{equation}
The Sobolev regularity of the true solution will be denoted by $m$, and
  we   will work with standard Hilbert spaces
$H^k(\Omega)$  and $H^{k-1/2}(\Gamma)$ with norms
$ \|u\|_{k,\Omega}$  and   $\|u\|_{k-1/2,\Gamma}$,  {respectively, for}
$k\leq m$. 

\bigskip
\begin{assumption}[%
 {Smoothness of domain and solution}]\label{ass:Omega}
We assume that the bounded domain $\Omega$ has a
 {piecewise}
$C^{m}$--boundary $\Gamma$ so that $\Omega$ is Lipschitz continuous
and satisfies  {an} interior cone condition.
Also, we assume that the functions $f$ and $g$ are smooth enough to admit  a \emph{classical} solution $u^* \in H^{m}(\Omega)$.
\qed
\end{assumption}
Now   the trace theorem \cite{Wloka-PartDiffEqua:87}
can be applied and we can define a trace operator:
\[
    \calt: H^m( \Omega ) \to H^{m-1/2}(\Gamma)
\mbox{ such that }  { \calt } u = u_{|\Gamma}
    \mbox{ for all } u \in C^{m}(\bar\Omega),
\]
for $m>1/2$, with a continuous right-inverse
linear extension operator $\mathcal{E}$ such that
\[
   \calt\circ \cale g  = g \mbox{ for all } g \in H^{m-1/2}(\Gamma).
\]
The smoothness assumption also allows {a}
partition of unity of the boundary, each  {part} of which can be
mapped to the unit ball in $\R^{d-1}$ by a $C^m$--diffeomorphism.
This allows us to define Sobolev norms on $\Gamma$ and apply
some Sobolev inequalities (i.e., kernel independent ones).

Let $\chi$ be any discrete set of $n_\chi$ points in $\Omega$.
For any $u\in H^m(\Omega)$,
 we define discrete norms on $\chi$ by 
\[
\|u\|_{\chi}= \|u\|_{0,\chi}  {=\|u\|_{\ell_2(\chi)}}, \quad
\|u\|_{k,\chi} := \Bigg(\sum_{|\alpha|\leq k} \|D^{\alpha} u\|_{\chi}^2
\Bigg)^{1/2},
 {\;0\leq k<m-d/2,\;}
\]
where $\alpha$ is some multi-{index}
and $D^\alpha u \in$  {$C(\Omega)$} are weak
derivatives of $u$.
The same notations  will also be used to denote
discrete norms on boundary for any discrete set $\chi\subset\Gamma$.

\bigskip
\begin{assumption}[Differential operator]\label{ass:L}
Assume  {that}
$\call$ as in \eref{eq:L} is a strongly elliptic operator with
coefficients belonging to $W_\infty^{m}(\Omega)$.
\qed\end{assumption}
Then, by results in \cite{Wendland:07}, $\call$ is a bounded
operator from $H^m(\Omega)$ to $H^{m-2}(\Omega)$ with
\begin{equation}\label{eq:Lu<uO}
  \|\call u\|_{m-k-2,\Omega} \leq C_{\Omega,\call} \|u\|_{m-k,\Omega},
  \;{0}\leq k \leq m-2,\, {k\in\N},
\end{equation}
for all $u\in H^m(\Omega)$.
Moreover,
the following boundary regularity estimate \cite{Jost-Partdiffequa:07} holds:
\begin{eqnarray}
  \| u\|_{{k+2},\Omega}
  &\leq& C_{\Omega,\call,k} \left( \|   {  \call u}  \|_ {k,\Omega}
+ \|  {u} \|_{k+1+1/2, {\Gamma}}  
\right),\; 0\leq k \leq m-2,
  \label{eq:bdyreg}
\end{eqnarray}
for  all $u\in H^m(\Omega)$  with $C_{\Omega,\call,k}$ depending on $\Omega$, the ellipticity constant of $\call$, and $k\geq 0$.


%

\begin{assumption}[{Kernel}]\label{ass:KRSnew}
Assume $\Phi$ is a reproducing kernel of $H^m(\Omega)$ for some integer $m\geq 2+\lceil \frac12(d+1)\rceil$.
More  {precisely},
we use a symmetric positive definite kernel $\Phi$  on $\R^d$ with smoothness $m$ that satisfies
\begin{equation}\label{eq:Fourier}
  c_{\Phi_m} ( 1+\|\omega\|_2^2)^{-m}
\leq \widehat{{\Phi_m}}(\omega) \leq C_{\Phi_m} ( 1+\|\omega\|_2^2)^{-m}
  \mbox{\quad { for all }$\omega \in \R^d$},
\end{equation}
for two constants $0< c_{\Phi_m} \leq C_{\Phi_m}$. \qed
\end{assumption}
For any $m>d/2$, its native space
$\caln_{\Omega,\Phi}$ on $\R^d$
\cite{Buhmann-Radibasifunc:03,Wendland-ScatDataAppr:05} is norm-equivalent to
$H^m(\R^d)$. This includes the  {standard}
 {{\em Whittle-Mat\'ern}}-\emph{Sobolev}
kernel with exact Fourier
transform $( 1+\|\omega\|_2^2)^{-m}$ takes the form
\[
    \Phi(x) := \|x\|_2^{m-d/2} \calk_{m-d/2}( {\|x\|_2} )\mbox{\quad for all $x\in\R^d$}, 
\]
where $\calk_\nu$ is the Bessel functions of the second kind.
The  {compactly supported piecewise polynomial}
Wendland functions \cite{Wendland-Erroestiintecomp:98} are another examples of
 {kernels} satisfying \eref{eq:Fourier}.

\begin{assumption}[Trial space]\label{ass:RSTS}
Let $Z=\{z_1,\ldots,z_{n_Z}\}$ be a discrete set of \emph{trial centers} in
$\Omega$.  {In analogy to \eref{KansaTrialZ}, but now with
  translation-invariance, we} define the  {finite-dimensional}
 trial space  {$\calu_Z$} as
\[
    \calu_Z = \calu_{Z,\Omega,\Phi} :=
\Span\{ \Phi(\dot - z_j):\, z_j\in Z\} \subset \caln_{\Omega,\Phi}.
\]
\qed
\end{assumption}
For describing the denseness of $Z\subset\Omega$, its \emph{fill distance}
 {for fixed $\Omega$} and \emph{separation distance} are defined as
\[
    h_Z := \sup_{\zeta\in\Omega}\min_{z\in Z} \| z-\zeta \|_{\ell_2(\R^d)}
\mbox{\quad and \quad}
    q_Z:=\frac12 \min_{ \scriptsize
    \begin{array}{c}
      z_i, z_j\in Z \\
      z_i\ne z_j
    \end{array}
    } \|z_i-z_j\|_{\ell_2(\R^d)},
\]
respectively, and the quantity $h_Z/q_Z=:\rho_Z$ is commonly referred as
the \emph{mesh ratio} of $Z$. For any $u$ in the native space
$\caln_{\Omega,\Phi}$ of $\Phi$, we denote $I_Zu$ to be the interpolant
of $u$ on $Z$ from the trial space $\calu_{Z}$.

\smallskip
\begin{assumption}[{Collocation points}]\label{ass:RSCP}
Let $X=\{x_1,\ldots,x_{n_X}\}$ be  {a} discrete set of
\emph{PDE collocation points} in $\Omega$ and $Y =\{y_1,\ldots,y_{n_Y}\}$
be a set of \emph{boundary collocation points} on $\Gamma$.
We assume the   {set $Z$ of} discrete trial centers  {to be}
sufficiently dense  with respect to $\Omega$, $\Phi$,
and $\call$ but independent of the solution, and the sets of points $X$
and $Z$  {to be} asymptotically quasi-uniform.  {That} is,
there exist constants  $\gamma_\chi>1$ such that
\begin{equation}\label{eq:qu}
   {\gamma_\chi}^{-1} q_{\chi} \leq h_{\chi}\leq\gamma_\chi q_{\chi}.
    \mbox{\quad { for} }\chi\in\{X,Z \}.
\end{equation}
{Note}
that the sets  {$X$ and $Y$} of collocation points 
{together} have to be as dense
as the trial centers in $Z$ to ensure stability. This paper will provide rigid sufficient conditions for this.
\qed
\end{assumption}

Imposing strong testing on \eref{eq:PDE} at collocation points in $X$ and $Y$
yields $n_X+n_Y > n_Z$ conditions, from which one can
 {hopefully} 
identify a numerical approximation from some trial spaces. The following
theorems summarize  {our}
convergence results  {for} three possible least-squares alternatives.
{The first concerns the case where we enlarge
the set $Z$ of trial points by adding
the set $Y$ of boundary
collocation points to it. Then, we can keep the numerical solution to be
exact on $Y$, and we  add this as a constraint.}

\smallskip
\begin{theorem}[{Constrained least squares}  {(CLS)}]%
\label{thm:CLS}
Suppose the Assumptions \ref{ass:Omega} to  {\ref{ass:RSCP}} hold. Let $u^*\in H^m(\Omega)$ denote the exact solution of the
elliptic PDE \eref{eq:PDE}.
In addition, the relative fill distances $ h_X/h_{Z\cup Y}$ and $ h_Z/h_{Z\cup
  Y}$
are sufficiently small and satisfy condition \eref{eqZYcond}.
Let $u^{CLS}_{X,Y}\in\calu_{Z\cup Y}$ be the constrained
least-squares solution defined as
\begin{equation}\label{eq:CLSdef}
    u^{CLS}_{X,Y}:= \underset{
                u\in\calu_{Z\cup Y}
    }{\arginf} \| \call u - f \|^2_{X}
    \mbox{\quad subject to } u_{|Y} = g_{|Y}.
\end{equation}
Then
the  error estimates
  \[
\|u^{CLS}_{X,Y}-u^*\|_{2,\Omega}\le
C_{\Omega,\Phi,\call,\gamma_X}
h_{Z\cup Y}^{m-d/2-2}\|u^*\|_{m,\Omega} \mbox{\quad for $m\geq 2+ \left\lceil\frac{d+1}2\right\rceil$},
  \]
and
  \[
\|u^{CLS}_{X,Y}-u^*\|_{2,\Omega}\le
C_{\Omega,\Phi,\call,\gamma_X}
h_{Z\cup Y}^{m-2}\|u^*\|_{m,\Omega} \mbox{\quad for $m> 3+\frac{d}2$},
  \]
hold for some constant $C_{\Omega,\Phi,\call,\gamma_X}$
that depends only on $\Omega$, $\Phi$, $\call$, and the
uniformity constant $\gamma_X$ of $X$.
\end{theorem}

\bigskip
 {The next case does not require exactness on $Y$ but still
keeps $Z\cup Y$
as the set of trial centers.}
\smallskip

\begin{theorem}[{Weighted least squares}  {(WLS)}]%
\label{thm:WLS}
Suppose all the assumptions in Theorem~\ref{thm:CLS} hold.
Let
$u^{WLS,\theta}_{X,Y,Z\cup Y}\in\calu_{Z\cup Y}$
be the weighted least-squares solution defined as
\begin{equation}\label{eq:WLSdef}
    u^{WLS,\theta}_{X,Y,Z\cup Y}:=
\underset{
                u\in\calu_{Z\cup Y}
    }{\arginf} \| \call u - f \|^2_{X} +
\left(\frac{h_Y}{h_X}\right)^{d\theta/2}h_Y^{-2\theta} \| u - g \|^2_{Y}
   \mbox{ \quad for  $ \theta\geq0$}.
\end{equation}
Then, for  $h_X \leq h_Y<1$ and $0\leq \theta \leq 2$, the error estimate
\begin{eqnarray*}
   \| u^{WLS,\theta}_{X,Y,Z\cup Y}- u^* \|_{2,\Omega}
   &\leq&
C_{\Omega,\Phi, \call,\gamma_X } \Big( 1+
  h_X^{\frac{(\theta-2)d}4} h_Y^{\frac{(\theta-2)(d-4)}4}\Big)
h_{Z \cup Y}^{m-d/2-2}\|u^*\|_{m,\Omega}
\end{eqnarray*}
{for $m\geq 2+ \left\lceil\frac{d+1}2\right\rceil$}, and
\begin{eqnarray*}
   \| u^{WLS,\theta}_{X,Y,Z\cup Y}- u^* \|_{2,\Omega}
   &\leq&
C_{\Omega,\Phi, \call,\gamma_X } \Big( 1+
  h_X^{\frac{(\theta-2)d}4} h_Y^{\frac{(\theta-2)(d-4)}4}\Big)
h_{Z \cup Y}^{m-2}\|u^*\|_{m,\Omega}
\end{eqnarray*}
for $m> 3+\frac{d}2$,
hold for some constant
$C_{\Omega,\Phi,\call}$ that depends only on $\Omega$, $\Phi$, and $\call$. For $2\leq \theta\leq\infty$, the estimates in Theorem~\ref{thm:CLS} remain valid.
\end{theorem}
\bigskip

 {Finally,   we go back to the case where $Z$ is the set of trial nodes,
independent of $X$ and $Y$}.
\smallskip
\begin{theorem}[WLS in a smaller trial space]\label{cor:WLS}
  Suppose the trial space of the weighted least-squares
approximation in Theorem
  \ref{thm:WLS} is restricted to $u^{WLS,\theta}_{X,Y,Z}\in\calu_Z$ instead of
  $\calu_{Z\cup Y}$.
Moreover, the relative fill distances $ h_X/h_{Z\cup Y}$ and $ h_Z/h_{Z\cup Y}$ are sufficiently small and satisfy condition \eref{eqZYcond2} instead of \eref{eqZYcond}.
Further assume that the sets $Y$ are asymptotically
quasi-uniform with constant $\gamma_Y>1$ as in \eref{eq:qu} and $h_Y\leq h_Z$.
  Then, for any $0\leq \theta \leq 2$,  the error estimates
  \begin{eqnarray*}
 && \| u^{WLS,\theta}_{X,Y,Z}- u^* \|_{2,\Omega}
 \leq
 C_{\Omega,\Phi, \call,\vec\gamma}
 \Big(1+
  h_X^{\frac{(\theta-2)d}4} h_Y^{\frac{(\theta-2)(d-4)}4}
    +
    h_Y^{-3/2} h_Z^2  \Big)
    h_Z^{m-d/2-2}\|u^*\|_{m,\Omega}
\end{eqnarray*}

for $m\geq 2+\left\lceil\frac{d+1}{2}\right\rceil$, and
  \begin{eqnarray*}
 && \| u^{WLS,\theta}_{X,Y,Z}- u^* \|_{2,\Omega}
 \leq
 C_{\Omega,\Phi, \call,\vec\gamma}
 \Big(1+
  h_X^{\frac{(\theta-2)d}4} h_Y^{\frac{(\theta-2)(d-4)}4}
    +
    h_Y^{-2} h_Z^2  \Big)
    h_Z^{m-2}\|u^*\|_{m,\Omega}
\end{eqnarray*}
for   $m>3+d/2$,  hold for  some constant
$C_{\Omega,\Phi,\call,\vec\gamma}$ that depends only on
$\Omega$, $\Phi$,  $\call$, and uniformity
constants $\vec\gamma=[\gamma_X,\gamma_Y,\gamma_Z]$.
\end{theorem}
\section{Optimal convergence rates for CLS}\label{sec:CLS}
We first prove some necessary inequalities essential to our proofs.

\smallskip
\begin{lemma}[Sampling Inequality of fractional order]\label{lem:sampling}
  Suppose $\Omega\subset\R^d$ is a bounded  {Lipschitz}
domain with a  {piecewise} $C^m$--boundary.
  Then, there exists positive constant $C_{\Omega, m, s} $
depending on $\Omega$, $m$ and $s$ such that the following holds:
  \[
  \|u\|_{s,\Omega}\le C_{\Omega, m, s}  \left(h_X^{m-s}\|u\|_{m,\Omega}
+ h_X^{d/2-s}\|u\|_{X}  \right) \mbox{\quad for $ 0\leq s\leq m$},
  \]
  and
  \[
  \|u\|_{s-1/2,\Gamma}\le C_{\Omega, m, s}  \left(h_Y^{m-s}\|u\|_{m,\Omega}
+ h_Y^{d/2-s}\|u\|_{Y}  \right) \mbox{\quad for $ 1/2\leq s\leq  m$},
  \]
  for any $u\in H^{m}(\Omega)$ with $m>d/2$ and any discrete sets
$X\subset \Omega$ and $Y\subset \Gamma$  with sufficiently small
mesh norm $h_X$ and $h_Y$.
\end{lemma}
\smallskip\pf. The interior sampling inequality for $X\subset\Omega$,
which only requires $\Omega$ be a bounded Lipschitz domain, is a
special case of a sampling inequality in \cite{Arcangeli:12}.
Applying the interior sampling inequality to the union of unit
balls in $\R^{d-1}$, which are images of the partition of unity of
$\Gamma$ under the {$C^m$}--diffeomorphism  in
Assumption~\ref{ass:Omega},  yields
  \[
  \|u\|_{s-1/2,\Gamma}\le C  \left(h_Y^{(m-1/2)-(s-1/2)}\|u\|_{m-1/2,\Gamma}
+ h_Y^{(d-1)/2-(s-1/2)}\|u\|_{Y}  \right),
  \]
  for all $1/2\leq s\leq  m$.
  Finally, by applying the trace theorem, the desired boundary
sampling inequality is obtained.
\qed

\begin{lemma}[Inverse Inequality]
\label{lem:invineq}
Let a kernel ${\Phi_m}\;:\;\R^d\times\R^d\to \R$ satisfying \eref{eq:Fourier}
with smoothness $m>d/2$ be given.
Suppose $\Omega\subset\R^d$ is a bounded Lipschitz domain satisfying
{an} interior cone condition.
Assume $0\leq \nu\leq m-d/2$ and {$d/2< \sigma\leq m-2\nu$}
for some {integers $\nu$ and positive $\sigma$.}
Then there is a constant
{$C_{\Omega,\Phi_m,  \sigma, \nu } $},
depending only on
{$\Omega, \Phi_m,  \sigma$, and $\nu$} such that
\begin{equation}\label{eq:invineq}
  \|u\|_{{\sigma+2\nu},\Omega}\le
{C_{\Omega,\Phi_m,  \sigma,\nu }}
 {\rho_Z^{m-\nu} h_Z^{-{\sigma}} \|u\|_{2\nu,\Omega}}
\mbox{ \quad for all $u \in \calu_Z$ }
\end{equation}
  holds  {in the trial
space of {$\Phi_m$} on  {all} sufficiently dense and
quasi-uniform sets $Z\subset \Omega$} with fill
distance $h_Z$ and mesh ratio $\rho_Z$.
\end{lemma}
%
\smallskip\pf.
The basic proof idea is to use an inverse inequality from
\cite[Eqn. 3.19]{Hangelbroek+NarcowichETAL-invetheocompLips:15}.
It has the $L_2(\Omega)$ norm on the right-hand side, but for
\eref{eq:invineq} we have to go over to derivatives there.
The idea is to push the derivatives into a new kernel.

Let $\Phi_m$ denote the given kernel satisfying
\eref{eq:Fourier} with parameter $m$.
For all $0\leq \nu<m-d/2$, we
define symmetric positive definite kernels
$\Psi_{m-\nu}:=(I-\Delta)^\nu \Phi_m$, whose Fourier transforms satisfy
\[
  c_1 ( 1+\|\omega\|_2^2)^{-(m-\nu)} \leq
\widehat{\Psi}_{m-\nu}(\omega)=(1+\|\omega\|_2^2)^{\nu}\widehat{\Phi}_{m}(\omega)
 \leq C_1 ( 1+\|\omega\|_2^2)^{-(m-\nu)},
\]
and hence, $\Psi_{m-\nu}$ has behavior like $\Phi_{m-\nu}$. This is like
applying the operator $( 1-\Delta)^{\nu/2}$ to both arguments of $\Phi_m$,
if $\Phi_m$ is written in difference form.

We use the notation
$u_{\beta,Z,\Phi_m}  := \sum_{z_j\in Z} \beta_j \Phi_m(\dot - z_j)$
to denote the
{functions} in the trial
space {$\calu_{Z,\Phi_m}$ spanned by translates}
of the kernel $\Phi_m$ on {the} trial centers
{in} $Z$ with coefficients {forming a vector}
{$\beta\in\R^{|Z|}$}.
Then
\[
\begin{array}{rcl}
(I-\Delta)^\nu u_{\beta,Z,\Phi_m}
&=&u_{\beta,Z,\Psi_{m-\nu}}
\end{array}
\]
holds, and these are the functions that we
use in \cite[Eqn. 3.19]{Hangelbroek+NarcowichETAL-invetheocompLips:15}.
This yields
\[
\begin{array}{rcl}
\big\|  u_{\beta,Z,\Psi_{m-\nu}}
\big\|_{\sigma,\Omega}
   &\leq&  C_{\Omega,\Phi_m,\mu,\sigma} \,\rho_Z^{m-\nu} h_Z^{-\sigma}
   \big\|  u_{\beta,Z,\Psi_{m-\nu}}  \big\|_{0,\Omega},
\end{array}
\]
for all $\beta\in \R^{|Z|}$ and $0\leq \nu<m-d/2,\;0\leq \sigma\leq m-\nu$.
The right hand side is what we want, because of
\begin{eqnarray*}
  \|u_{\beta,Z,\Psi_{m-\nu}}\|_{0,\Omega}^2 &=&  \int_\Omega
\left|(I-\Delta)^\nu u_{\beta,Z,\Phi_{m}} \right|^2 \;dx\\
  &\leq& C_{\Omega,\nu} \sum_{|\alpha|\leq 2\nu }  \int_\Omega
\left|  D^{\alpha} u_{\beta,Z,\Phi_{m}} \right|^2 \;dx\\
  &=&  C_{\Omega,\nu}  \|u_{\beta,Z,\Phi_{m}}  \|_{2\nu,\Omega}^2,
\end{eqnarray*}
for $2\nu\leq m$.
We are now left with
\[
\begin{array}{rcl}
\big\|  u_{\beta,Z,\Psi_{m-\nu}}
\big\|_{\sigma,\Omega}^2
&=&
\|(I-\Delta)^\nu u_{\beta,Z,\Phi_{m}}\|_{\sigma,\Omega}^2\\
\end{array}
\]
and our goal is to bound this from
below by $\|u_{\beta,Z,\Phi_{m}}\|^2_{\sigma+2\nu,\Omega}$.
By definition via Fourier transforms,
\[
{H}^k_2(\R^d) =
\left\{v\;: \; (1+\|\omega\|_2^2)^{k/2}\widehat{v}(\omega) \in L^2(\R^d)\right\}.
\]
{We use the sloppy notation $\Phi_k$ to denote kernels in the family of $\Phi_m$ with smoothness $d/2<k\leq m$.}
Then, ${H}^k_2(\R^d)$ is norm equivalent to the native space
$N_{\Phi_k}(\R^d)$ of ${\Phi_k}$ on $\R^d$
and equal as sets; i.e., 
\[
 c_{\Phi_m,k}\|u\|_{N_{\Phi_k}(\R^d)} \leq \|u\|_{H_2^k(\R^d)}\leq C_{\Phi_m,k}\|u\|_{N_{\Phi_k}(\R^d)}
\]
for some constants $0<c_{\Phi_m,k}\leq C_{\Phi_m,k}$.
Before we go over to subdomains, we note that the Fourier transform of the
operator $I-\Delta$ is $1+\|\omega\|_2^2$, and this extends to arbitrary
non-integer powers. Thus,
\[
(I-\Delta)^{\nu}{H}^{k{+2\nu}}_2(\R^d)={H}^{k}_2(\R^d)
\mbox{\quad for $0\leq 2\nu\leq m-k$},
\]
and
\[
\|(I-\Delta)^{\nu}u\|_{H_2^{k}(\R^d)}= \|u\|_{H_2^{k+2\nu}(\R^d)}
\mbox{\quad for all $u\in H_2^{k+2\nu}(\R^d), \;0\leq
2\nu\leq m-k$.}
\]

For a domain $\Omega\subset\R^d$ with
Lipschitz boundary, we also have
that $N_{\Phi_k}(\Omega)$ and ${H}^k_2(\Omega)$
are equal as sets and the norms are equivalent
\cite[Cor 10.48]{Wendland-ScatDataAppr:05}, where $H_2^k(\Omega)$ has the
standard definition via weak derivatives and $N_{\Phi_k}(\Omega)$
has the standard definition via
a closure of ${\Phi_k}$-translates. We use the sloppy notation
\[
 c_{\Omega,\Phi_m,k}\|u\|_{N_{\Phi_k}(\Omega)}
 \leq \|u\|_{H_2^k(\Omega)}\leq
 C_{\Omega,\Phi_m,k}\|u\|_{N_{\Phi_k}(\Omega)}
\]
for some constants $0<c_{\Omega,\Phi_m,k}\leq C_{\Omega,\Phi_m,k}$. By
Theorem 10.47 there, the
restriction operator $R_{{\Phi_k},\Omega}\;:\;N_{\Phi_k}(\R^d)\to N_{\Phi_k}(\Omega)$
is well-defined  and
satisfies
\[
\|R_{{\Phi_k},\Omega}v\|_{N_{\Phi_k}(\Omega)}\leq \|v\|_{N_{\Phi_k}(\R^d)} \mbox{\quad for all $ v\in N_{\Phi_k}(\R^d)$}.
\]
Going the other way, there is an isometric extension operator
$E_{{\Phi_k},\Omega}\;:\;N_{\Phi_k}(\Omega)\to N_{\Phi_k}(\R^d)$
\cite[Th. 10.46]{Wendland-ScatDataAppr:05}.
Most of this can already be found in \cite{schaback:1997-3}.
Another extension operator is
$
E_{k,\Omega}\;:\;H_2^k(\Omega) \to H_2^k(\R^d),
$
and it is bounded. This takes into account that the global space is
defined via Fourier transforms, while the local one has $L_2$ integrals
over weak derivatives.

It is known that the Sobolev extension operators do not commute
with general derivatives. On the trial functions, we do have  $E_{\Phi_\sigma,\Omega}(I-\Delta)^\nu
=(I-\Delta)^\nu E_{\Phi_{\sigma+2\nu},\Omega}$ for $m-\nu>d/2$ and  $0\leq \sigma+2\nu \leq m$ since
\begin{equation}\label{eq:ExtId}
\begin{array}{rcl}
E_{\Phi_\sigma,\Omega}(I-\Delta)^\nu u_{\beta,Z,\Phi_{m}}
&=&
E_{\Phi_\sigma,\Omega}u_{\beta,Z,\Psi_{m-\nu}}\\
&=&
u_{\beta,Z,\Psi_{m-\nu}}\\
&=&
(I-\Delta)^\nu u_{\beta,Z,\Phi_{m}}\\
&=&
(I-\Delta)^\nu E_{\Phi_{\sigma+2\nu},\Omega}u_{\beta,Z,\Phi_{m}}\\
\end{array}
\end{equation}
hold if the functions lie in the correct
space, i.e.,
\[
\begin{array}{rcl}
E_{\Phi_\sigma,\Omega}u_{\beta,Z,\Psi_{m-\nu}}
=
u_{\beta,Z,\Psi_{m-\nu}},
&\mbox{ and }&
E_{\Phi_{\sigma+2\nu},\Omega}u_{\beta,Z,\Phi_{m}}
=u_{\beta,Z,\Phi_{m}}.
\end{array}
\]
The functions
are globally defined anyway, and thus they coincide with their
extension if the global norms are bounded. Thus, we need that
\[
\begin{array}{rcl}
u_{\beta,Z,\Psi_{m-\nu}} \in N_{\Phi_\sigma}(\R^d),
&\mbox{ and }&
u_{\beta,Z,\Phi_{m}} \in N_{\Phi_{\sigma+2\nu}}(\R^d).
\end{array}
\]
The condition for the first case
is $2(m-\nu)-\sigma>d/2$.
The second case requires a finite $d$-variate integral over
\[
|\widehat \Phi_m(\omega)|^2(1+\|\omega\|_2^2)^{\sigma+2\nu}
=(1+\|\omega\|_2^2)^{-2m+\sigma+2\nu},
\]
that yields the same condition.  
%
%
Using the native space extension operators, we get
\[
\begin{array}{rcl}
\|(I-\Delta)^\nu u_{\beta,Z,\Phi_{m}}\|_{H_2^\sigma(\Omega)}
&\geq&
C_{\Omega,\Phi_m,\sigma}
\|(I-\Delta)^\nu u_{\beta,Z,\Phi_{m}}\|_{N_{\Phi_\sigma}(\Omega)}\\
&= &
C_{\Omega,\Phi_m,\sigma}
\|E_{\Phi_\sigma,\Omega}(I-\Delta)^\nu u_{\beta,Z,\Phi_{m}}\|_{N_{\Phi_\sigma}(\R^d)}\\
&\geq &
C_{\Omega,\Phi_m,\sigma}^\prime
\|E_{\Phi_\sigma,\Omega}(I-\Delta)^\nu u_{\beta,Z,\Phi_{m}}\|_{H_2^\sigma(\R^d)},
\end{array}
\]
%

and, by the extension identity \eref{eq:ExtId},
\[
\begin{array}{rcl}
\|(I-\Delta)^\nu u_{\beta,Z,\Phi_{m}}\|_{H_2^\sigma(\Omega)}
&\geq&
C_{\Omega,\Phi_m, \sigma}^\prime
\|(I-\Delta)^\nu E_{\Phi_{\sigma+2\nu},\Omega}u_{\beta,Z,\Phi_{m}}\|_{H_2^\sigma(\R^d)}\\
&=&
C_{\Omega,\Phi_m, \sigma}^\prime
\|E_{\Phi_{\sigma+2\nu},\Omega}u_{\beta,Z,\Phi_{m}}\|_{H_2^{\sigma+2\nu}(\R^d)}.
\end{array}
\]
We now go local by
\[
\begin{array}{rcl}
\|(I-\Delta)^\nu u_{\beta,Z,\Phi_{m}}\|_{H_2^\sigma(\Omega)}
&\geq &
C_{\Omega,\Phi_m, \sigma,\nu}
\|E_{\Phi_{\sigma+2\nu},\Omega}u_{\beta,Z,\Phi_{m}}\|_{N_{\Phi_{\sigma+2\nu}}(\R^d)}\\
&= &
C_{\Omega,\Phi_m, \sigma,\nu}
\|u_{\beta,Z,\Phi_{m}}\|_{N_{\Phi_{\sigma+2\nu}}(\Omega)}\\
&\geq &
C_{\Omega,\Phi_m, \sigma,\nu}^\prime
\|u_{\beta,Z,\Phi_{m}}\|_{H_2^{\sigma+2\nu}(\Omega)}.\\
\end{array}
\]
We completed proving a local Bernstein inequality; note
that  the weaker norm on the right hand side must take an
even order, whereas the global counterpart
\cite{Narcowich+WardETAL-SoboErroEstiBern:06} allows
any nonnegative integer orders.
\qed
\begin{lemma}[$H^2$--Stability]\label{lem:CLSstab}
Let a kernel $\Phi$ as in  \eref{eq:Fourier}
with smoothness
{$m\geq 2$ and $m>d/2$ be given.}
Suppose $\Omega\subset\R^d$ is a bounded Lipschitz
domain satisfying  {an} interior cone condition.
If the elliptic operator $\call$ satisfies all
assumptions to allow regularity \eref{eq:bdyreg},
then  there exists a constant $C_{\Omega, \call} $,
depending only on  $\Omega$, $\Phi$, and $\call$ such that
\[
   \| u \|_{2,\Omega} \leq C_{\Omega, \call}
\big( h_X^{d/2} \|\call u\|_{X} + h_Y^{d/2-2} \|u\|_{Y} \big)
\]
holds in two circumstances:
\begin{itemize}
 \item for all $u\in\calu_{Z\cup Y}$ under the condition
\begin{equation}\label{eqZYcond}
C_{\Omega, \call,m}(h_X^{m-2}+h_Y^{m-2})h_{Z\cup Y}^{-m+2}\rho_{Z\cup Y}^{m-1}<1/2,
\end{equation}
\item or for all $u\in\calu_{Z}$ under the condition
\begin{equation}\label{eqZYcond2}
C_{\Omega, \call,m}(h_X^{m-2}+h_Y^{m-2})h_Z^{-m+2}\rho_Z^{m-1}<1/2.
\end{equation}
\end{itemize}
\end{lemma}
\par\smallskip\noindent
Note that the factor at the second term
in the assertion is not $h_Y^{(d-1)/2}$ as one
  would expect. This might be connected to the fact that the natural norm on the
boundary is the $L_\infty$ norm,
due to the Maximum Principle.
\smallskip\pf.
We apply the first inequality of Lemma~\ref{lem:sampling} (for $u=\call u,\,s=0,m=m-2)$
to get
$$
  \|\call u\|_{0,\Omega}\le C_{\Omega, m-2, 0}  \left(
h_X^{m-2}\|\call u\|_{m-2,\Omega}
+ h_X^{d/2}\|\call u\|_{X}  \right) \fa u\in H^m(\Omega)
$$
and, by \eref{eq:Lu<uO},
$$
  \|\call u\|_{0,\Omega}\le C_{\Omega, m-2, 0}  \left(
h_X^{m-2}\|u\|_{m,\Omega}
+ h_X^{d/2}\|\call u\|_{X}  \right) \fa u\in H^m(\Omega).
$$
Using the second inequality of Lemma~\ref{lem:sampling} (for $s=2,m$), we get
$$
  \|u\|_{1+1/2,\Gamma}\le C_{\Omega, m, 2}  \left(h_Y^{m-2}\|u\|_{m,\Omega}
+ h_Y^{d/2-2}\|u\|_{Y}  \right) \fa u\in H^m(\Omega)
$$
and this combines with the
$H^2$ regularity estimate \eref{eq:bdyreg}:
\[
    \|u\|_{2,\Omega} \leq  C_{\Omega, \call}( \|\call u\|_{0,\Omega}
    + \| u\|_{1+1/2,\Gamma})
\]
into
$$
\|u\|_{2,\Omega} \leq  C_{\Omega, \call,m}
\left(
(h_X^{m-2}+h_Y^{m-2})\|u\|_{m,\Omega}
+ h_X^{d/2}\|\call u\|_{X}
+ h_Y^{d/2-2}\|u\|_{Y}
\right).
$$
Up to here, we are still in full Sobolev space. Now we use
the inverse inequality, whatever the trial space is.
If we only take $Z$ nodes like in the lemma, then
$$
\|u\|_{2,\Omega} \leq  C_{\Omega, \call,m}
\left(
(h_X^{m-2}+h_Y^{m-2})h_Z^{-m+2}\rho_Z^{m-1}\|u\|_{2,\Omega}
+ h_X^{d/2}\|\call u\|_{X}
+ h_Y^{d/2-2}\|u\|_{Y}
\right)
$$
for all $u\in \calu_{Z}$ and the $H^2$ stability is
$$
\|u\|_{2,\Omega} \leq  C_{\Omega, \call,m}
\left(h_X^{d/2}\|\call u\|_{X}
+ h_Y^{d/2-2}\|u\|_{Y}
\right)
$$
for all $u\in \calu_{Z}$ under the condition
$$
C_{\Omega, \call,m}(h_X^{m-2}+h_Y^{m-2})h_Z^{-m+2}\rho_Z^{m-1}<1/2.
$$
If we now take $Z\cup Y$ nodes, then
$$
\|u\|_{2,\Omega} \leq  C_{\Omega, \call,m}
\left(
(h_X^{m-2}+h_Y^{m-2})h_{Z\cup Y}^{-m+2}\rho_{Z\cup Y}^{m-1}\|u\|_{2,\Omega}
+ h_X^{d/2}\|\call u\|_{X}
+ h_Y^{d/2-2}\|u\|_{Y} \right)
$$
for all $u\in \calu_{Z\cup Y}$ and the $H^2$ stability is
$$
\|u\|_{2,\Omega} \leq  C_{\Omega, \call,m}
\left(h_X^{d/2}\|\call u\|_{X}
+ h_Y^{d/2-2}\|u\|_{Y}
\right)
$$
for all $u\in \calu_{Z\cup Y}$ under the condition
$$
C_{\Omega, \call,m}(h_X^{m-2}+h_Y^{m-2})h_{Z\cup Y}^{-m+2}\rho_{Z\cup Y}^{m-1}<1/2.
$$
\qed

\begin{lemma}[Consistency]\label{lem:CLSconsis}
If the elliptic operator $\call$ satisfies Assumption~\ref{ass:L} and if the kernel satisfies Assumption~\ref{ass:KRSnew}, we have
\[
    \min_{
    \scriptsize
    \begin{array}{c}
            v\in\calu_{Z\cup Y}  \\
            v_{|Y}={u^*}_{|Y}
          \end{array}
    } {\|\call v- \call u^*\|_X} \leq
    C_{\Omega,\Phi, \call}
    {n_X^{1/2}}
    h_{Z\cup Y}^{m-2-d/2} \|u^*\|_{m,\Omega}
    \mbox{\quad for $m\geq 2+\left\lceil \frac{d+1}2\right\rceil$},
\]
and
\[
    \min_{
    \scriptsize
    \begin{array}{c}
            v\in\calu_{Z\cup Y}  \\
            v_{|Y}={u^*}_{|Y}
          \end{array}
    } {\|\call v- \call u^*\|_X} \leq
    C_{\Omega,\Phi, \call}
    {n_X^{1/2}\rho_X^{d/2}}
    h_{Z\cup Y}^{m-2} \|u^*\|_{m,\Omega}
    \mbox{\quad for $m> 3+d/2$},
\]
for any $u^*\in H^m(\Omega)$.
\end{lemma}


\smallskip\pf.
By comparing the minimizer $v^*\in\calu_{Z\cup Y}$ of the optimization problem with the interpolant $I_{Z\cup Y}u^* \in \calu_{Z\cup Y}$ that also satisfies the constraints at $Y$, we turn the problem into an error estimate for radial basis function interpolation:
\begin{eqnarray*}
{\|\call v^*- \call u^*\|_X}
    &\leq&  \|\call I_{Z\cup Y}u^* - \call u^* \|_{X}.
\end{eqnarray*}
The first error estimate can be derived based on native space error estimates \cite[Thm.11.9]{Wendland-ScatDataAppr:05} and upper bounds of power functions \cite[Sec.15.1.2]{Fasshauer-Meshapprmethwith:07}.
For $m\geq 2+\left\lceil \frac{d+1}2\right\rceil$, we have
\begin{eqnarray*}
    \|\call I_{Z\cup Y}u^* - \call u^* \|_{X}
    &\leq& n_X^{1/2} \|\call I_{Z\cup Y}u^* - \call u^* \|_{L^\infty(\Omega)}
\\
  &\leq& C_{\Omega, \call}  n_X^{1/2} \max_{|\alpha|\leq 2} |D^{\alpha} I_{Z\cup Y}u^* - D^{\alpha}u^*|
\\
&\leq& C_{\Omega,\Phi, \call} n_X^{1/2} h_{Z\cup Y}^{m-d/2-2} \|u^*\|_{m,\Omega}.
\end{eqnarray*}

If we employ kernels with a higher smoothness parameter $m>3+d/2$, we can use the estimates for functions with scattered zeros.
Applying \cite[Prop.3.3]{Narcowich:05}
to our
Hilbert space setting and taking care of the
definitions of discrete {norms}
yield the desired  error bound.
\qed

To prove Theorem~\ref{thm:CLS}, suppose Assumptions \ref{ass:Omega}--
{\ref{ass:RSCP}} hold so that all lemmas in this section can be applied.
Let $u^{CLS}_{X,Y}\in\calu_{Z\cup Y}$  be the CLS approximation of \eref{eq:PDE}, defined as in \eref{eq:CLSdef}.
Moreover, let $I_{Z\cup Y}u^*$ denote the unique interpolant of
the exact solution
{$u^*\in  H^m(\Omega)$
from the  trial space $\calu_{Z\cup Y}\subset\caln_{\Omega,\Phi}= H^m(\Omega)$.}
Assume the condition \eref{eqZYcond} holds;
 we shall show that the CLS solution converges to the interpolant in $\calu_{Z\cup Y}$.
\begin{eqnarray}
    \|u^{CLS}_{X,Y}- u^* \|_{2,\Omega}
    &\leq& \|u^{CLS}_{X,Y}-   I_{Z\cup Y}u^* \|_{2,\Omega} + \|  I_{Z\cup Y}u^* - u^* \|_{2,\Omega} \nonumber
    \\
    &\leq& \|u^{CLS}_{X,Y}-   I_{Z\cup Y}u^* \|_{2,\Omega} + C_{\Omega, \Phi, \call} h_{Z\cup Y}^{m-2} \|u^*\|_{m,\Omega}, \nonumber
\end{eqnarray}
since the stability result in
Lemma~\ref{lem:CLSstab} only applies to functions
in the trial space.
The last inequality (by Lemma~\ref{lem:sampling})
suggests that we can focus on the difference
$u^{CLS}_{X,Y}-I_{Z\cup Y}u^* \in \calu_{Z\cup Y}$,
which has zeros at nodes $Y$. Using the boundary regularity in \eref{eq:bdyreg},  we have
\begin{eqnarray*}
    \|u^{CLS}_{X,Y}-I_{Z\cup Y}u^* \|_{2,\Omega}
    &\leq&
C_{\Omega,\Phi, \call}\big( h_X^{d/2} \|\call u^{CLS}_{X,Y}-\call I_{Z\cup Y}u^*\|_{X}
   { + 
   0
   } \big).
\end{eqnarray*}
Applying Lemma~\ref{lem:CLSconsis} yields
\begin{eqnarray*}
    \| u^{CLS}_{X,Y} - u^* \|_{2,\Omega}
    &\leq&  C_{\Omega,\Phi,\call} ( h_X^{d/2}{n_X^{1/2}} + 1)  h_{Z\cup Y}^{m-2-d/2}\|u^*\|_{m,\Omega},
\end{eqnarray*}
and
\begin{eqnarray*}
    \| u^{CLS}_{X,Y} - u^* \|_{2,\Omega}
    &\leq&  C_{\Omega,\Phi,\call} ( h_X^{d/2}{n_X^{1/2}\rho_X^{d/2}} + 1)  h_{Z\cup Y}^{m-2}\|u^*\|_{m,\Omega},
\end{eqnarray*}
for $m\geq 2+\left\lceil \frac{d+1}2\right\rceil$ and $m>3+d/2$ respectively.
The bracketed factor
is bounded if $X$ is uniformly distributed.
\section{Convergence for WLS}\label{sec:WLS}
Instead of a specific weight, we will consider a class of weighted least-squares formulations by a simple inequality.

\begin{lemma}\label{lem:abineq}
  Let $a,b>0$, $0<\epsilon<1$, and $0\leq\theta\leq2$. Then the following inequalities hold:
  \[
    (\epsilon a+b)^2
    \leq 2\big( \epsilon^\theta\, a^2  +b^2 \big).
  \]
\end{lemma}
\smallskip\pf.  Consider $0\leq \theta/2 \leq 1$. From $(\epsilon a+b)^2 \leq 2(\epsilon^2 a^2+b^2)$ and $\epsilon \leq \epsilon^{\theta/2}$, we have $(\epsilon a+b)^2 \leq 2(\epsilon^\theta a^2+b^2)$.
\qed

\begin{lemma}[$H^2$--Stability]\label{lem:WLSstab}
Suppose the assumptions in Lemma~\ref{lem:CLSstab} hold under the condition \eref{eqZYcond}.
If {$h_X \leq h_Y <1$}, then there exists a constant $C_{\Omega,\Phi, \call}$, depending only on  $\Omega$, $\Phi$, and $\call$ such that
\[
   \| u \|_{2,\Omega}  \leq C_{\Omega,\Phi, \call} 
   h_Y^{d/2-2}\left( \left(\frac{h_Y}{h_X}\right)^{d\theta/2}h_Y^{-2\theta} \right)^{-1/2}
   \left( \|\call u\|_{X}^2 + \left(\frac{h_Y}{h_X}\right)^{d\theta/2}h_Y^{-2\theta} \|u\|_Y^2 \right)^{1/2}
\]
holds for all $0\leq \theta\leq2$ and all $u\in {\calu_{Z\cup Y}}$ for any
finite sets $X\subset\Omega$ and  $Y\subset\Gamma$.
\end{lemma}
\smallskip\pf.
The CLS stability in Lemma~\ref{lem:CLSstab} has to be further modified to suit
the need of WLS. With the denseness requirement
{\eref{eqZYcond}},
let us start with
\begin{equation}\label{eq:stab}
   \| u \|_{2,\Omega} \leq C_{\Omega,\Phi, \call} (h_X^{d/2} \|\call u\|_{X} + h_Y^{d/2-2}\|u\|_Y)
\end{equation}
{for all $u \in \calu_{Z\cup Y}$}.
We want to obtain a stability estimate with discrete sum of squares. Rewrite \eref{eq:stab} as
\[
   \| u \|_{2,\Omega}^2 \leq C_{\Omega,\Phi, \call} h_Y^{d-4} \big( \epsilon \|\call u\|_{X} + \|u\|_Y )^2 \quad\mbox{with } \epsilon = (h_X /h_Y)^{d/2}h_Y^2.
\]
Note that having $\epsilon<1$ is a very mild requirement, for example {$h_X \leq h_Y <1$}, and will not be an obstacle between theories and practice.
By Lemma~\ref{lem:abineq}, we have
\begin{eqnarray*}
   \| u \|_{2,\Omega}
   &\leq& \left( C_{\Omega,\Phi, \call} h_Y^{d-4} \big( \epsilon^\theta \|\call u\|_{X}^2 + \|u\|_Y^2 ) \right)^{1/2}
   \\
   &\leq& C_{\Omega,\Phi, \call} h_Y^{d/2-2} \epsilon^{\theta/2} \big(  \|\call u\|_{X}^2 + \epsilon^{-\theta} \|u\|_Y^2 )^{1/2},
\end{eqnarray*}
for any $0\leq\theta\leq 2$.   {Substituting}
 $\epsilon$ back yields
\[
   \| u \|_{2,\Omega} \leq C_{\Omega,\Phi, \call} h_Y^{d/2-2} ((h_X /h_Y)^{d/2}h_Y^2)^{\theta/2} \big(  \|\call u\|_{X}^2 + ((h_X /h_Y)^{d/2}h_Y^2)^{-\theta} \|u\|_Y^2 )^{1/2},
\]
and we obtain the desired WLS stability after simplification. \qed

\begin{lemma}[Consistency]\label{lem:WLSconsis}
For any $W>0$, define a functional $J_W:H^m(\Omega)\to\R$ by $J_W(u) := \big( \|\call u\|_{X}^2 + W\|u\|_Y^2 \big)^{1/2}$. Suppose the assumptions in Lemma~\ref{lem:CLSconsis} hold. Then, the error estimates in Lemma~\ref{lem:CLSconsis} also hold if the left-handed sides are replaced by $\displaystyle\min_{v\in\calu_{Z\cup Y}}J_W(v-u^*)$ for any $W>0$.
\end{lemma}
\pf. Again, we compare the minimizer $v^*$ with the interpolant $I_{Z\cup Y}u^*$ in $\calu_{Z\cup Y}$:
\begin{eqnarray*}
   J_W^2(v^*-u^*) &\leq& J_W^2( I_{Z\cup Y}u^*-u^*)
    \\
    &=& \| \call   I_{Z\cup Y}u^* - \call u^*  \|_{X}^2 + W\| I_{Z\cup Y}u^* - u^*\|_Y^2,
\end{eqnarray*}
where the last term vanishes {due} to the zeros of $I_{Z\cup Y}u^* - u^*$ at $Y$. 
\qed

With both consistency and stability results, we can now prove the convergence of a class of  WLS solutions defined by \eref{eq:WLSdef}. By similar arguments used in Section~\ref{sec:CLS}, we only need to show that the WLS solution converges to the interpolant $I_{Z\cup Y}u^*$ of the exact solution $u^*$ from the trial space $\calu_{Z\cup Y}$.
For $0\leq\theta\leq2$, consider the functional
\begin{equation}\label{eq:Jtheta}
  J_{W(\theta)}(u):= \big( \|\call u\|_{X}^2 + W(\theta) \|u\|_Y^2 \big)^{1/2} \quad\mbox{with }
  W(\theta):=(h_Y/ h_X)^{d\theta/2}h_Y^{-2\theta}.
\end{equation}
Applying the results of Lemmas \ref{lem:WLSstab} and \ref{lem:WLSconsis}, we
have the WLS solution convergence within the trial space; for simplicity, let
$\tau$ be $d/2$ if $m\geq 2+\left\lceil\frac{d+1}{2}\right\rceil$
and  zero if $m>3+d/2$. Then,
\begin{eqnarray*}
   \| u^{WLS,\theta}_{X,Y,Z\cup Y}- I_{Z\cup Y} u^* \|_{2,\Omega}
   &\leq& C_{\Omega,\Phi, \call}
   h_Y^{d/2-2}W^{-1/2}
   J_{W(\theta)}(u^{WLS,\theta}_{X,Y,Z\cup Y}- I_{Z\cup Y} u^* )
   \\
   &\leq& C_{\Omega,\Phi, \call,\gamma_X} \left(\frac{h_X}{h_Y}\right)^{d\theta/4} h_X^{-d/2}h_Y^{d/2-(2-\theta)}
    h_{Z\cup Y}^{m-2-\tau}\|u^*\|_{m,\Omega}.
\end{eqnarray*}
The last holds because
$h_X^{d/2}n_X^{1/2}\rho_X^{d/2}$ can be bounded by some $C_{\gamma_X}$.
Now we can compare {the}
difference between {the} WLS solution and the exact solution,
\begin{eqnarray*}
   &&\| u^{WLS,\theta}_{X,Y,Z\cup Y}- u^* \|_{2,\Omega}
   \leq \|  I_{Z\cup Y} u^* - u^*\|_{2,\Omega}+ \| u^{WLS,\theta}_{X,Y,Z\cup Y}- I_{Z\cup Y} u^* \|_{2,\Omega}
   \\
   &&\qquad\leq C_{\Omega,\Phi, \call ,\gamma_X}(h_{Z\cup Y}^{m-2}\|u^*\|_{m,\Omega}+h_Y^{-(2-\theta)(4-d)/4} h_X^{d\theta/4-d/2}
    h_{Z\cup Y}^{m-2-\tau}\|u^*\|_{m,\Omega})
    \\
   &&\qquad\leq C_{\Omega,\Phi, \call,\gamma_X }
  \big(1+ h_X^{-d(2-\theta)/4}   h_Y^{(2-\theta)(d-4)/4}   \big)
    h_{Z\cup Y}^{m-2-\tau}\|u^*\|_{m,\Omega},
\end{eqnarray*}
for any $0\leq \theta \leq 2$.
{The constant $1$, coming from $\|  I_{Z\cup Y} u^* - u^*\|_{2,\Omega}$,
is absolutely necessary or else the error bound
will allow arbitrarily fast convergence with respect to $h_X\to0$ for some $\theta$ and $d$.
}
It is obvious that $\theta=2$ maximizes the convergence rate:
\begin{eqnarray*}
   \| u^{WLS,2}_{X,Y,Z\cup Y}- I_{Z\cup Y} u^* \|_{2,\Omega}
&\leq&
   C_{\Omega,\Phi, \call,\gamma_X}
    h_{Z\cup Y}^{m-2-\tau}\|u^*\|_{m,\Omega},
\end{eqnarray*}
where $\tau=d/2$ for $m\geq 2+\left\lceil\frac{d+1}{2}\right\rceil$ or $\tau=0$ for $m>3+d/2$.
The CLS and the optimal WLS($\theta=2$) formulation share convergence estimates of the same form.
They both match the convergence estimate of the interpolant exactly for $m>3+d/2$, that in turn confirms their optimality. To complete proving Theorem~\ref{thm:WLS}, we consider the stability for $\theta=2$ and Lemma~\ref{lem:WLSstab} gives
\begin{eqnarray*}
    \| u \|_{2,\Omega}
    &\leq&  C_{\Omega,\Phi, \call} h_X^{d/2}
    (\|\call u\|_{X}^2 + (h_Y/h_X)^{d}h_Y^{-4}\|u\|_Y^2)^{1/2}
    \\
    &\leq&  C_{\Omega,\Phi, \call} h_X^{d/2}
    (\|\call u\|_{X}^2 + (h_Y/h_X)^{d}h_Y^{-2\theta}\|u\|_Y^2)^{1/2},
\end{eqnarray*}
for any $\theta\geq2$ as long as $h_Y<1$. {We extend} 
the definition of functional $J_{W(\theta)}$ to $\theta\geq2$ by the same definition as in \eref{eq:Jtheta}.
Then, for any $u\in H^m(\Omega)$, we have
\[
     J_{\theta_1}(u) \leq    J_{\theta_2}(u), \quad\mbox{for  $2\leq \theta_1 \leq \theta_2\leq \infty$}.
\]
Since the CLS formulation is equivalent to the WLS with $\theta=\infty$, for any $\theta\geq 2$,  we have
\[
    \min_{v\in\calu_{Z\cup Y}} J_{W(\theta)}(v-u^*) \leq J_{W(\theta)}(u^{CLS}_{X,Y}-u^*) \leq J_{W(\infty)}(u^{CLS}_{X,Y}-u^*),
\]
where the last term is minimal by the definition of CLS solution.
Theorem~\ref{thm:WLS} can now be concluded based on Theorem~\ref{thm:CLS}. %
We have to take \eref{eqZYcond} into account in both theorems.

{\noindent \textbf{Remark:}
Theorem~\ref{thm:CLS} suggests that how fine the boundary collocation $Y$ should be.
In particular, for $0\leq \theta \leq 2$, we need
$  h_Y \leq h_X^{ 1-4/d} $
in order to ensure the optimal convergence rate and the importance of boundary collocation increases with dimensions.
}

\section{Optimal WLS weighting revisited}\label{sec:UZ}

Hu, Chen and \etal \cite{Hu+ChenETAL-Weigradibasicoll:07} showed by scaling analysis that the optimal weighting for overdetermined Kansa methods is $n_Z^2$
{for bounded $\Omega\subset \R^2$.
Since the fill distances $h_X$, $h_Y$, and $h_Z$ are of the same magnitude in \cite{Hu+ChenETAL-Weigradibasicoll:07}, 
this weighting corresponds to $\theta=1$  in our notation. 
To ensure our proven theories are consistent with the previous findings,
we must extend our theories to a smaller  trial space  $\calu_Z$  used in
\cite{Hu+ChenETAL-Weigradibasicoll:07}.
In the rest of this section, we will focus on the WLS
convergence in this smaller trial space and prove Theorem~\ref{cor:WLS}.}



To begin, let us return to the proof for WLS consistency (Lemma \ref{lem:WLSconsis}) but restrict the  approximation in the smaller trial space $\calu_Z$, within which the stability result in Lemma \ref{lem:WLSstab} remains valid.
However,
we can only compare the minimizer $v^*\in \calu_Z$ with the interpolant $I_Z u^*\in \calu_Z$ to the exact solution  $u^*\in H^m(\Omega)$:
\begin{eqnarray*}
   \min_{v\in\calu_Z} J_W^2(v-u^*)
   &\leq& J_W^2(I_Z u^*-u^*)
   \\
   &=& \|\call I_Z  u^*-\call  u^*\|_{X}^2 + W\|I_Z u^*-u^*\|_Y^2.
\end{eqnarray*}
The PDE residual on $X$ is exactly the same as that in the previous section.
Without $Y$ in the trial centers to annihilate the boundary collocation, we simplify need to identify the extra terms associated to boundary error on $Y$.
Following the ideas in the proof of Lemma~\ref{lem:CLSconsis},
for $m\geq 2+\left\lceil\frac{d+1}2\right\rceil$, we can bound the boundary term by
\begin{eqnarray*}
\|I_Z u^*-u^*\|_Y^2&\leq& n_Y^{1/2} \|I_Z u^*-u^*\|_{L^\infty(\Gamma)}
\\
&\leq& n_Y^{1/2} \|I_Z u^*-u^*\|_{L^\infty(\Omega)}
\\
&\leq& C_{\Omega,\Phi, \call} n_Y^{1/2} h_{Z}^{m-d/2} \|u^*\|_{m,\Omega}. 
\end{eqnarray*}
Hence, the error estimate for WLS on $\calu_Z$ contains an extra term
\begin{eqnarray}
    &&
    h_Y^{d/2-2}W^{-1/2}\left( W^{1/2}C_{\Omega,\Phi, \call} n_Y^{1/2} h_{Z}^{m-d/2} \|u^*\|_{m,\Omega} \right)
    \nonumber\\&&\qquad\leq
    C_{\Omega,\Phi, \call,\gamma_Y} h_Y^{d/2-2}h_Y^{-(d-1)/2}h_{Z}^{m-d/2} \|u^*\|_{m,\Omega}
    \nonumber\\&&\qquad\leq
    C_{\Omega,\Phi, \call,\gamma_Y} (h_Y^{-3/2}h_Z^2 ) h_{Z}^{m-d/2-2} \|u^*\|_{m,\Omega}.
    \label{smallCLS1}
\end{eqnarray}

For the other case when $m>3+d/2$,
we want to bound the $\ell_2(Y)$ norm on the boundary by some $\ell_2(\widetilde{Z}\cup Y)$ norm in the domain  (like the trace theorem does). For any subset $\widetilde{Z}\subseteq Z$,  we have
\begin{eqnarray*}
    \|I_Z u^*-u^*\|_Y
    &{=}& \|I_Z u^*-u^*\|_{\widetilde{Z}\cup Y}
    \\
    &\leq& C_{\Omega,\Phi}
    n_{\widetilde{Z}\cup Y}^{1/2} \rho_{\widetilde{Z}\cup Y}^{d/2}
    h_Z^{m}\|u^*\|_{m,\Omega}.
\end{eqnarray*}
%
We want to select $\widetilde{Z}$ so that $\rho_{\widetilde{Z} \cup Y}$ can be bounded by the denseness measures of $Y$ and $Z$.
We already assumed $Z$ is quasi-uniform  in Assumption \ref{ass:RSCP}.
Let us further assume that $Y$ is also quasi-uniform
on $\Gamma$ and  satisfies \eref{eq:qu} with some constant
$\gamma_Y>1$. Moreover, the set $Y$ is sufficiently dense with respect to $Z$
and $\Omega$ so that $h_Y\leq h_Z$ and $q_{Y,\Gamma} < q_{Y,\Omega}$
(see \cite[Thm.6]{Fuselier+Wright-ScatDataInteEmbe:12}).


Consider the following subset that excludes all points in $Z$ that are within distance $h_Z$ to the boundary:
\[
    \widetilde{Z} := \Big\{ z\in Z \cap
    \{ \zeta\in\Omega \,:\; \mbox{dist}(\zeta-\Gamma)> h_Z\}
    \Big\} \subseteq Z.
\]
Then,  we have
\begin{eqnarray*}
  \min(q_Z,q_Y)
  \leq  q_{\widetilde{Z}\cup Y}
  \leq h_{\widetilde{Z}\cup Y}
  &\leq&
  \Big(\sup_{\zeta\in\Omega_{h_Z}} + \sup_{\zeta\in{\Omega\setminus\Omega_{h_Z}}}\Big)
            \min_{z \in {\widetilde{Z}\cup Y} }\|z-\zeta\|_{\ell_2(\R^d)}
  \\
  &\leq& h_Z +(h_Z + h_{Y}).
\end{eqnarray*}
It is now clear that the set ${\widetilde{Z}\cup Y}$ is also quasi-uniform with respect to some parameter $\gamma_{\widetilde{Z}\cup Y}$ that depends on $\gamma_Z$ and $\gamma_Y$.
Hence, we can bound $\rho_{\widetilde{Z}\cup Y}$ by some generic constant  $C_{\gamma_Y,\gamma_Z}$.
To control the term $n_{\widetilde{Z}\cup Y}$, consider
\begin{eqnarray*}
  n_{\widetilde{Z}\cup Y} &\leq&  n_Z + n_Y \\
  &\leq & C_{\Omega,\gamma_Z}   h_Z^{-d} + C_{\Omega,\gamma_Y}  h_Y^{-(d-1)}
  \\  &\leq &
  C_{\Omega,\gamma_Y,\gamma_Z}\big( h_Y^{-d} + h_Y^{-d+1}\big).
\end{eqnarray*}
Since we assumed $h_Y<1$, we have $n_{\widetilde{Z}\cup Y} \leq C_{\Omega,\gamma_Y,\gamma_Z} h_Y^{-d}$.
Together, we have
\[
    \|I_Z u^*-u^*\|_Y \leq
     { C_{\Omega,\Phi,\gamma_Y,\gamma_Z} h_Y^{-d/2} }
    h_Z^{m}\|u^*\|_{m,\Omega},
\]
and the extra term associated with boundary error on $Y$ is
\begin{eqnarray}
    &&
    h_Y^{d/2-2}W^{-1/2}\left( W^{1/2}{ C_{\Omega,\Phi,\gamma_Y,\gamma_Z} h_Y^{-d/2} }
    h_Z^{m}\|u^*\|_{m,\Omega} \right)
    \nonumber\\&&\qquad\leq
    C_{\Omega,\Phi, \call,\gamma_Y,\gamma_Z} (h_Y^{-2}h_Z^{2})h_{Z}^{m-2} \|u^*\|_{m,\Omega}.
    \label{smallCLS2}
\end{eqnarray}
Adding the corresponding boundary errors in  \eref{smallCLS1} and \eref{smallCLS2}, for $m\geq 2+\left\lceil\frac{d+1}2\right\rceil$ and $m>3+d/2$ respectively, into the estimates in Theorem~\ref{thm:WLS} completes the proof of Theorem~\ref{cor:WLS}.
And, these boundary errors do not affect the convergence rates as long as  $h_Y \lesssim h_Z$.
Moreover, since both \eref{smallCLS1} and \eref{smallCLS2} are
independent {of} $\theta$, our error estimate allows the least-squares weighting
suggested in \cite{Hu+ChenETAL-Weigradibasicoll:07} to be optimal, but we have to take \eref{eqZYcond2} and other requirements on the trial centers and collocation points into account.

\section{Numerical demonstrations}\label{sec:num}
We test the proposed formulations in $\Omega=[-1,1]^2$,
Discretization is done  by using  regular  $Z$ with $n_Z = 11^2,16^2,\ldots,36^2$. 
collocation points $X$ are either regular or scattered.
For the regular cases,
collocation points $X$ (strictly in the interior) and $Y$ are constructed similarly with $h_X =\delta_i h_Z$ and $h_Y = \delta_b h_Z $ with $\delta_i=1,1/2,1/3$ and $\delta_b=1,1/2$ such that $Z\subseteq X$ and $(Z\cap \Gamma)\subseteq Y$ respectively.
All reported errors, either in $L^2$ or $H^2$, are absolute error and are approximated by
{using}
a fixed set  of $100^2$ regular points, which is denser than the collocation $X\cup Y$ sets in all tests.

In matrix form, collocation conditions for the PDE and boundary condition can be written as
\[
    K_{\call,X} \lambda = f_{|X} \mbox{\quad and \quad} K_{\calb,Y} \lambda = g_{|Y},
\]
respectively, with entries $[K_{\call,X}]_{ij} = \call\Phi( x_i-z_j)$ and $[K_{\calb,Y}]_{ij} = \Phi( y_i-z_j)$ for $x_i\in X$ and $y_i\in Y$.
Both resultant matrices have $n_Z+n_Y$ (and  $n_Z$) columns for trial space $\calu_{Z\cup Y}$ (and $\calu_Z$) corresponding to each $z_j$ from the trial space.
In the CLS approach \eref{eq:CLSdef}, the constraints at $Y$ are enforced using the null space matrix of the boundary collocation matrix, denoted by $\caln_{\calb,Y}:=\mbox{null}( K_{\calb,Y} )$, as in \cite{Ling+Hon-ImprnumesolvKans:05}, so that the unknown coefficient is expressed in the form
\[
    \lambda = \caln_{\calb,Y}\gamma + K_{\calb,Y}^\dag g_{|Y},
\]
for some new unknown $\gamma$, which can be found by solving
\[
 K_{\call,X} \caln_{\calb,Y}\gamma  = f_{|X} - K_{\call,X} K_{\calb,Y}^\dag g_{|Y}.
\]

In all WLS($\theta$) formulations, with LS weighting specified by $W(\theta)$  in \eref{eq:WLSdef}, the unknown coefficient $\lambda$ is obtained by solving the following overdetermined  system
\[
    \left[
      \begin{array}{c}
        K_{\call,X} \\
       W(\theta) K_{\calb,Y} \\
      \end{array}
    \right]
    \lambda
    =
    \left[
      \begin{array}{c}
        f_{|X}  \\
       W(\theta) g_{|Y} \\
      \end{array}
    \right]
\]
with the Matlab function \verb"mldivide" in the least-squares sense.
For all computations in this section, we did not employ any technique to deal with the problem of ill-conditioning unless specified otherwise, i.e., when we employ the stable RBF-QR decomposition for the Gaussian basis. To deal with the numerical instability, readers are referred to our trial subspace selection techniques \cite{Ling+Schaback-imprsubsselealgo:09}.

\begin{figure}
  \centering
\begin{tabular}{ccc}
  \begin{overpic}[width=0.32\textwidth,trim= 20 20 25 25, clip=true,tics=10]{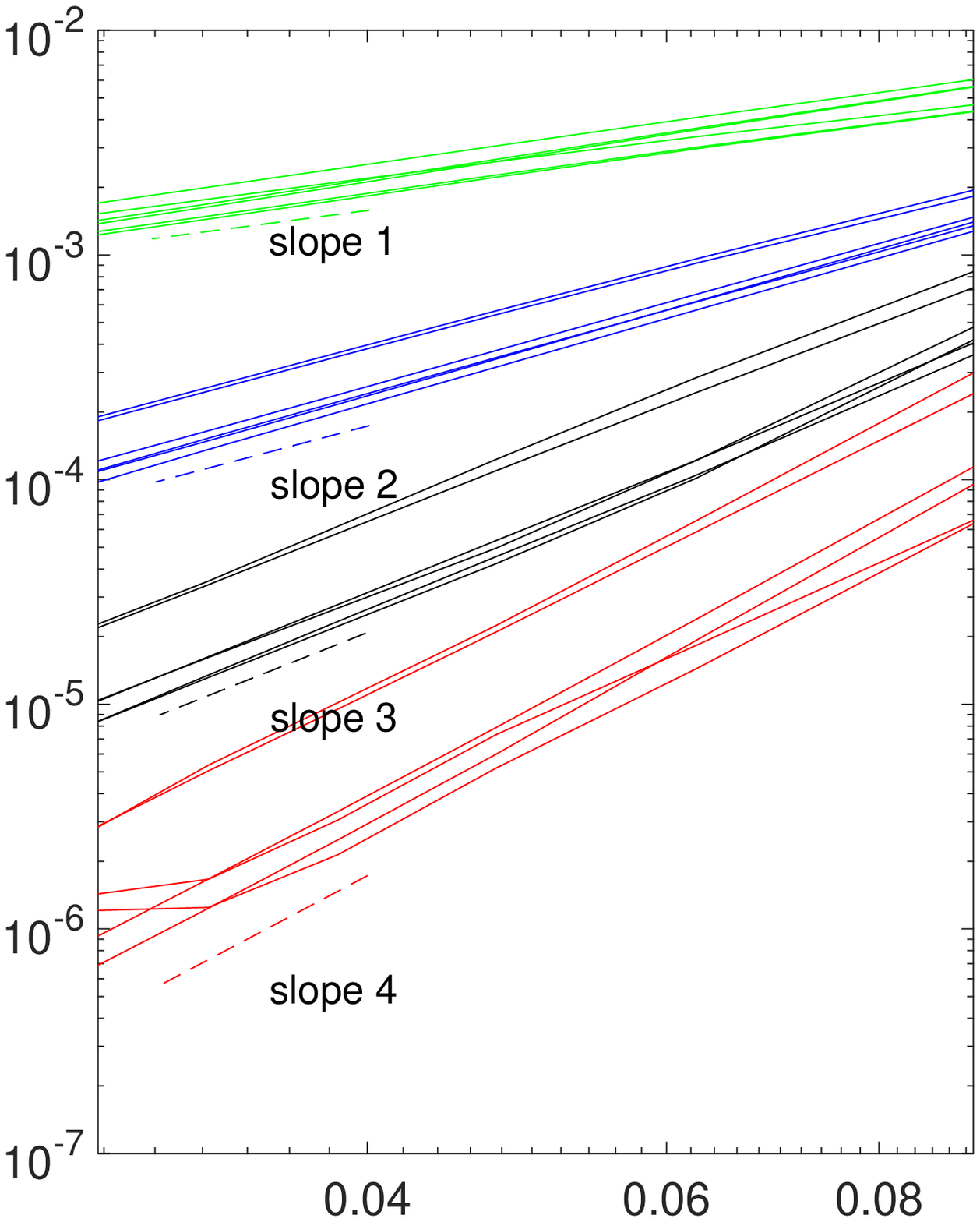}
  \put (8,100) {\scriptsize $u^*=\sin(\pi x/2)\cos(\pi y/2)$ }
  \put (35,12) {CLS in $\calu_{Z\cup Y}$}
  \end{overpic}
  &
  \begin{overpic}[width=0.32\textwidth,trim= 20 20 25 25, clip=true,tics=10]{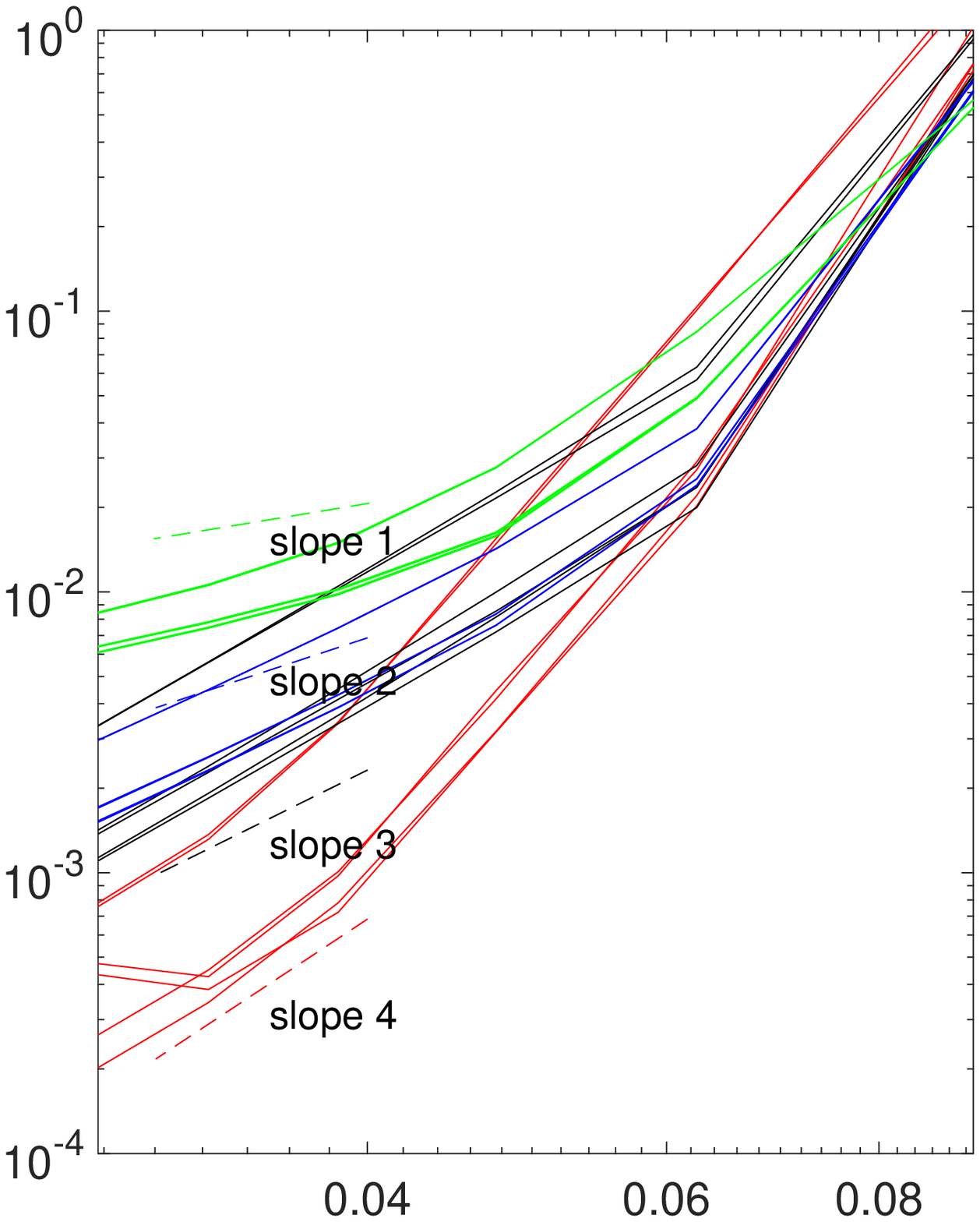}
  \put (8,100) {\scriptsize $u^*=\peaks(3x,3y)$ }
  \end{overpic}
  &
  \begin{overpic}[width=0.32\textwidth,trim= 20 20 25 25, clip=true,tics=10]{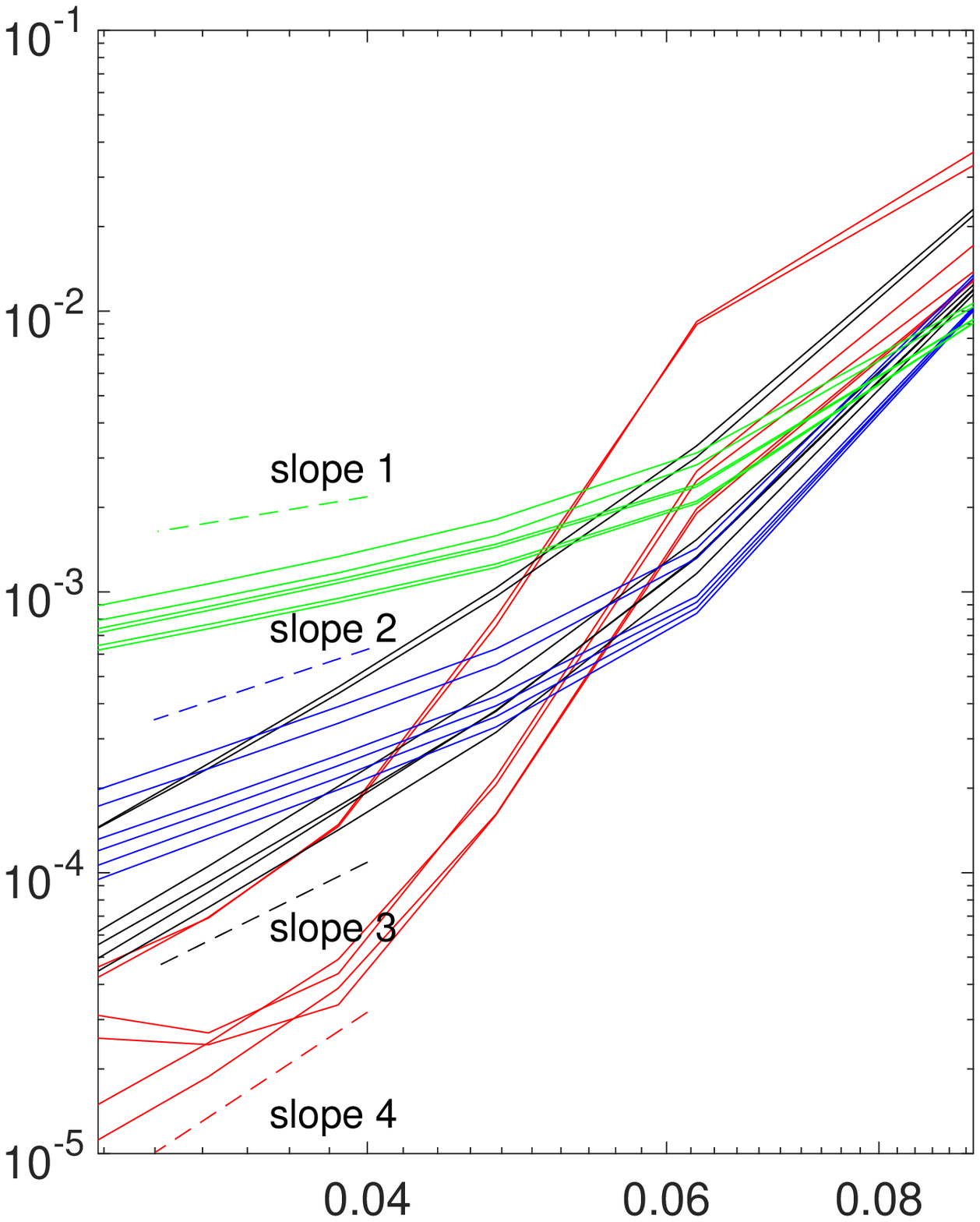} 
  \put (8,100) {\scriptsize $u^*=\franke(2x-1,2y-1)$ }
  \end{overpic}
\end{tabular}
  \caption{Example~\ref{examp1}: $H^2(\Omega)$ error profiles for casting the CLS formulation in $\calu_{Z\cup Y}$ with Whittle-Mat\'{e}rn-Sobolev kernels of order $m=3,\ldots,6$ (green, blue, black, and red) to solve $\Delta u = f$ with different exact solution $u^*$. }\label{fig:CLSUZYconv}
%
%
%
\bigskip
\begin{tabular}{ccc}
  \begin{overpic}[width=0.32\textwidth,trim= 20 20 25 25, clip=true,tics=10]{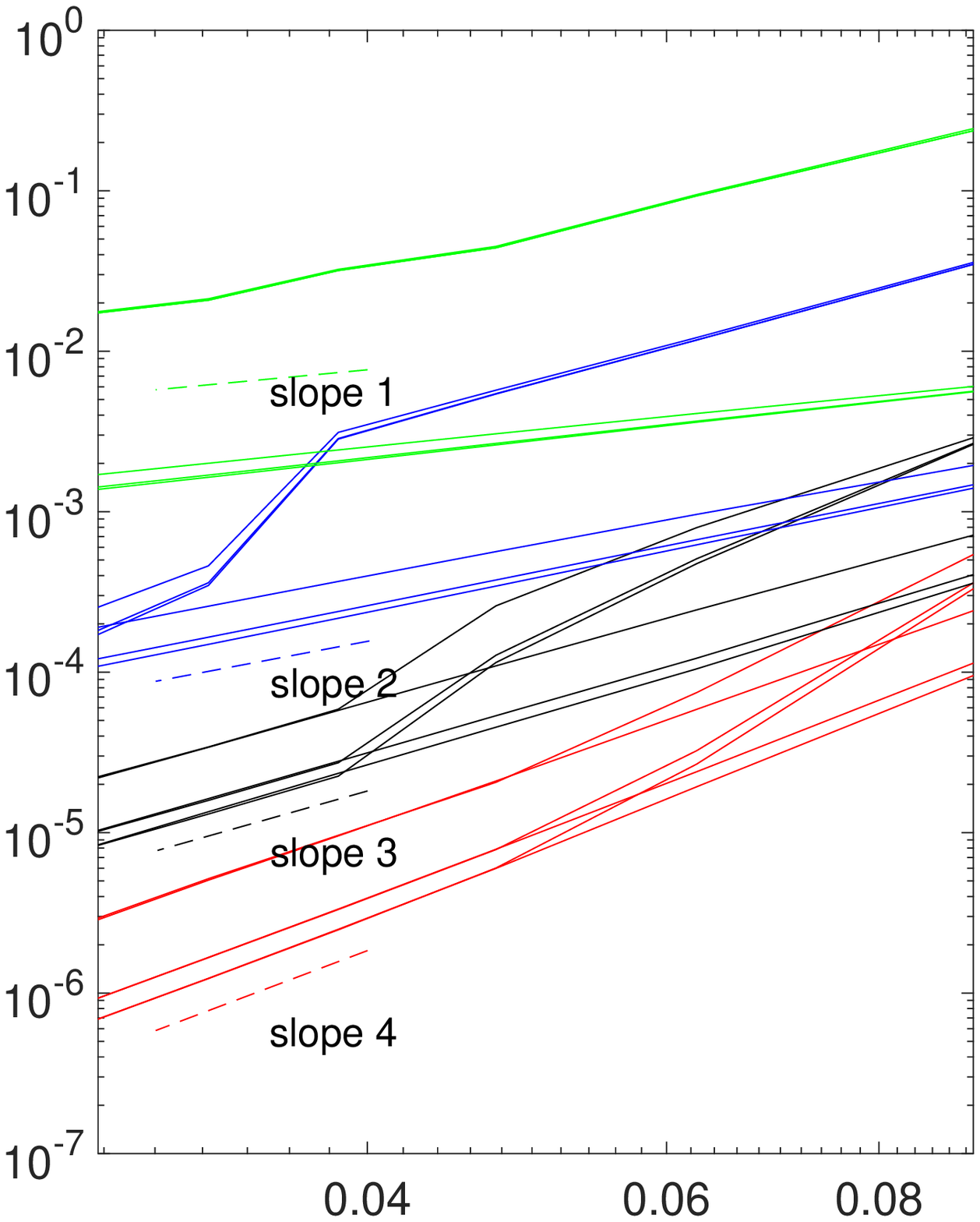}
  \put (42,12) {CLS in $\calu_{Z}$}
  \end{overpic}
  &
  \begin{overpic}[width=0.32\textwidth,trim= 20 20 25 25, clip=true,tics=10]{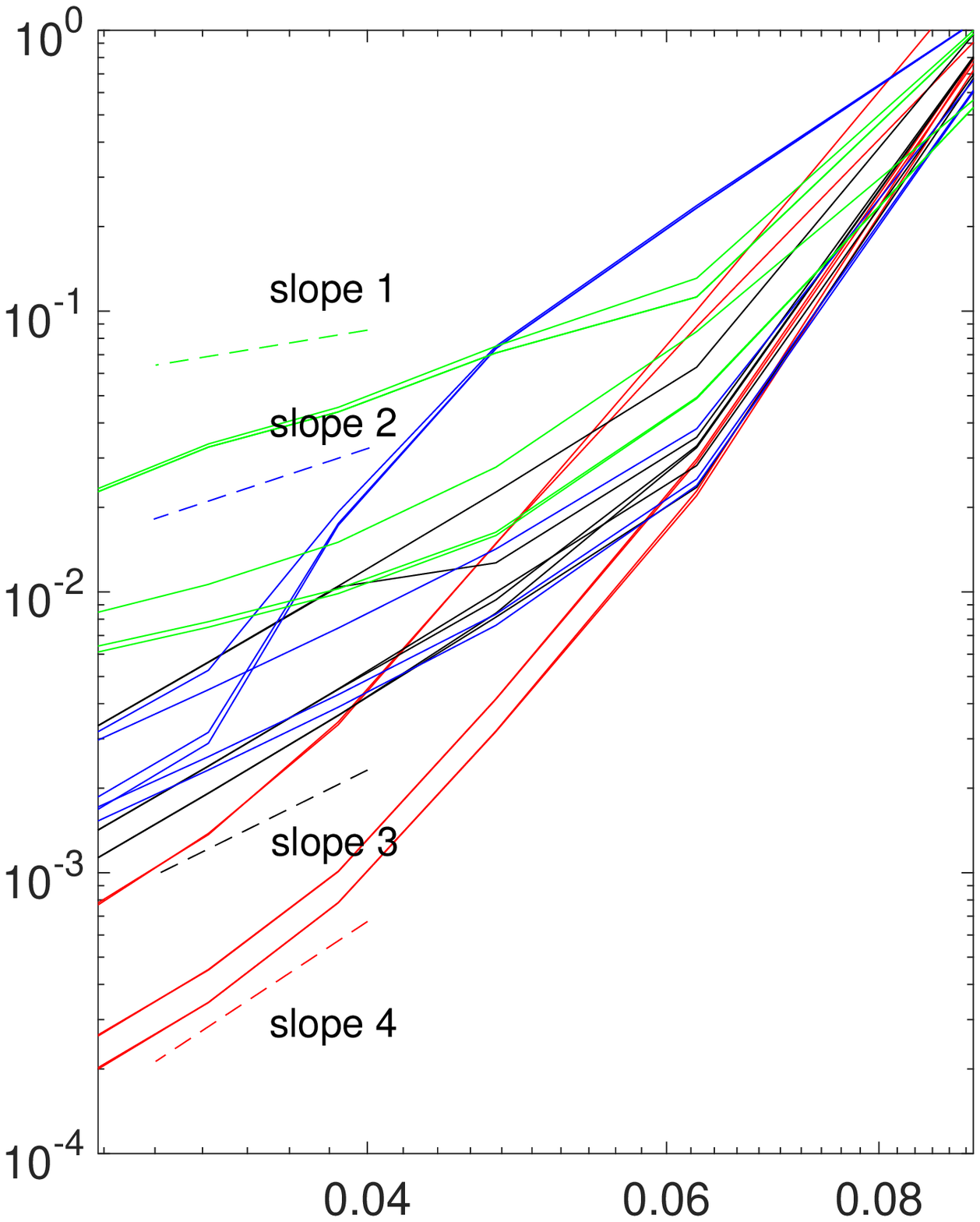}
  \end{overpic}
  &
  \begin{overpic}[width=0.32\textwidth,trim= 20 20 25 25, clip=true,tics=10]{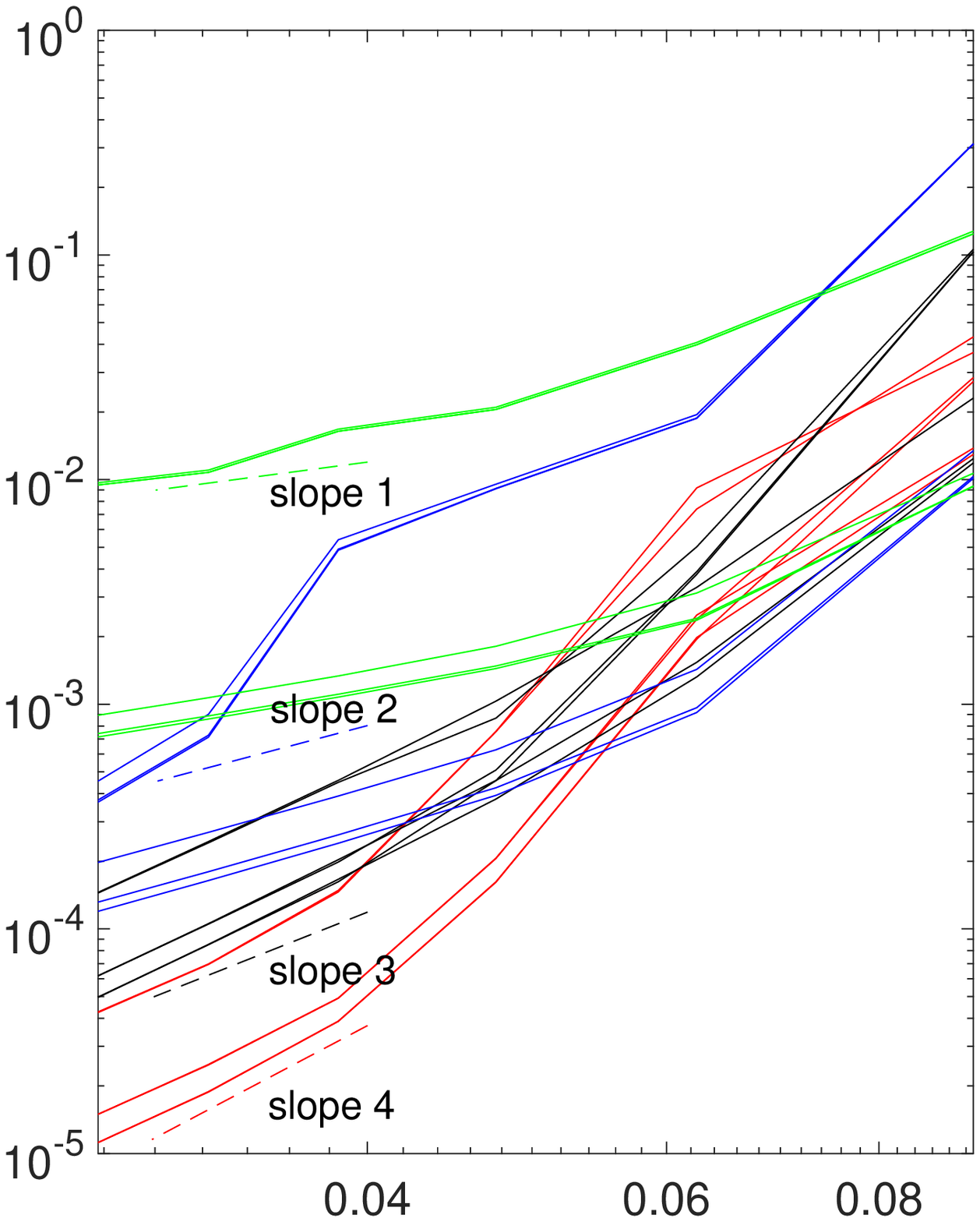} 
  \end{overpic}
\end{tabular}
  \caption{Example~\ref{examp2}: $H^2(\Omega)$ error profiles for casting the CLS formulation in $\calu_Z$ to the same settings as in  Figure~\ref{fig:CLSUZYconv}.}\label{fig:CLSUZconv}
\end{figure}

\smallskip
\begin{example}\label{examp1}
  \textbf{How dense is dense?}
\end{example}
First, we consider a Poisson problem with 
Dirichlet boundary value generated from three different exact solutions $u^*= \sin(\pi x/2)\cos(\pi y/2)$, $\peaks(3x,3y)$, and $\franke(2x-1,2y-1)$ by the corresponding functions in \Matlab.
%
We cast the CLS formulation \eref{eq:CLSdef} using
{unscaled Whittle-Mat\'{e}rn-Sobolev}
kernels that reproduce $H^m(\Omega)$ with $m=3,\ldots,6$.
Note that our proven $H^2$--convergence theories require {$m\geq 2+ \lceil1.5\rceil=4 $ for $\Omega\subset\R^2$.}
To see the effect of ``over-testing'', all sets in this example are regular and we tested $h_X=\{h_Z,h_Z/2,h_Z/3\}$ and $h_Y=\{h_Z,h_Z/2\}$.
Figure~\ref{fig:CLSUZYconv} compactly shows all convergence
profiles in $H^2(\Omega)$ with respect to $h_Z$ (instead of $h_{Z\cup Y}$ for  ease of comparison to the results in the next example) at a glance.

To begin, let us focus on the $H^2(\Omega)$ errors for $u^*= \sin(\pi x/2)\cos(\pi y/2)$ in
Figure~\ref{fig:CLSUZYconv}. 
Generally speaking,  all collocation
{settings} demonstrate an $m-2$ convergence {rate} for all tested smoothness $m$;
this also includes the original Kansa formulation with $Z=X\cup Y$.
It is obvious that the error profiles for each
tested $m$ {are} split into two groups.
The least accurate groups (i.e., the group above)
correspond to $h_X=h_Z$. Without over-testing the PDE, this setting would
probably fail the denseness requirement \eref{eqZYcond}
but yet allow convergence at the optimal rate.
All errors reduce at a rather constant rate,
except that we can see two unstable profiles
in the cases of $m=6$.  These numerical instabilities
correspond to the two cases with large numbers of boundary collocations; $(h_X,h_Y)=(h_Z/2,h_Z/2)$ and $(h_X,h_Y)=(h_Z/3,h_Z/2)$.
In comparison, the other two tested solutions $u^*=\peaks(3x,3y)$, and $\franke(2x-1,2y-1)$ are more oscillatory. We can see the CLS convergence rates slow down and approach the optimal $m-2$ order.
We can clearly see that higher smoothness typically suffers
{more from} the effects of ill-conditioning in all tested $u^*$.  Error
reduction reaches a valley
{as $Z$ is being refined and then increases again.}

We omit the $L^2(\Omega)$ error profiles,
which  show exactly two extra orders as one would expect  and achieve an
$m$-order convergence before numerical instability kicks in.

\smallskip
\begin{example}\label{examp2}
  \textbf{CLS convergence in trial space $\calu_Z$.}
\end{example}
Putting the theoretical requirement aside, we are interested in the numerical performance of casting the CLS in the smaller and more practical trial spaces $\calu_{Z}$. 
Elementary linear algebra says that if $n_Z < n_Y$, then we may not be able to find nontrivial functions from $\calu_Z$ with zeros at $Y$.
However, one can observe numerically that the CLS formulation hardly runs into trouble when it is cast in this smaller trial space.
Numerically, as $h_Y\to0$, the rank of the boundary collocation matrix is bounded; for example, for $n_Z=21^2$ with finer and finer $Y$, we can see from Table~\ref{tab:rank} that the rank of the boundary matrix is numerically bounded.

\begin{table}
\centering
\begin{tabular}{|c|c||c|c|}
  \hline
  \multicolumn{2}{|c||}{$\calu_{Z\cup Y}$} & \multicolumn{2}{|c|}{$\calu_{Z}$} \\ \hline
   $~~~~n_Y~~~~$ & rank($K_{\calb,Y}$) & $~~~~n_Y~~~~$  & rank($K_{\calb,Y}$) \\    \hline
    80 & 80  & 80 & 80 \\
    160& 104 & 160& 96 \\
    244& 108 & 244& 96 \\
    328& 100 & 328& 96 \\
  \hline
\end{tabular}
\caption{Numerical ranks of boundary collocation matrices resulting from different $n_Y$.}\label{tab:rank}
\end{table}

Figure~\ref{fig:CLSUZconv} shows the $H^2(\Omega)$ error profiles for the CLS performance in $\calu_Z$ with all other settings identical to those in Example~\ref{examp1}. Comparing the CLS convergence rates in the two trial spaces, we observe that optimal convergence is also possible in $\calu_Z$ but only for small enough $h_Z$. The CLS accuracy can ``catch up'' when the numerical rank of $K_{\calb,Y}$ is relatively insignificant compared to $n_Z$. Therefore, the larger the rank($K_{\calb,Y}$) the longer CLS takes to achieve optimal convergence.
By using a smaller trial space, we not only gain computational efficiency but suffer less ill-conditioning. In all tested cases and parameters, we see no accuracy drop due to ill-conditioning in Figure~\ref{fig:CLSUZconv}.
In our next demonstration, we will see that using $\calu_Z$ as trial space only makes the CLS convergence lag behind, but not to any other WLS formulations.

\begin{figure}
  \centering
\begin{tabular}{cc}
  \begin{overpic}[width=0.4\textwidth,tics=10]{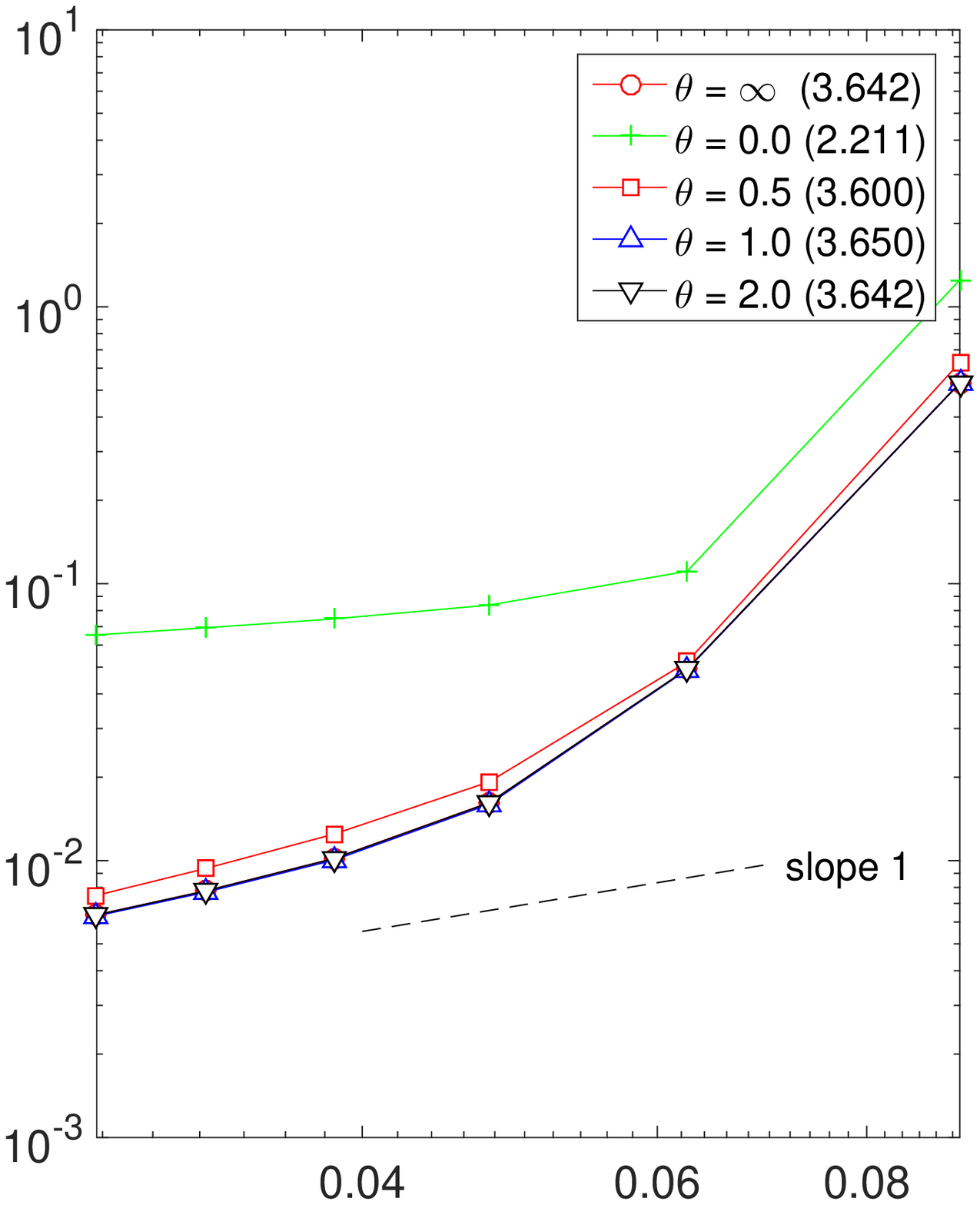} 
  \put (15,85) {$\calu_{Z\cup Y}$ }
  \put (15,95) {$m=3$}
  \end{overpic}
  &
  \begin{overpic}[width=0.4\textwidth,tics=10]{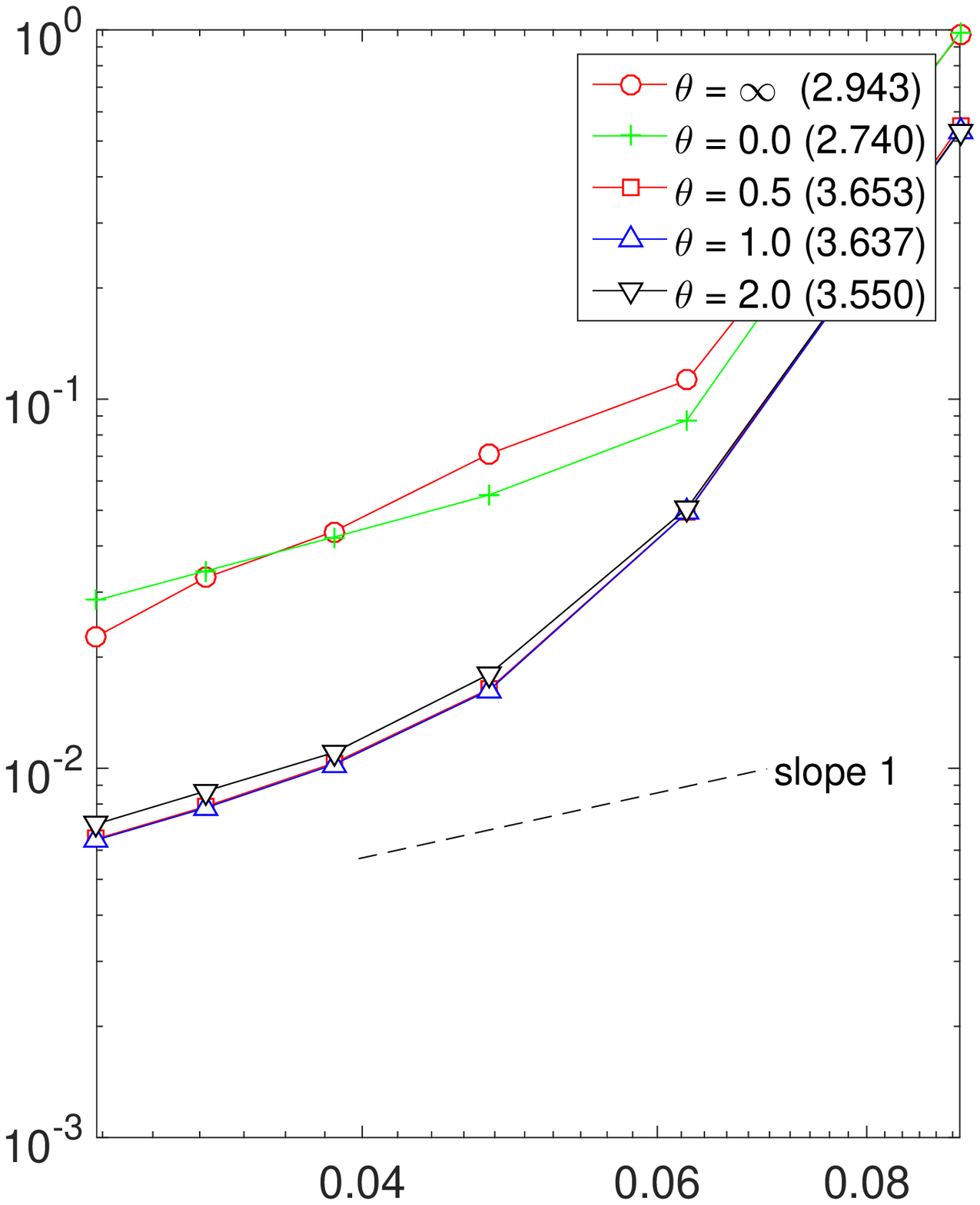} 
  \put (15,85) {$\calu_{Z}$ }
  \end{overpic}
\\
  \begin{overpic}[width=0.4\textwidth,tics=10]{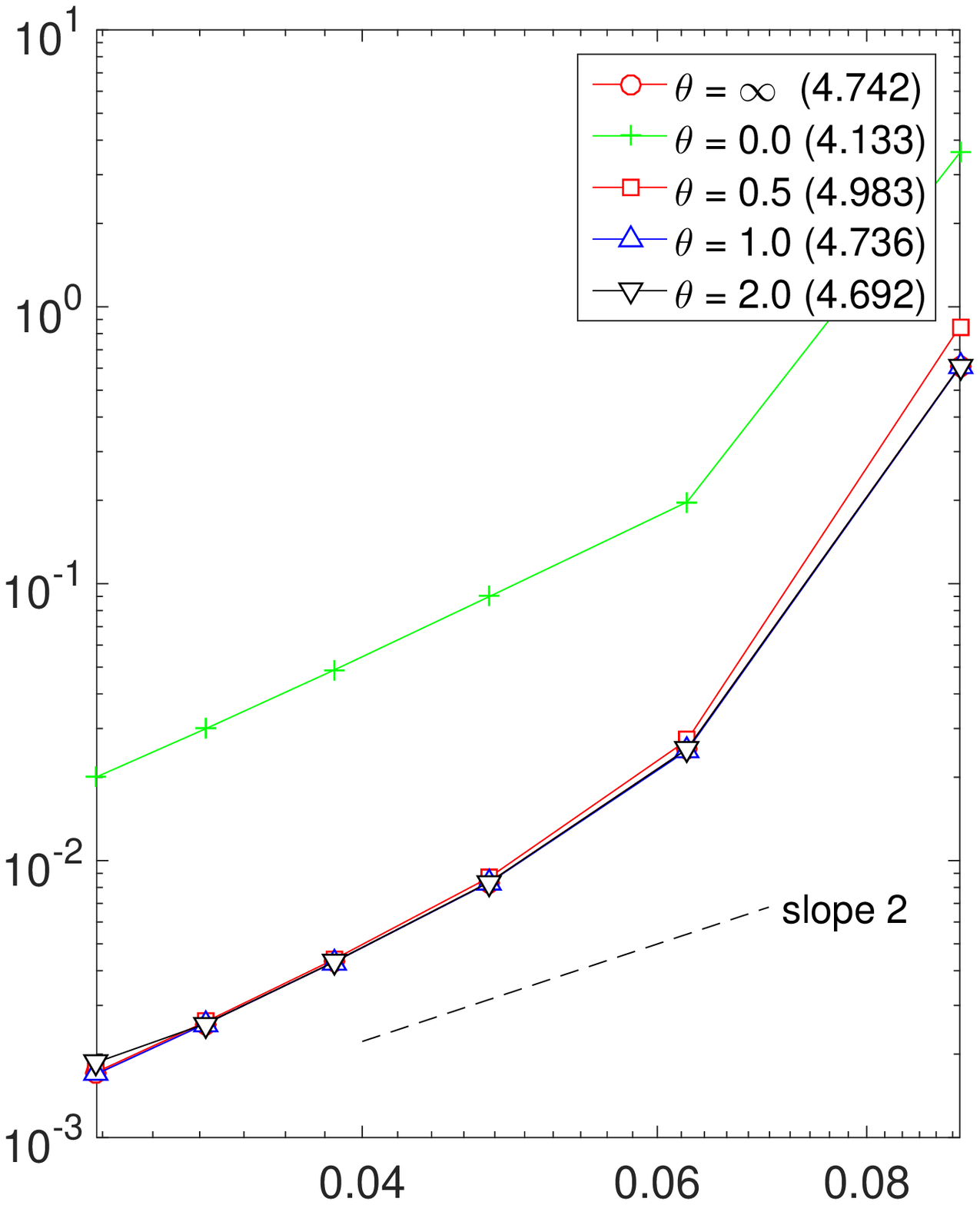} 
  \put (15,85) {$\calu_{Z\cup Y}$ }
  \put (15,95) {$m=4$}
  \end{overpic}
  &
  \begin{overpic}[width=0.4\textwidth,tics=10]{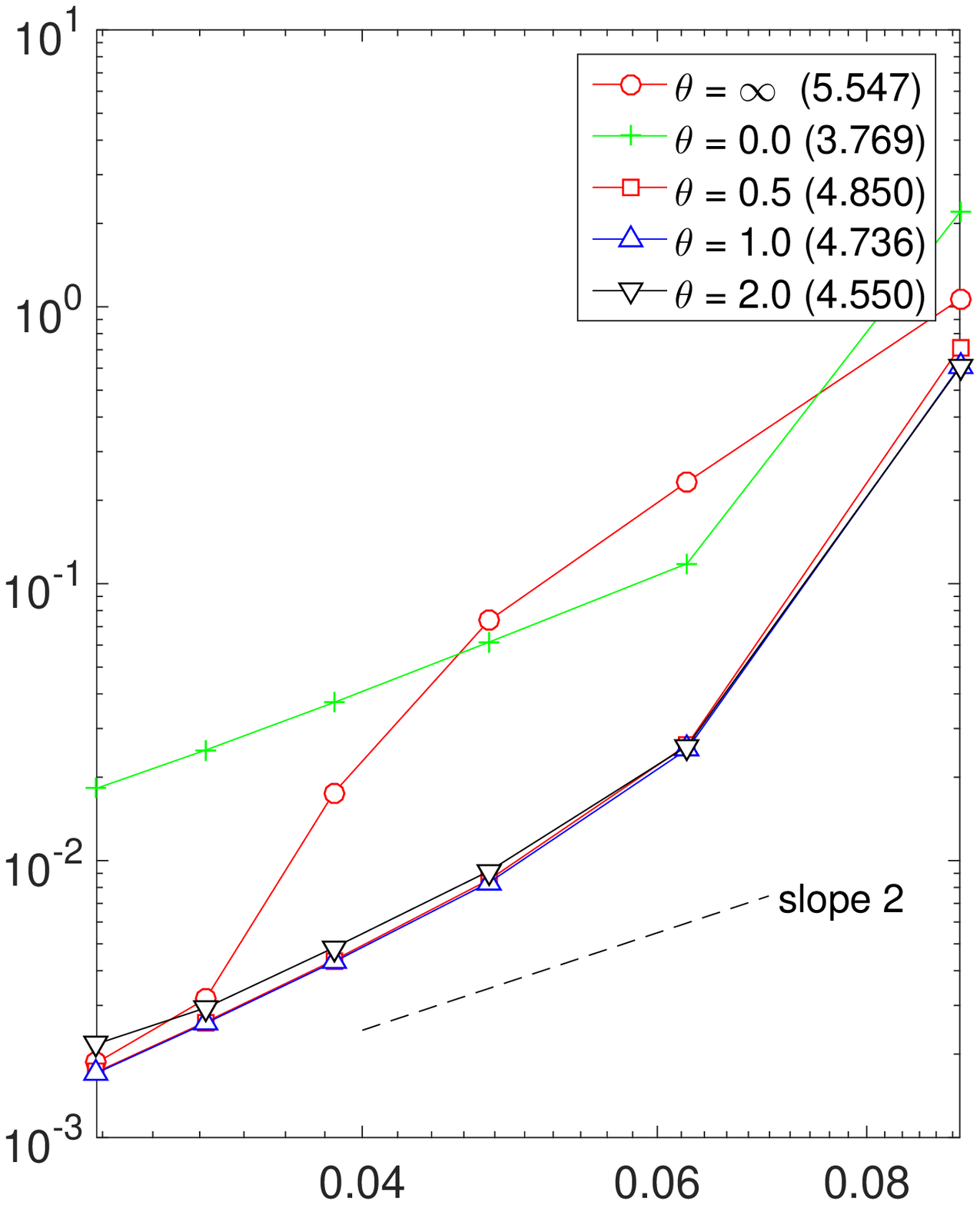} 
  \put (15,85) {$\calu_{Z}$ }
  \end{overpic}
\\
  \begin{overpic}[width=0.4\textwidth,tics=10]{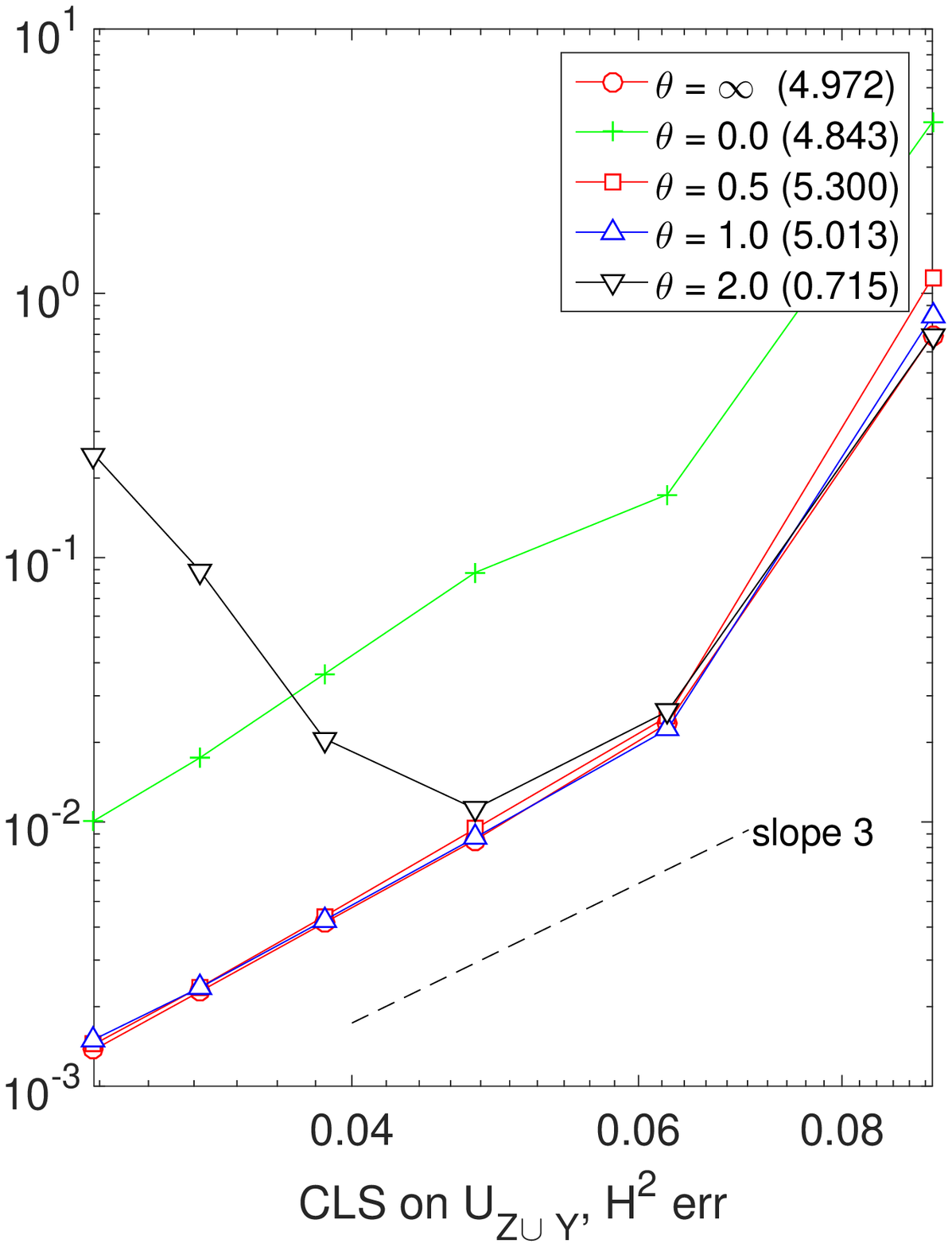} 
  \put (15,85) {$\calu_{Z\cup Y}$ }
  \put (15,95) {$m=5$}
  \end{overpic}
  &
  \begin{overpic}[width=0.4\textwidth,tics=10]{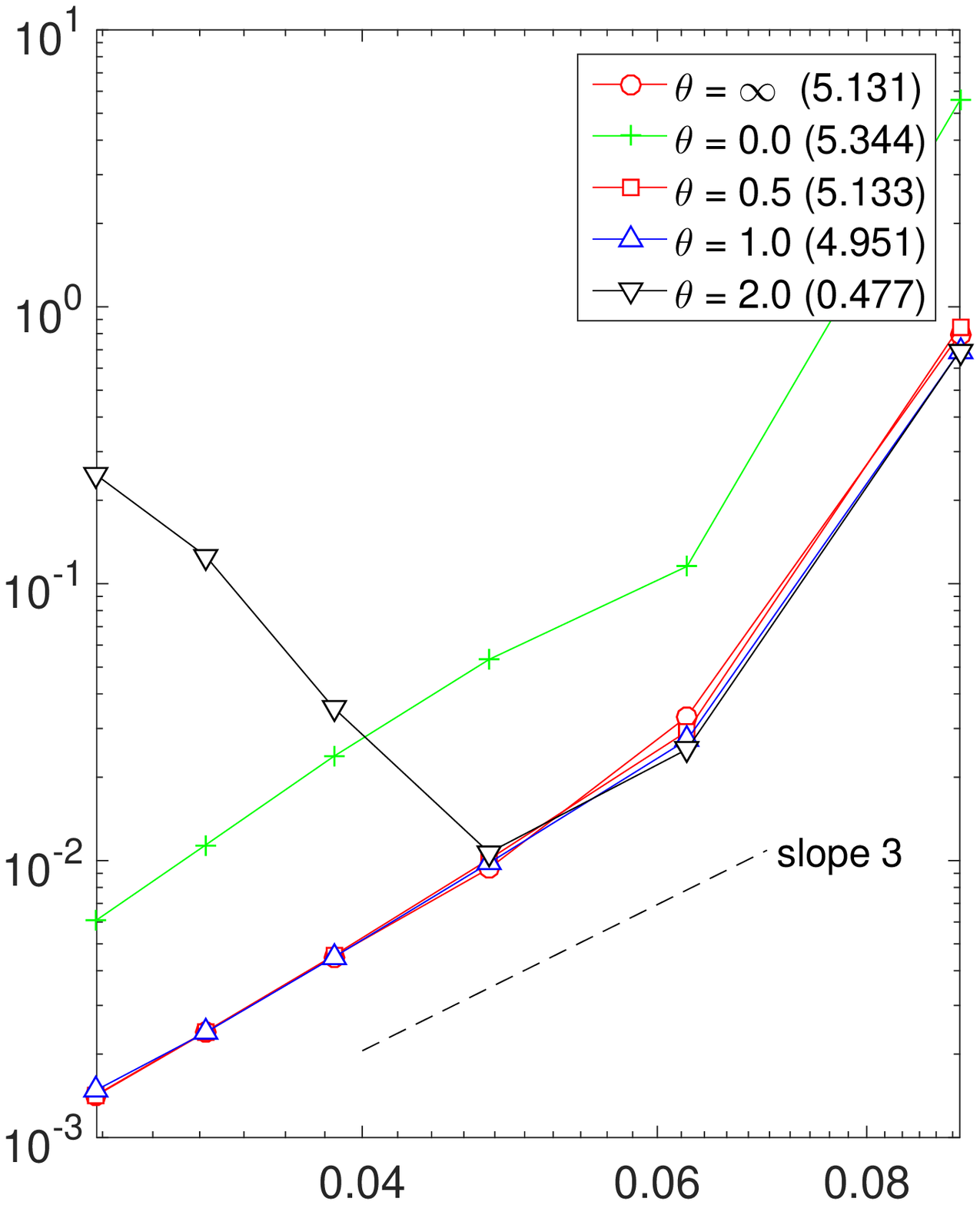} 
  \put (15,85) {$\calu_{Z}$ }
  \end{overpic}
\end{tabular}
  \caption{Example~\ref{examp3}: $H^2(\Omega)$ error profiles for casting the WLS$(\theta)$ formulation in $\calu_{Z\cup Y}$ and $\calu_Z$ with Whittle-Mat\'{e}rn-Sobolev kernels of order $m=3,4,5$ to solve $\Delta u = f$ with exact solution $u^*=\peaks(3x,3y)$.}\label{fig:3apeaks}
\end{figure}

\smallskip
\begin{example}\label{examp3}
  \textbf{Numerically optimal weight for WLS.}
\end{example}
Now we consider the WLS($\theta$) formulation in \eref{eq:WLSdef} with $\theta\in\{\infty, 0, 0.5,1,2\}$.  We begin with the same set up as in Example~\ref{examp1} and set $h_X = h_Y = h_Z/2$ to solve the Poisson problem in $\Omega=[-1,1]^2$.
The WLS weighting in this test are $W(\theta) =1, \,n_Z, \,n_Z^2$ and $n_Z^4$, and WLS($\infty$) is equivalent to the CLS.
Figure~\ref{fig:3apeaks} show the $H^2(\Omega)$ error resulting from various $WLS(\theta)$ formulations associated with  $u^*=\sin(\pi x/2,\pi y/2)$ and  $u^*=\peaks(3x,3y)$ respectively.
The estimated convergence rates shown in the legends are obtained from least-squares fitting to all data; if the convergence profile is not ``straight'' enough,   the corresponding estimate is not trustworthy.

\begin{figure}
  \centering
  \includegraphics[width=0.6\hsize]{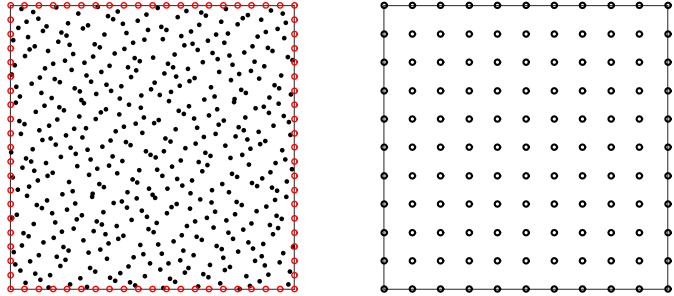}\\
  \caption{Schematic point sets, collocations $X$ and $Y$ (left), and trial $Z$ (right), used to solve various PDEs.}\label{fig:PDEpt.eps}
\end{figure}

From these figures, we immediately see that there is no benefit at all (in terms
of both efficiency and accuracy) to go for the unweighted WLS(0)
formulation. Unlike the CLS in $\calu_Z$, 
{all tested WLS($\theta$) with $\theta>0$ 
do not have a lag} in convergence rate but may 
suffer ill-conditioning for large $\theta$.
By comparing all tested cases, we see that $\theta=1$ {allows}
good accuracy and optimal convergence rate in both trial spaces.

\begin{figure}
  \centering
\begin{tabular}{cc} 
  \begin{overpic}[width=0.4\textwidth,tics=10]{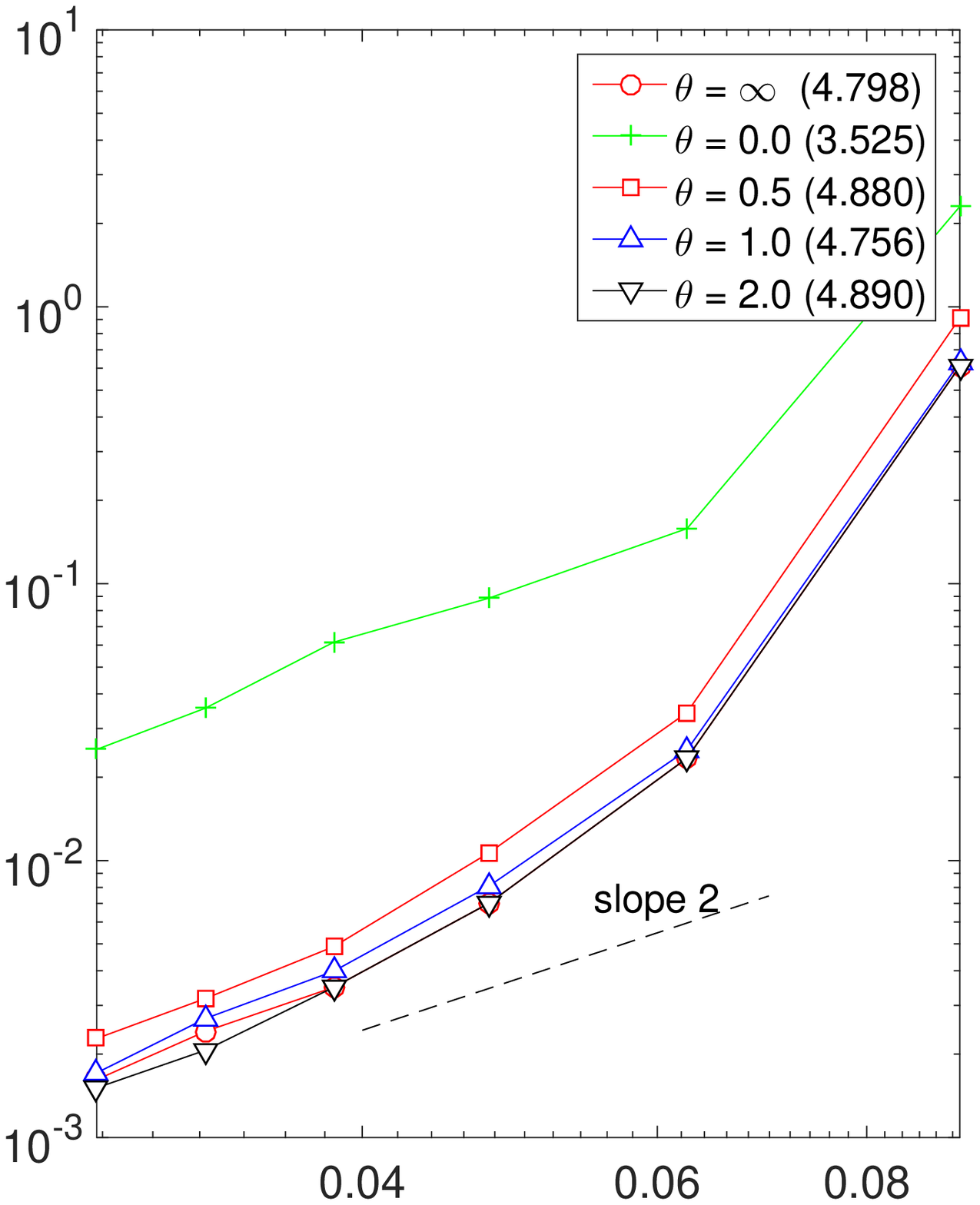} 
  \put (15,85) {$\calu_{Z\cup Y}$ }
  \put (10,95) {$\call u =  \triangle u + [2,\,3]^T\nabla u -4 u$}
  \end{overpic}
  &
  \begin{overpic}[width=0.4\textwidth,tics=10]{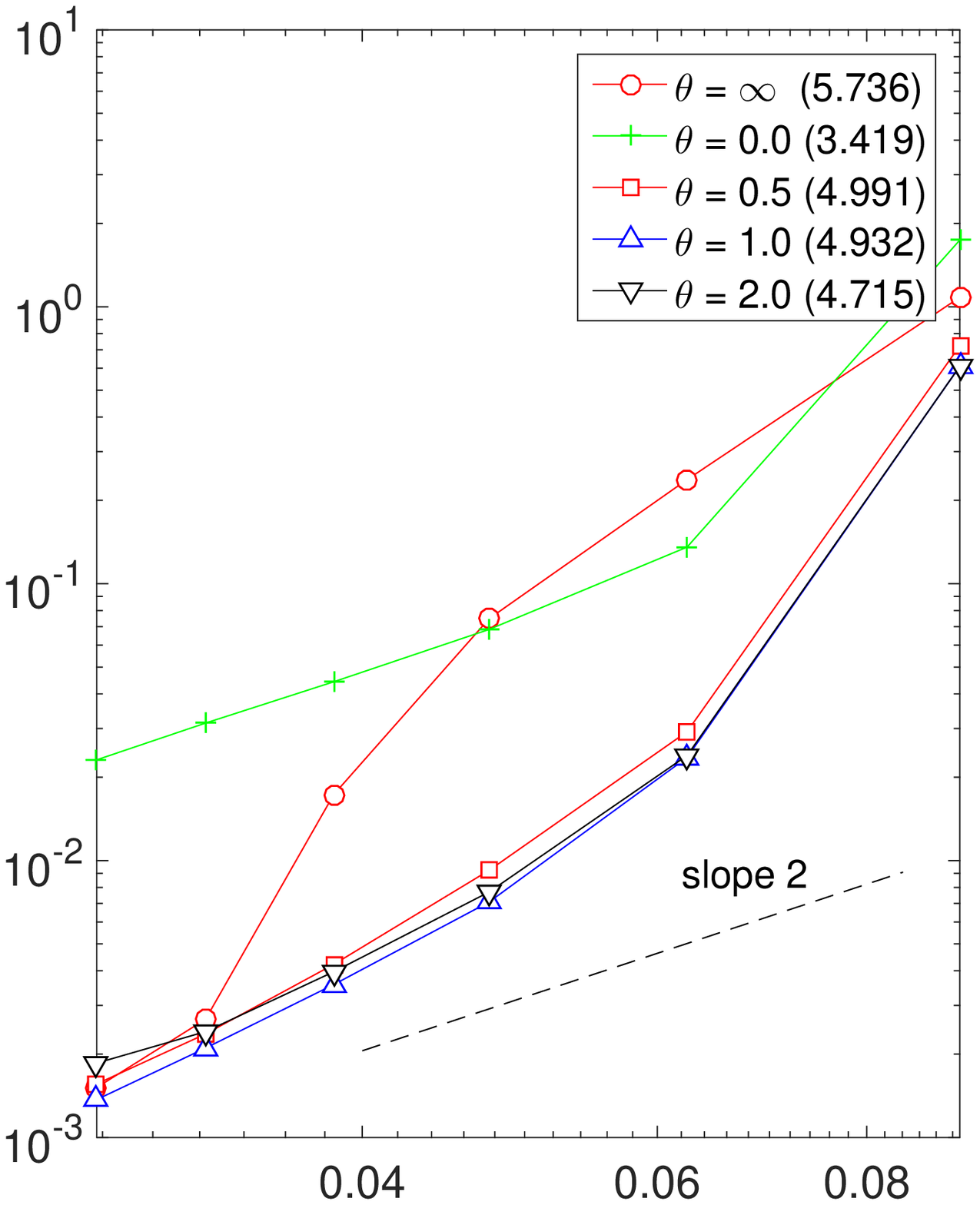} 
  \put (15,85) {$\calu_{Z}$ }
  \end{overpic}
\\
  \begin{overpic}[width=0.4\textwidth,tics=10]{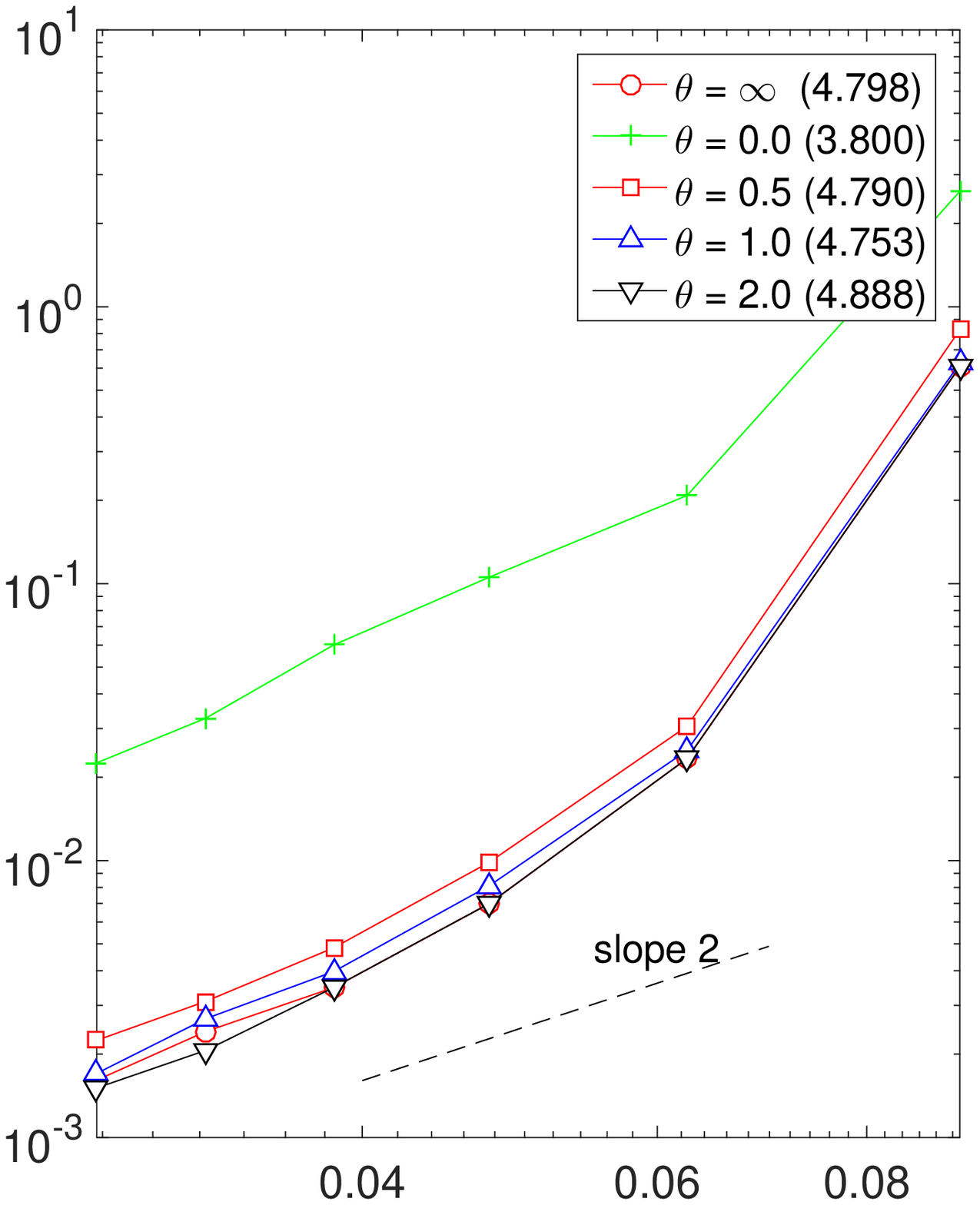} 
  \put (15,85) {$\calu_{Z\cup Y}$ }
  \put (10,95) { $\call u =  \triangle u + (x^2+1)u $}
  \end{overpic}
  &
  \begin{overpic}[width=0.4\textwidth,tics=10]{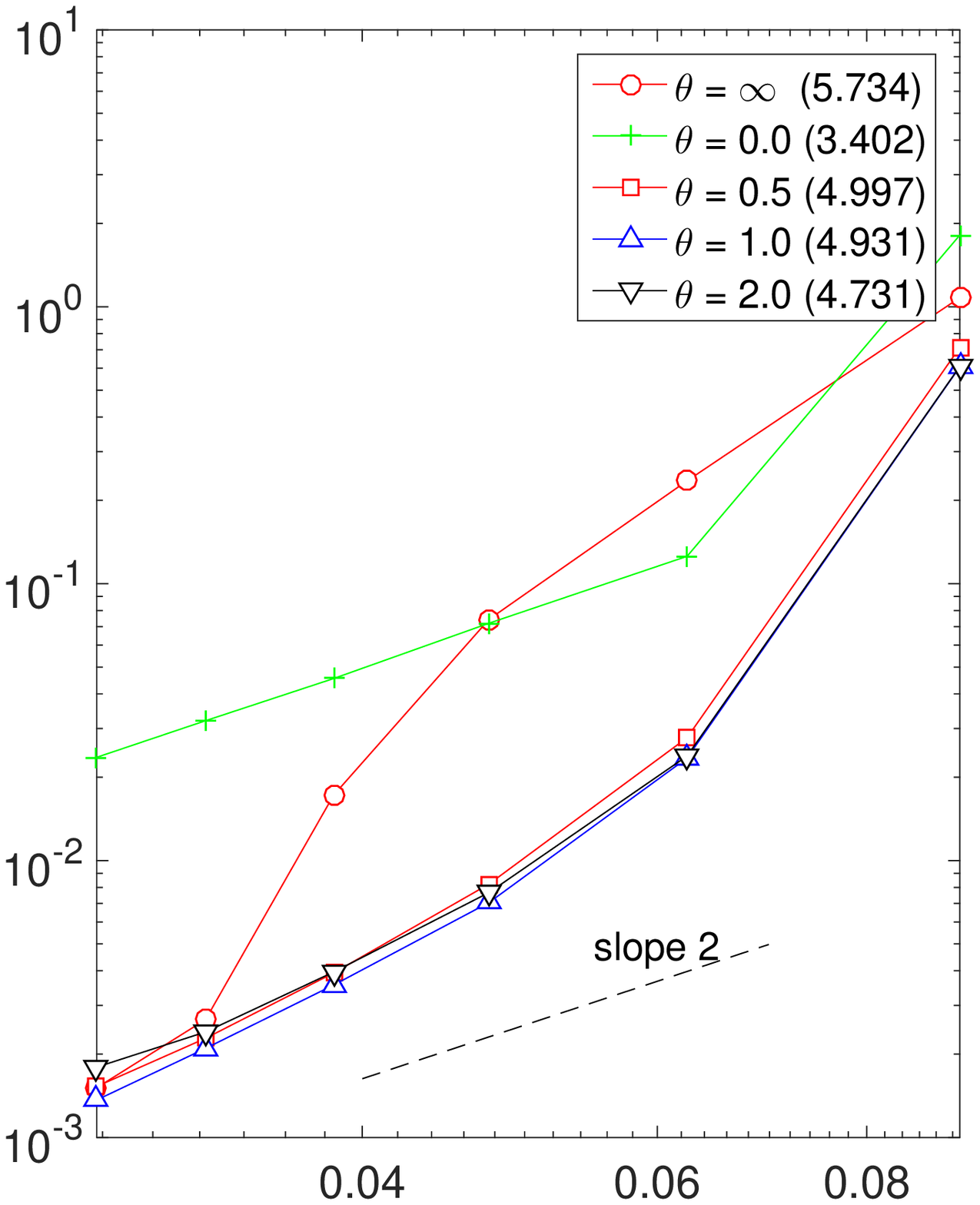} 
  \put (15,85) {$\calu_{Z}$ }
  \end{overpic}
\\
  \begin{overpic}[width=0.4\textwidth,tics=10]{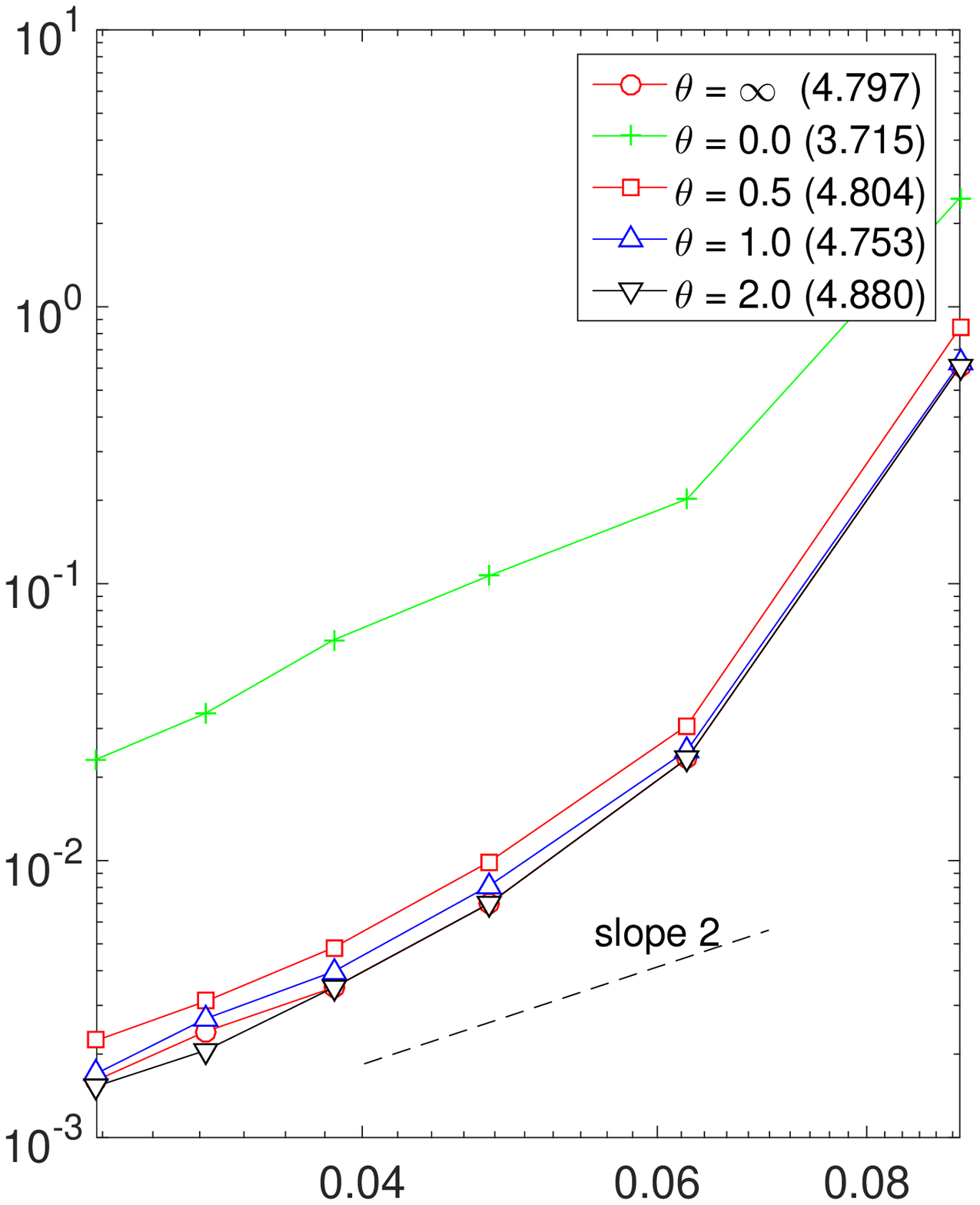} 
  \put (15,85) {$\calu_{Z\cup Y}$ }
  \put (10,95) { $\call u =  \triangle u + xu $}
  \end{overpic}
  &
  \begin{overpic}[width=0.4\textwidth,tics=10]{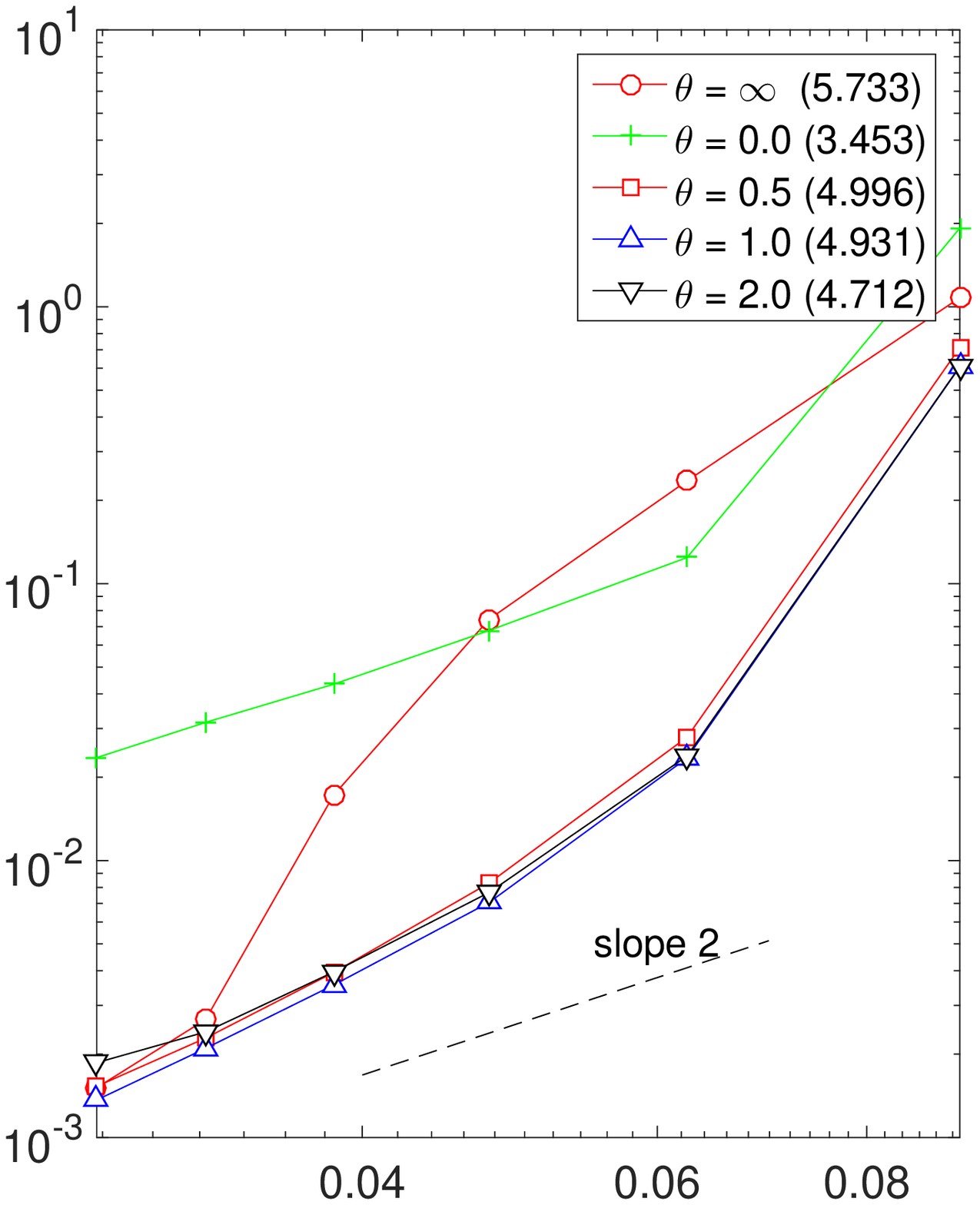} 
  \put (15,85) {$\calu_{Z}$ }
  \end{overpic}
\\
\end{tabular}
  \caption{Example~\ref{examp3}: $H^2(\Omega)$ error profiles for casting the WLS$(\theta)$ formulation in $\calu_{Z\cup Y}$ and $\calu_Z$ with Whittle-Mat\'{e}rn-Sobolev kernels of order $m=4$ to solve various PDEs with exact solution $u^*=\peaks(3x,3y)$.}\label{fig:3bpeaks}
\end{figure}

To further verify these observations, we now use sets of $n_X$ scattered collocation points (generated by the Halton sequence; see Figure~\ref{fig:PDEpt.eps}) to solve three different PDEs. 
Boundary collocation points remain regular with $h_Y=h_Z/2$. We present the numerical result for $m=4$ in Figure~\ref{fig:3bpeaks}. All PDEs have \verb"peaks(3x,3y)" as exact solution. These results should be compared to the results of $m=4$ in Figure~\ref{fig:3apeaks} and the convergence patterns are very similar. Yet, there are some minor, but notable differences.
Based on these test results, the WLS formulation in $\calu_Z$ with  moderate weighting (i.e. $\theta=0.5$ or $1$), which agrees with that in \cite{Hu+ChenETAL-Weigradibasicoll:07}, is relatively  stable in the numerical sense and is as accurate as the CLS formulation in $\calu_{Z\cup Y}$.

\begin{figure}
\centering
\begin{tabular}{cc}
  \begin{overpic}[width=0.4\textwidth,tics=10]{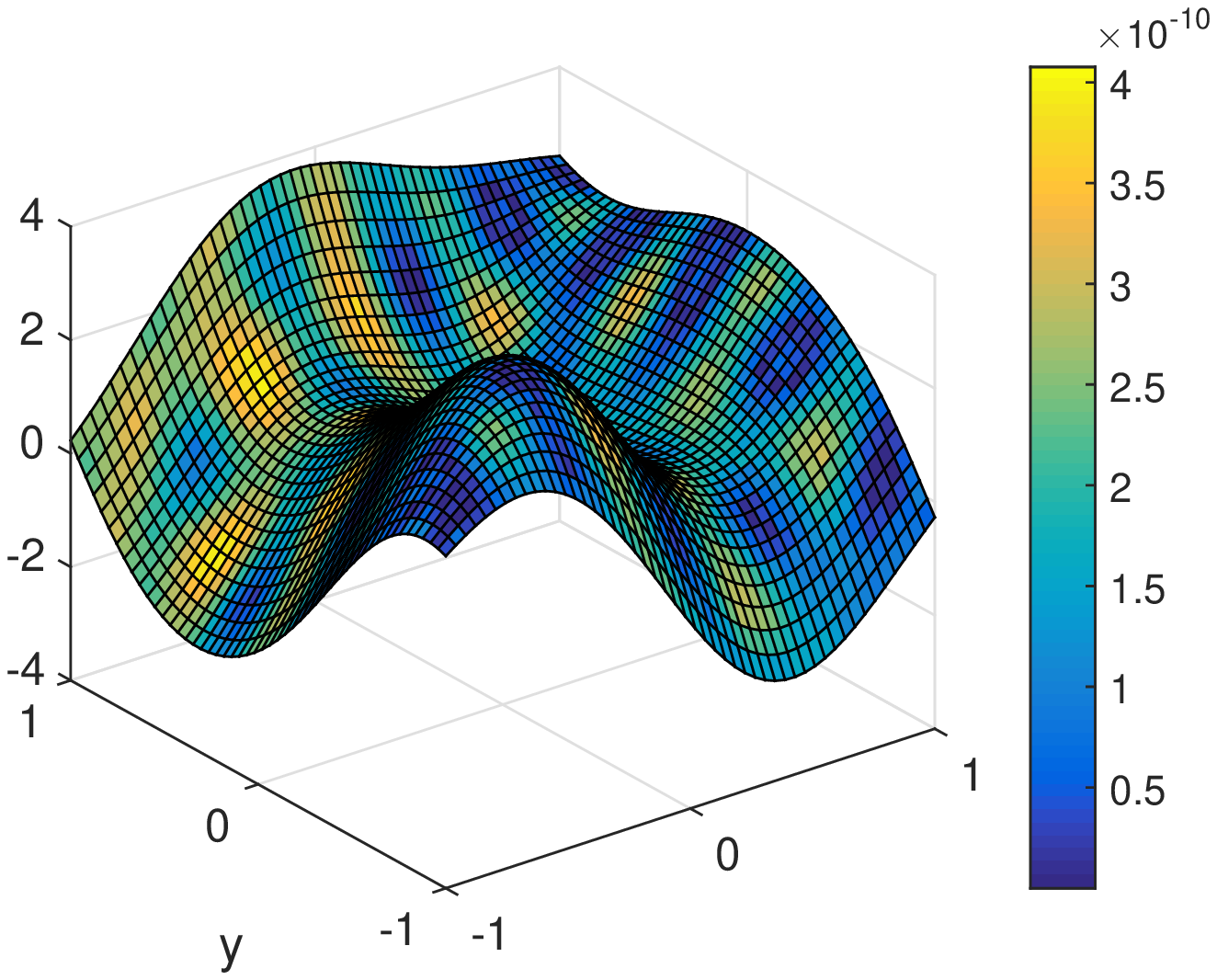} 
  \put (10,72) {$u^*=peaks(x,y)$ }
  \end{overpic}
  &
  \begin{overpic}[width=0.4\textwidth,tics=10]{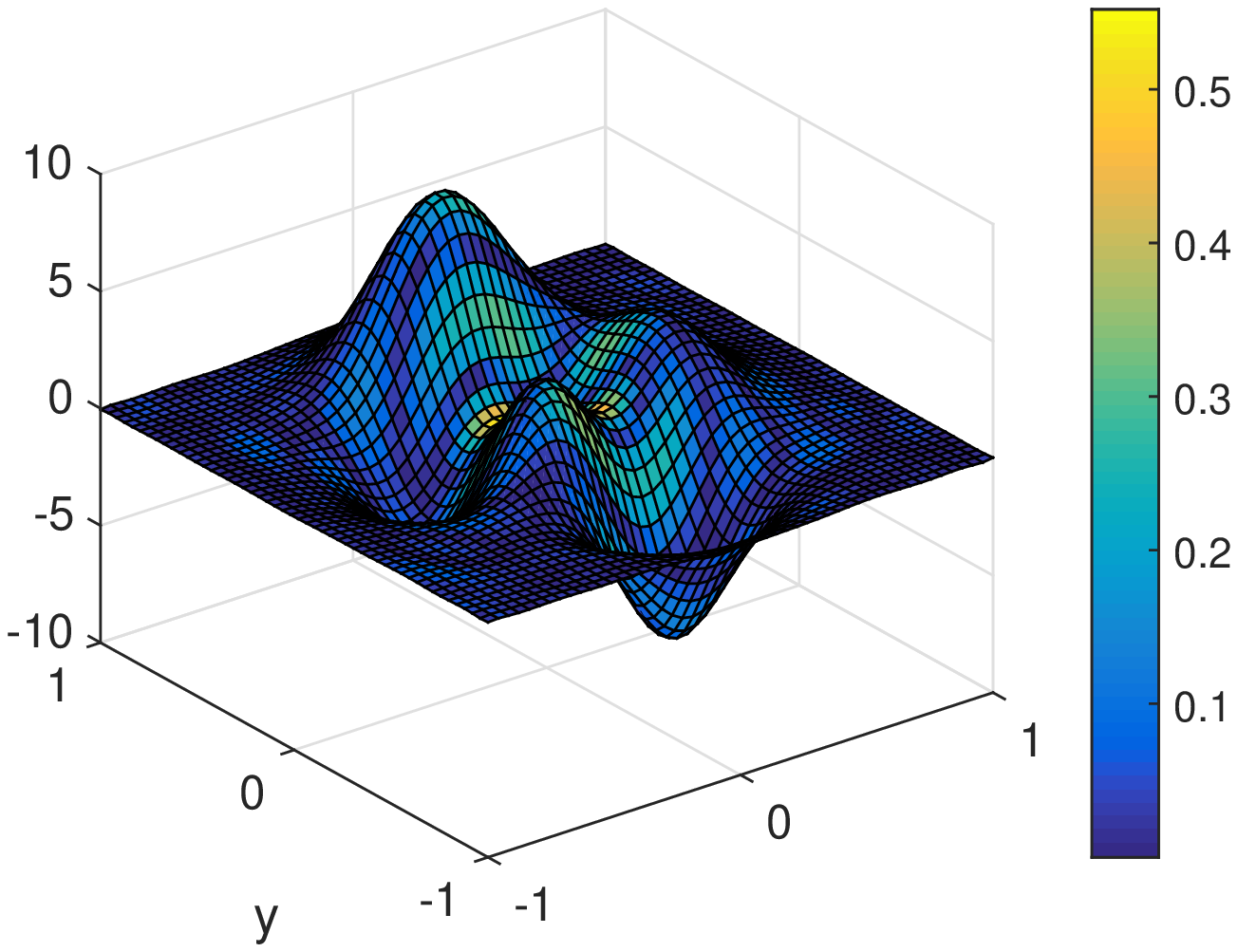} 
  \put (10,72) {$u^*=peaks(3x,3y)$}
  \end{overpic}
\end{tabular}
  \caption{Numerical approximations from CLS in $\calu_Z$ using the unscaled Gaussian kernel for solving a Poisson problem with exact solutions zoom-in peak (left) and full peak (right).}\label{fig:GACLSzoompeak}
\end{figure}

\begin{figure}
\centering
\begin{tabular}{cc}
  \begin{overpic}[width=0.4\textwidth,tics=10]{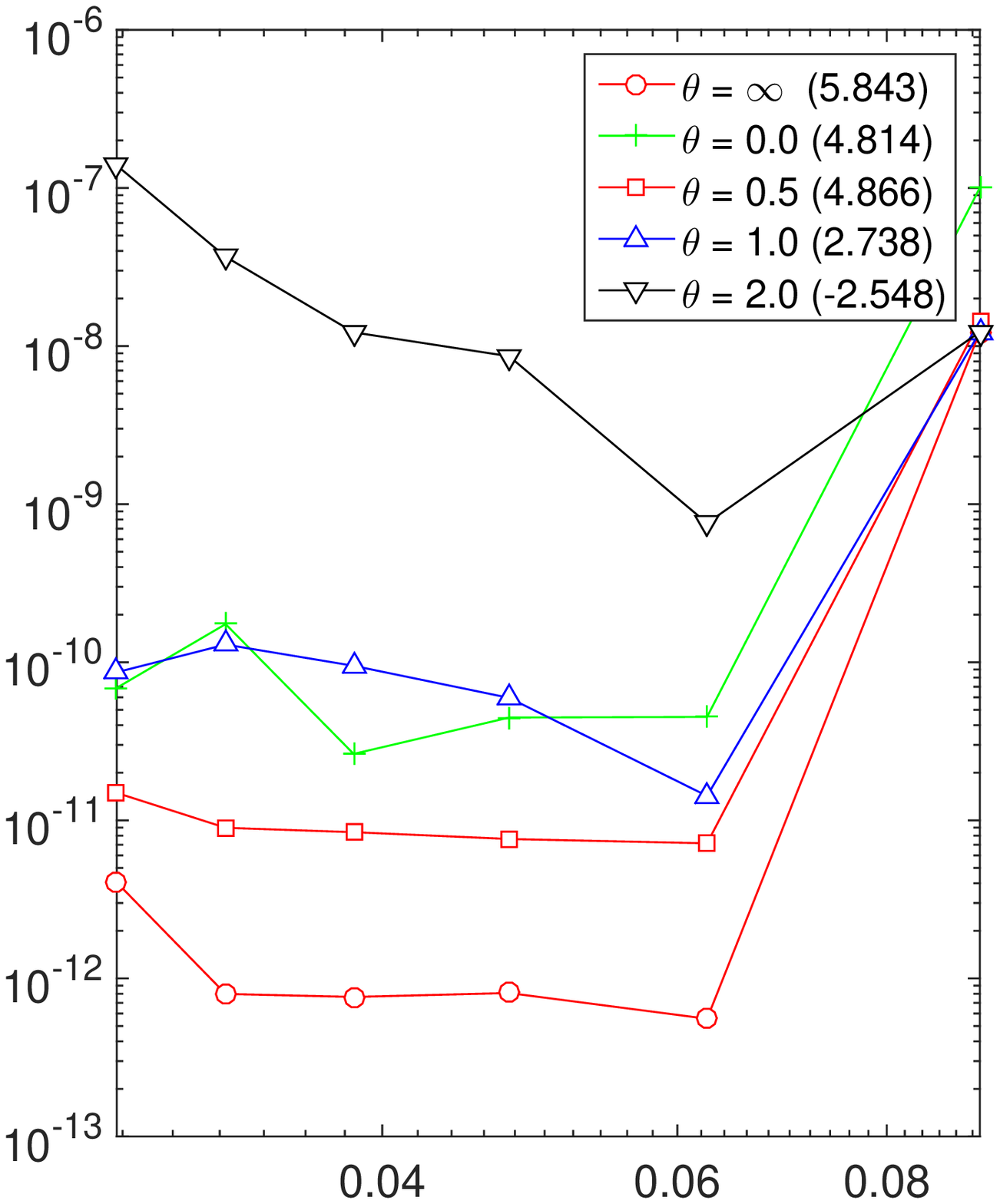} 
  \put (10,95) {$u^*=peaks(x,y)$ }
  \put (65,15) {GA}
  \end{overpic}
  &
  \begin{overpic}[width=0.4\textwidth,tics=10]{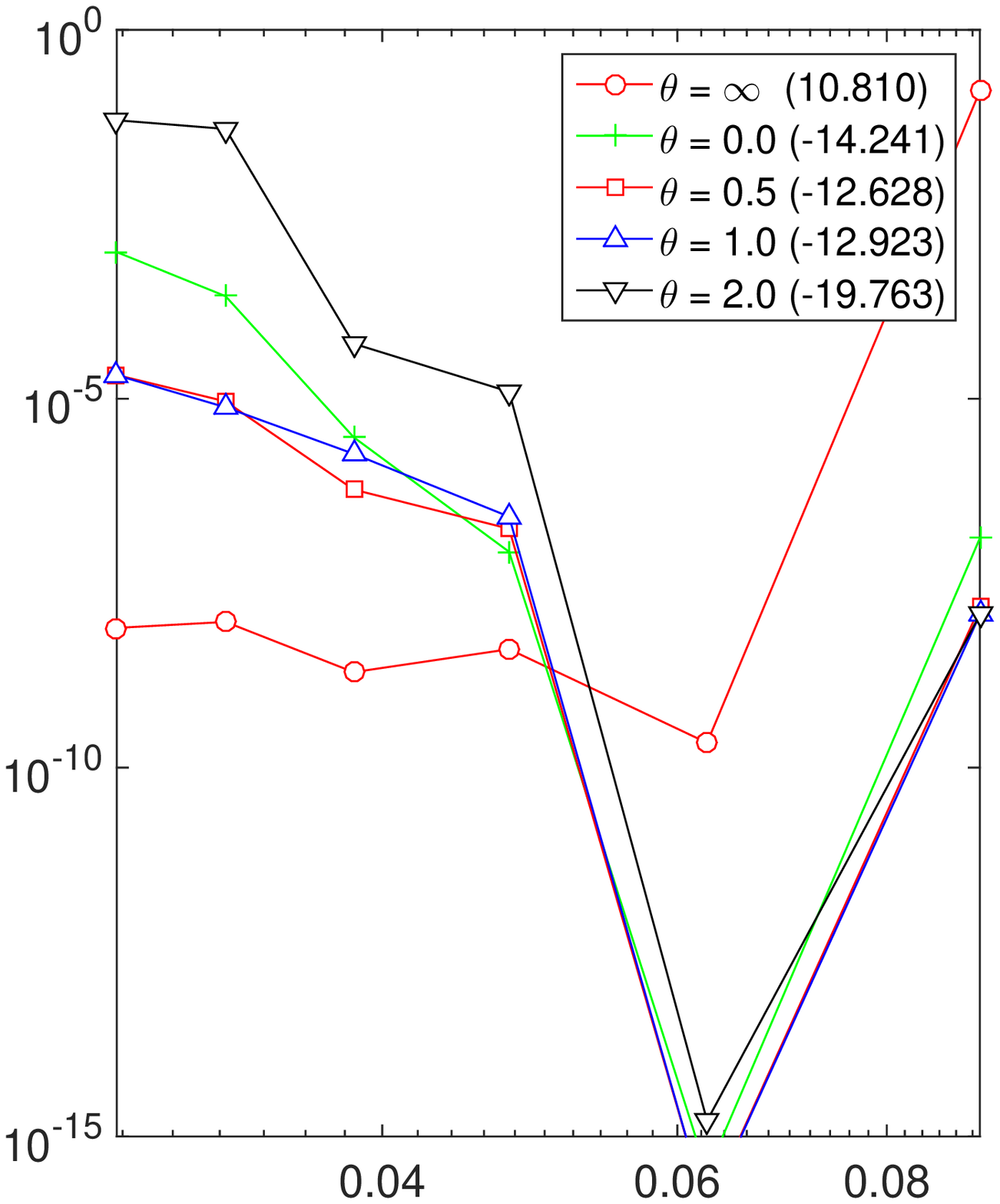} 
  \put (50,15) {RBF-QR}
  \end{overpic}
\\
  \begin{overpic}[width=0.4\textwidth,tics=10]{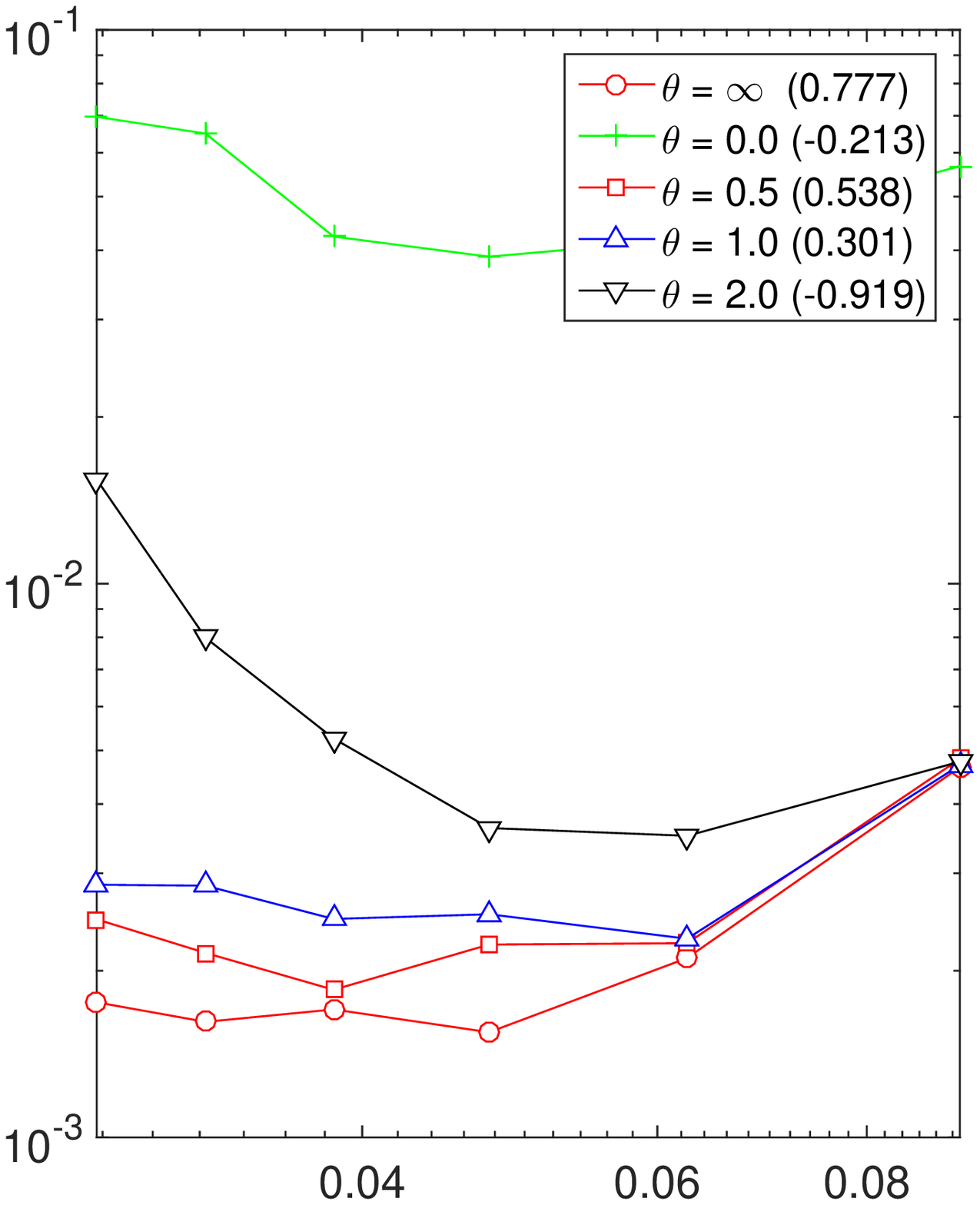} 
  \put (10,95) {$u^*=peaks(3x,3y)$ }
  \put (65,15) {GA}
  \end{overpic}
  &
  \begin{overpic}[width=0.4\textwidth,tics=10]{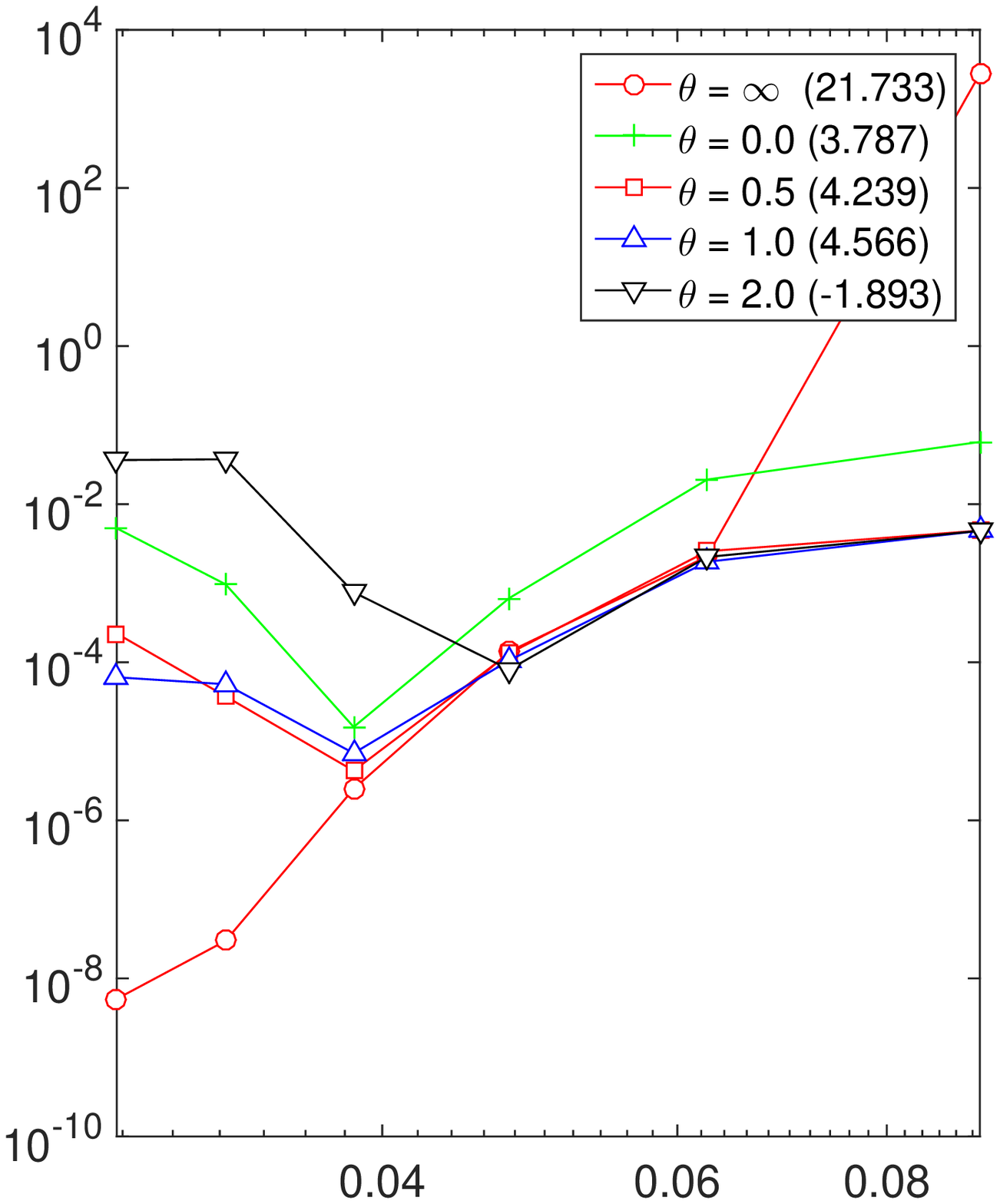} 
  \put (50,15) {RBF-QR}
  \end{overpic}
\end{tabular}
  \caption{Example~\ref{examp4}: $L^2(\Omega)$ error profiles for casting the CLS formulation in  $\calu_Z$ with unscaled GA kernels to solve $\Delta u=f$ with different exact solution.}\label{fig:GAQR}
\end{figure}

\smallskip
\begin{example}\label{examp4}
  \textbf{Some observations for GA and MQ.}
\end{example}
Let us begin with two CLS approximations obtained by using the unscaled Gaussian kernel, see Figure~\ref{fig:GACLSzoompeak}. 
We generate the trial space $\calu_Z$ using $n_Z=36^2$ centers. Collocations are enforced on regular $X$ and $Y$ with $h_X=h_Y=h_Z/2$. We then solve two Poisson problems with different exact solutions $u^*=peaks(x,y)$ and $u^*=peaks(3x,3y)$. Both systems share the same Kansa matrix  and only differ by the right-hand vectors $f_{|X}$ and $g_{|Y}$.
Although both systems have exactly the same condition number, the resulting
accuracy has a huge difference: 8.1E-12 and 3.5E-3 $L^2(\Omega)$ error
respectively, for the truncated and full peaks.
One could argue that the Gaussian native space is relatively small and it may not contain the full peak function that causes low accuracy. To know for sure, a change of stable basis (a.k.a. RBF-QR algorithm \cite{Fasshauer+McCourt-StabEvalGausRadi:12,Fornberg+LarssonETAL-StabCompwithGaus:11}) can give us a clearer picture. Using some downloadable Matlab  codes \cite{Larsson-RADIBASIFUNCINTE}, we can cast the CLS formulation completely in the new stable basis and yield new $L^2(\Omega)$ errors,  5.7E-8 and 6.3E-9 respectively for the two exact solutions.
For the zoom-in peak, accuracy drops, which can be explained by the truncation error within the RBF-QR algorithm. Whereas, the huge accuracy improvement suggests there is a good candidate in the native space of Gaussian to approximate the full peak function. In Figure~\ref{fig:GAQR}, we show the error profiles for direct Gaussian (GA) and stable basis (RBF-QR) for both solutions.
When solving the zoom-in peak, GA can provide highly accurate approximation and its error stagnates for large $n_Z$; adding RBF-QR (as is without modification) to the algorithm will introduce numerical instability. The situation is very different in the full peak when GA fails to converge; having RBF-QR brings convergence back in the game.
We remark that the RBF-QR algorithm prefers {a} 
flatter (than our unscaled) Gaussian basis for both accuracy and efficiency.
Moreover, some of the highly accurate approximations in RBF-QR suggest that it is highly possible to properly truncate the Gaussian expansion in order to couple with the CLS and WLS approaches. This is out of the scope of this work, but worth some further investigation.

\begin{figure}
\centering
\begin{tabular}{cc}
  \begin{overpic}[width=0.4\textwidth,tics=10]{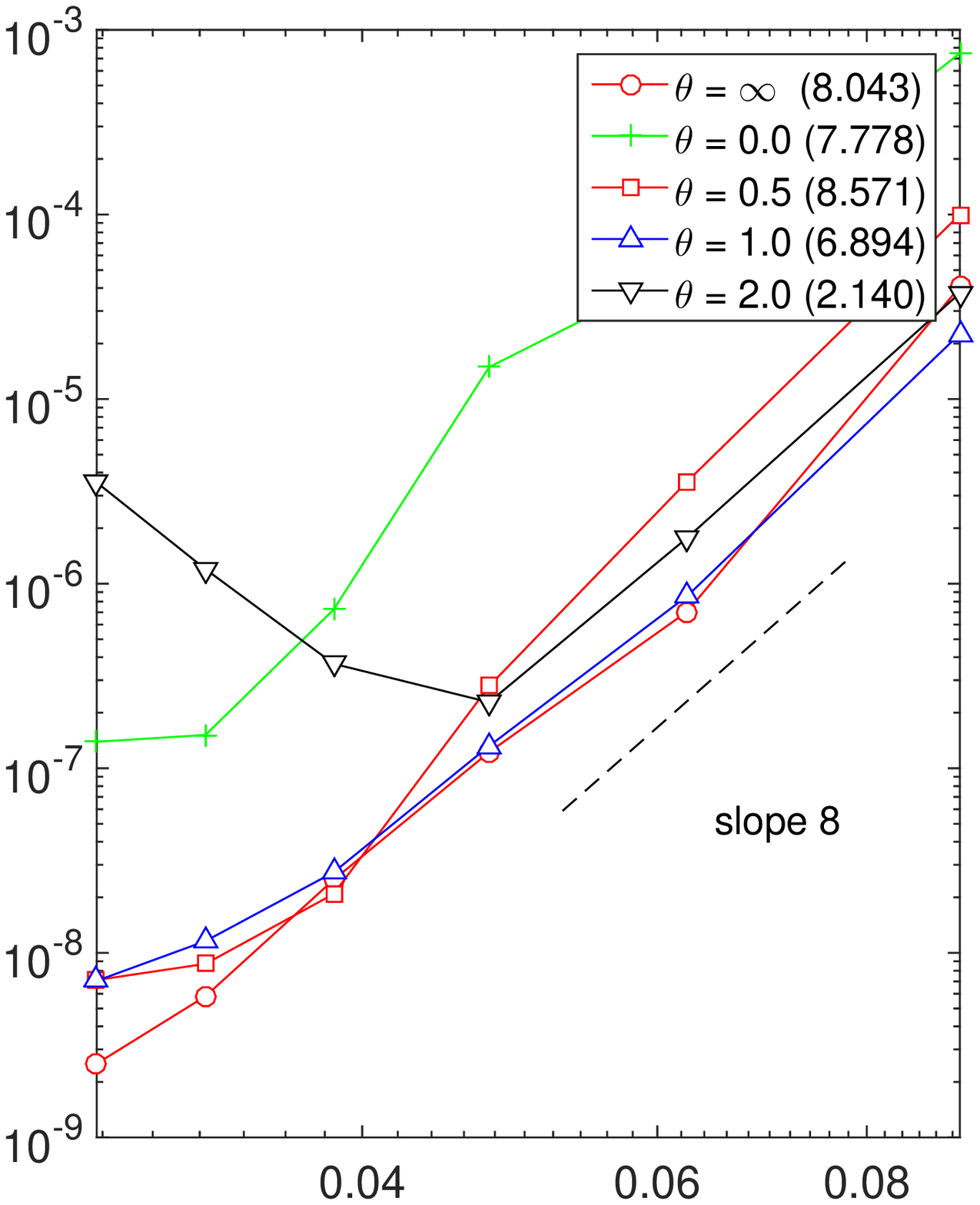} 
  \put (10,95) {$u^*=peaks(x,y)$ }
  \put (40,15) {MQ in $\calu_{Z\cup Y}$}
  \end{overpic}
  &
  \begin{overpic}[width=0.4\textwidth,tics=10]{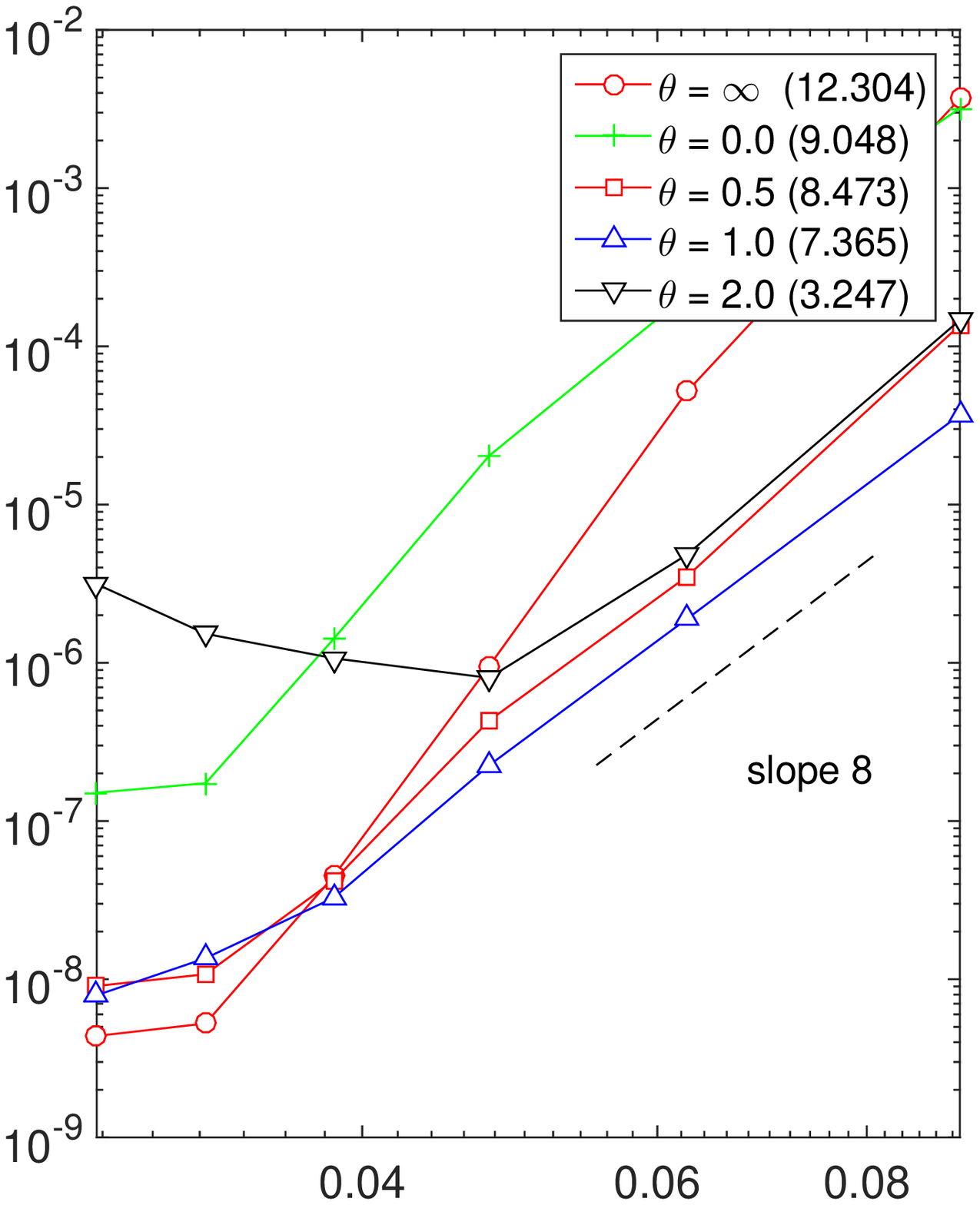} 
  \put (45,15) {MQ in $\calu_{Z}$}
  \end{overpic}
\\
  \begin{overpic}[width=0.4\textwidth,tics=10]{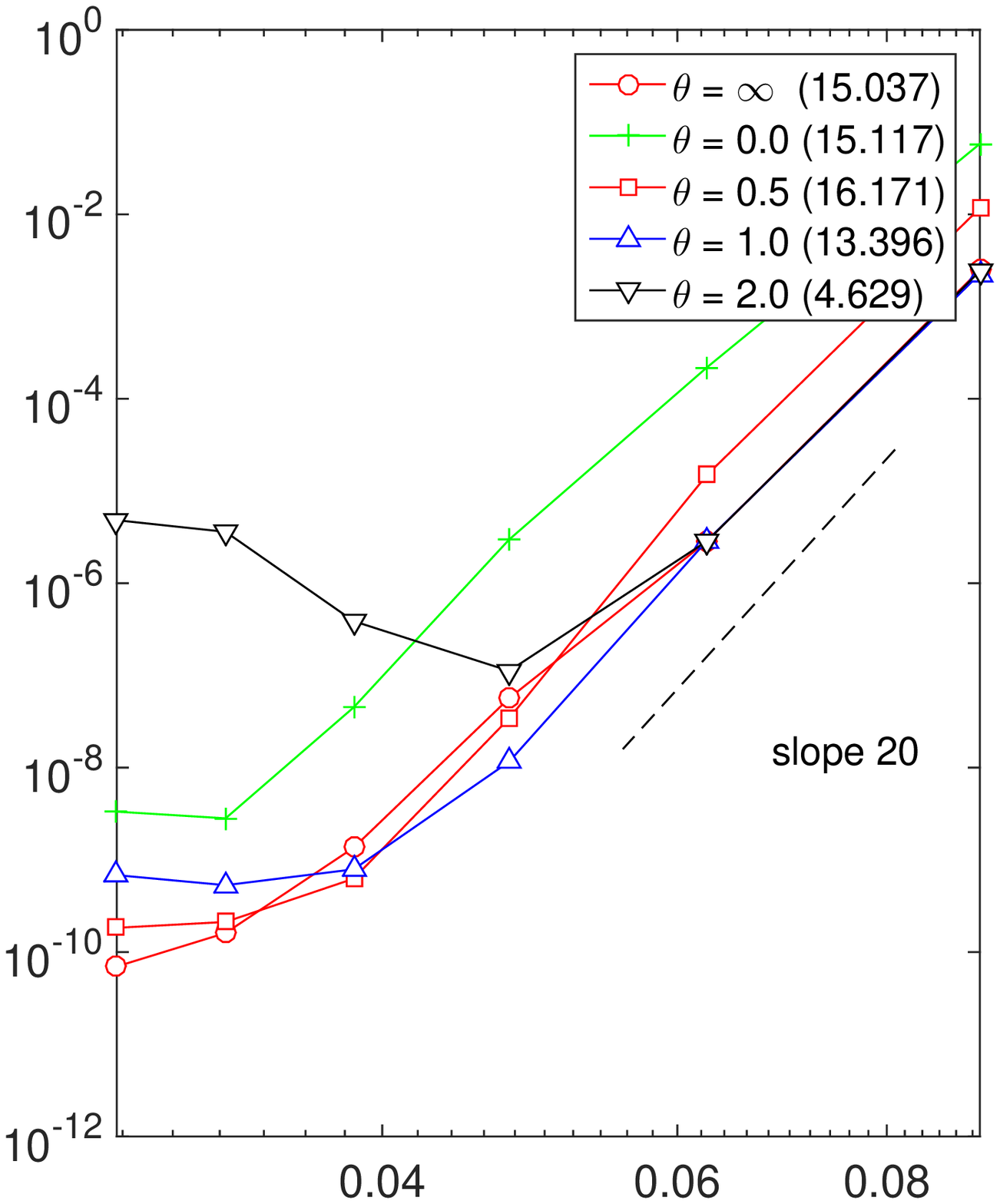} 
  \put (10,95) {$u^*=peaks(3x,3y)$ }
  \put (40,15) {MQ in $\calu_{Z\cup Y}$}
  \end{overpic}
  &
  \begin{overpic}[width=0.4\textwidth,tics=10]{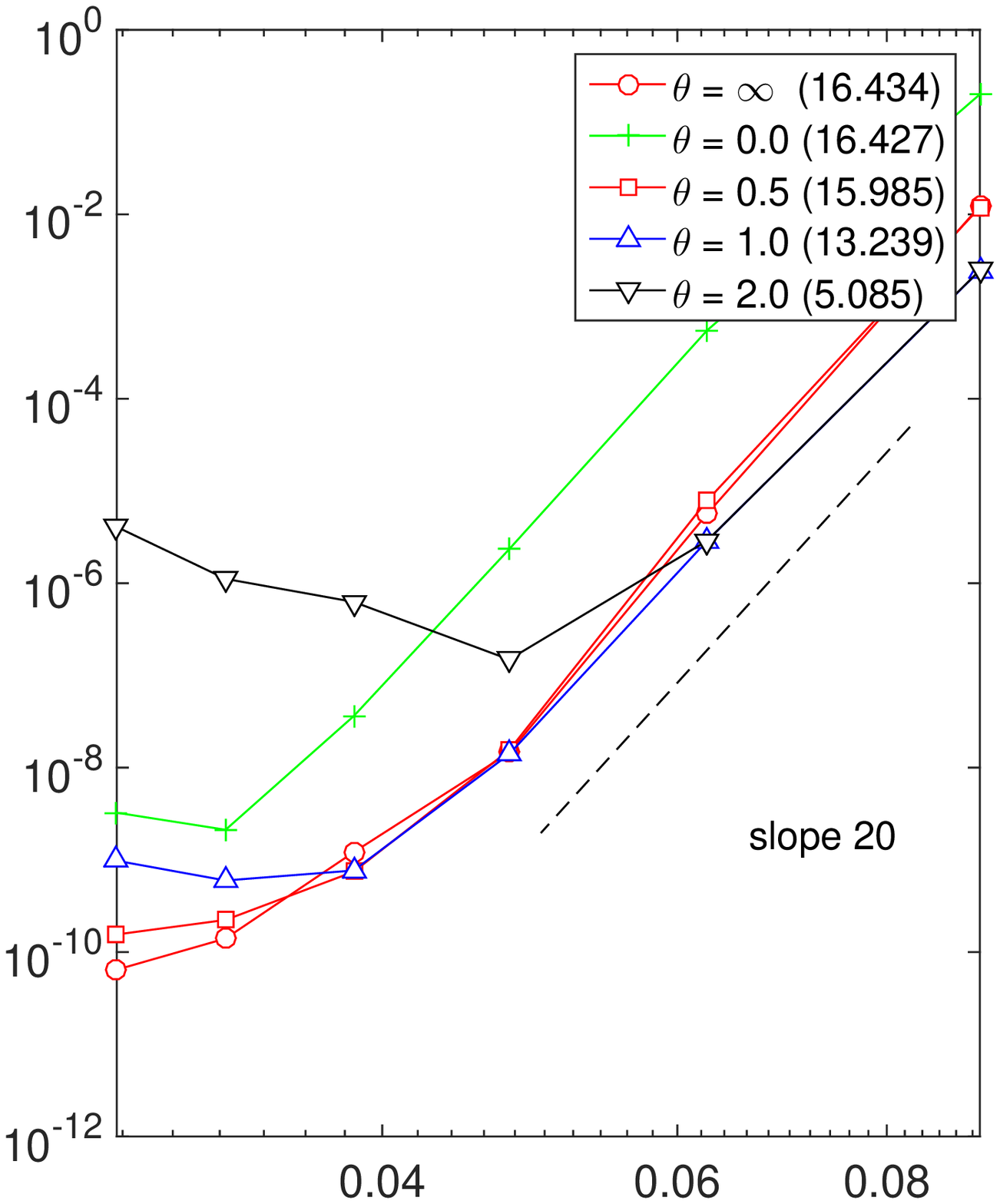} 
  \put (45,15) {MQ in $\calu_{Z}$}
  \end{overpic}
\end{tabular}
  \caption{Example~\ref{examp4}: $:L^2(\Omega)$ error profiles for casting the CLS formulation in $\calu_Z$ with unscaled MQ kernels to the same settings as in  Figure~\ref{fig:GAQR}.}\label{fig:MQ}
\end{figure}

Our demonstration will end with the results of the unscaled multiquadrics (MQ) kernel in Figure~\ref{fig:MQ}. While the full peak function  is giving the Gaussian kernel trouble, the multiquadrics is doing very well. 
Despite so, all CLS, WLS($0.5$) and WLS($1$) with the MQ basis converge to the full peaks at an estimated rate of 20. Also, the convergence-lag for the CLS in $\calu_Z$ is not at all noticeable.
Turning to the less varying zoom-in peak,  the same MQ-PDE resultant matrix only yields an 8th-order  convergence.
This example again confirms  that {the} 
condition number alone cannot be used to
predetermine the  accuracy or convergence rate of Kansa related methods using 
the GA
or MQ basis.



\section*{Conclusion}

We prove some error estimates for a constraint least-squares and a class of weighted least-squares strong-form RBF collocation formulations for solving general second order elliptic problem with nonhomogenous Dirichlet boundary condition. All analysis is carried out in Hilbert spaces so that both PDE and RBF theories apply.
We show that the CLS  and WLS formulations using kernels, 
{which reproduce $H^m(\Omega)$}, with sufficient smoothness  can converge at the optimal ${m-2}$ rate in $H^2(\Omega)$.
Besides some standard smoothness assumptions for high order convergence, the
sets of  collocation points have to 
satisfy some denseness conditions for the convergence theories to hold.

We verify by numerical examples that {there are}
many convergent formulations for $\Omega\subset\R^2$ 
{that} enjoy the optimal convergence rate.
We thoroughly study the numerical performance of Whittle-Mat\'{e}rn-Sobolev kernels in two trial spaces. The larger space that includes all boundary collocation points as trial centers is more theoretically sound
{(in the sense of the range of optimal weighting),}
whereas the small one is computationally more efficient. Taking both accuracy
and efficiency into consideration, casting WLS in the small trial space with a
moderate weight 
consistently yields competitive
accuracy and numerical stability. This recommendation  extends to
the commonly used Gaussian and multiquadrics kernels, which do not reproduce
Sobolev space as the theories require. We also provide a brief demonstration for
using the RBF-QR algorithm on our formulations to hint at possibilities for
future {research}.

\section*{Acknowledgements}
This work was supported by a GRF Grant from the Hong Kong Research Grant Council and an FRG Grant from Hong Kong Baptist University.

\bibliographystyle{amsplain} 
\def\cprime{$'$} \def\cprime{$'$} \def\cprime{$'$}
\providecommand{\bysame}{\leavevmode\hbox to3em{\hrulefill}\thinspace}
\providecommand{\MR}{\relax\ifhmode\unskip\space\fi MR }
\providecommand{\MRhref}[2]{%
  \href{http://www.ams.org/mathscinet-getitem?mr=#1}{#2}
}
\providecommand{\href}[2]{#2}

\newpage

\begin{figure}
  \centering
\begin{tabular}{cc}
  \begin{overpic}[width=0.4\textwidth,tics=10]{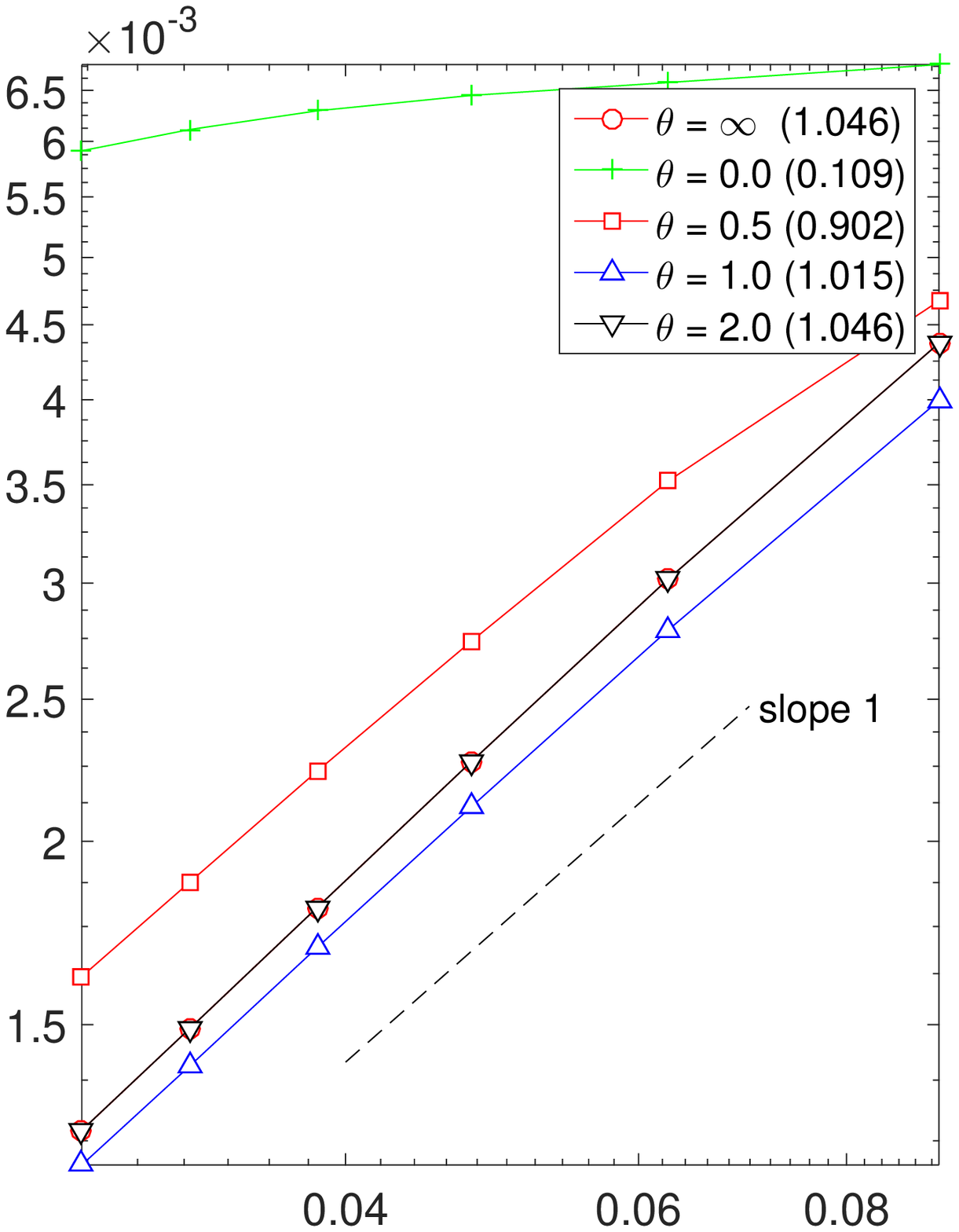} 
  \put (15,85) {$\calu_{Z\cup Y}$ }
  \put (15,100) {$m=3$}
  \end{overpic}
  &
  \begin{overpic}[width=0.4\textwidth,tics=10]{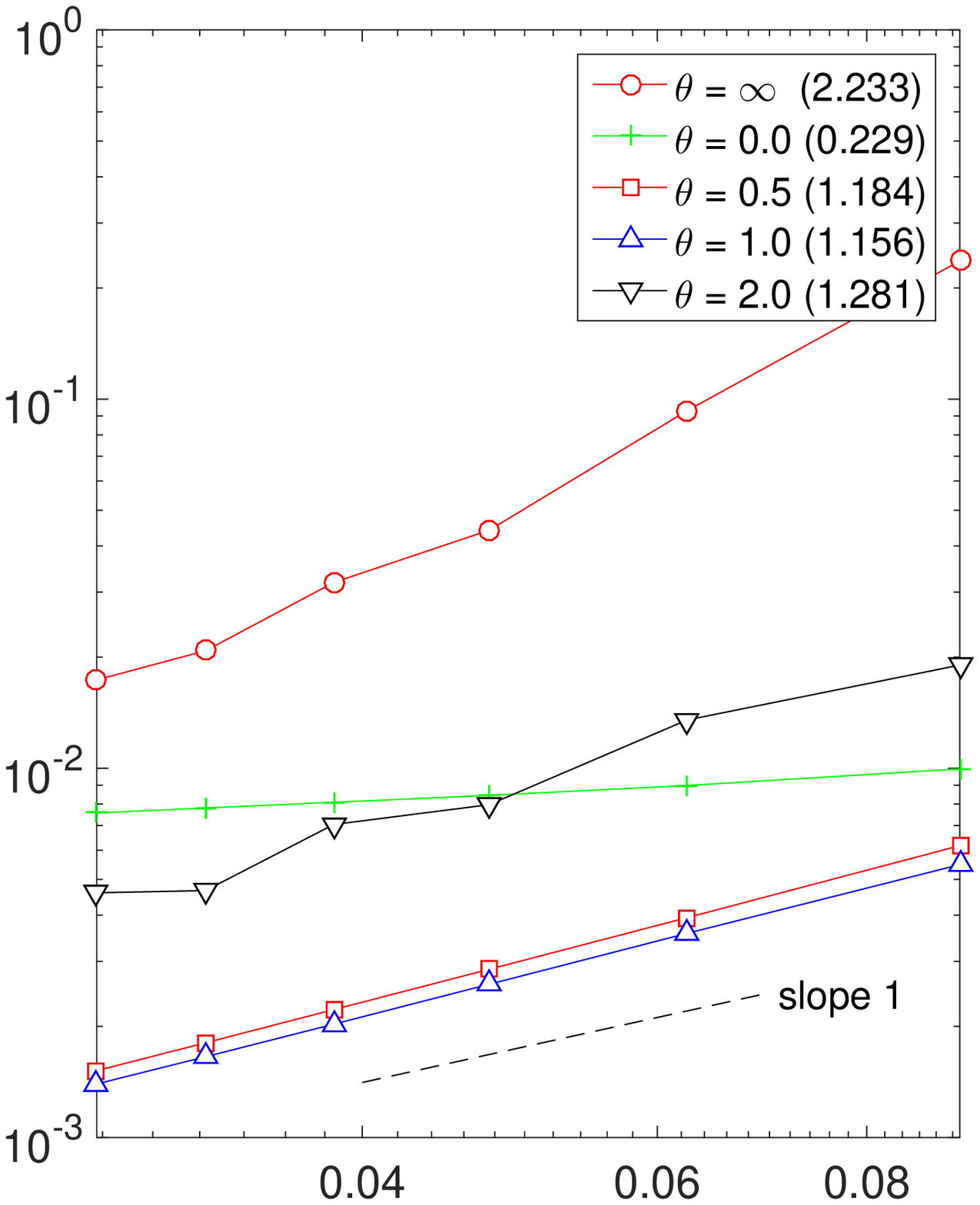} 
  \put (15,85) {$\calu_{Z}$ }
  \end{overpic}
\\
  \begin{overpic}[width=0.4\textwidth,tics=10]{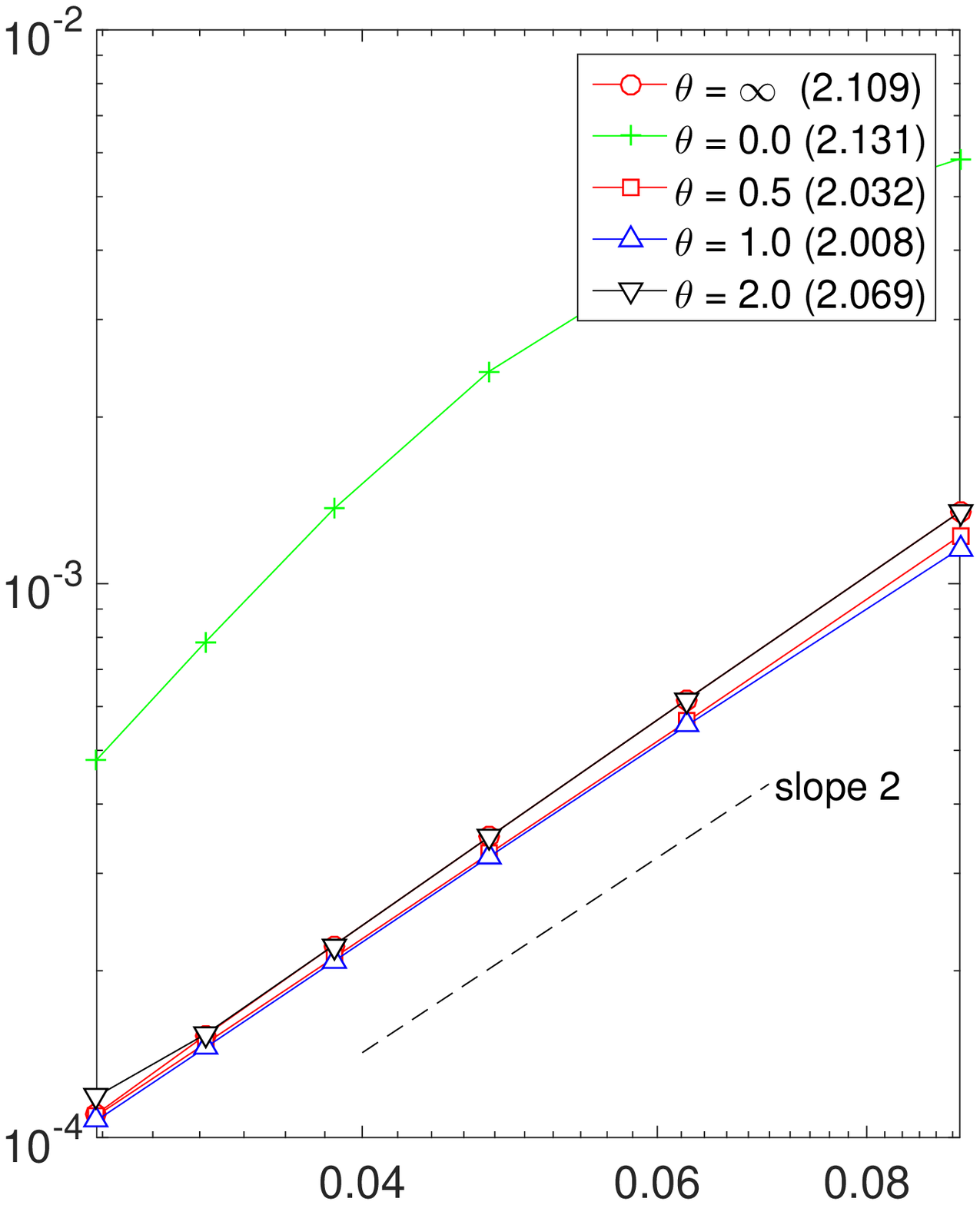} 
  \put (15,85) {$\calu_{Z\cup Y}$ }
  \put (15,95) {$m=4$}
  \end{overpic}
  &
  \begin{overpic}[width=0.4\textwidth,tics=10]{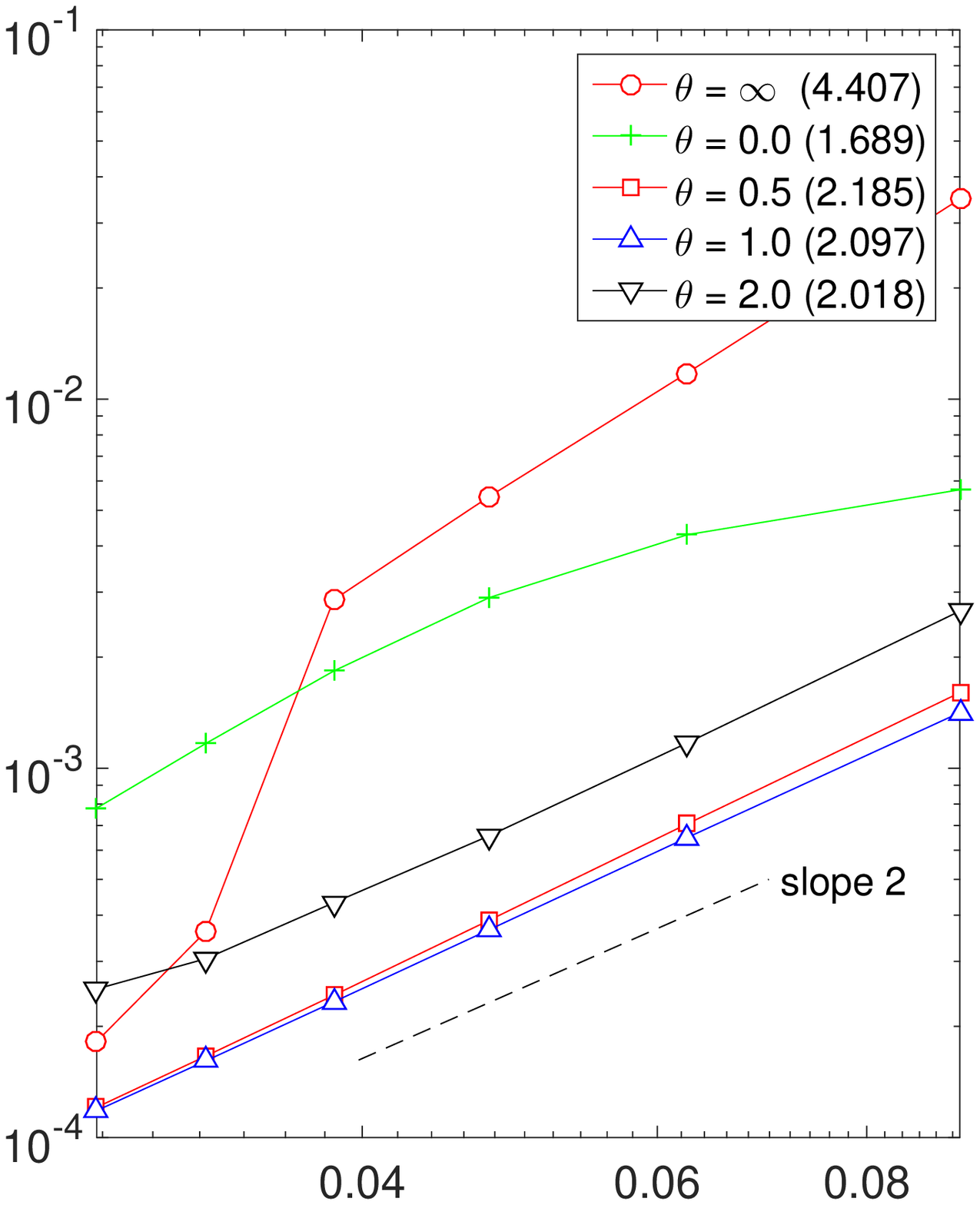} 
  \put (15,85) {$\calu_{Z}$ }
  \end{overpic}
\\
  \begin{overpic}[width=0.4\textwidth,tics=10]{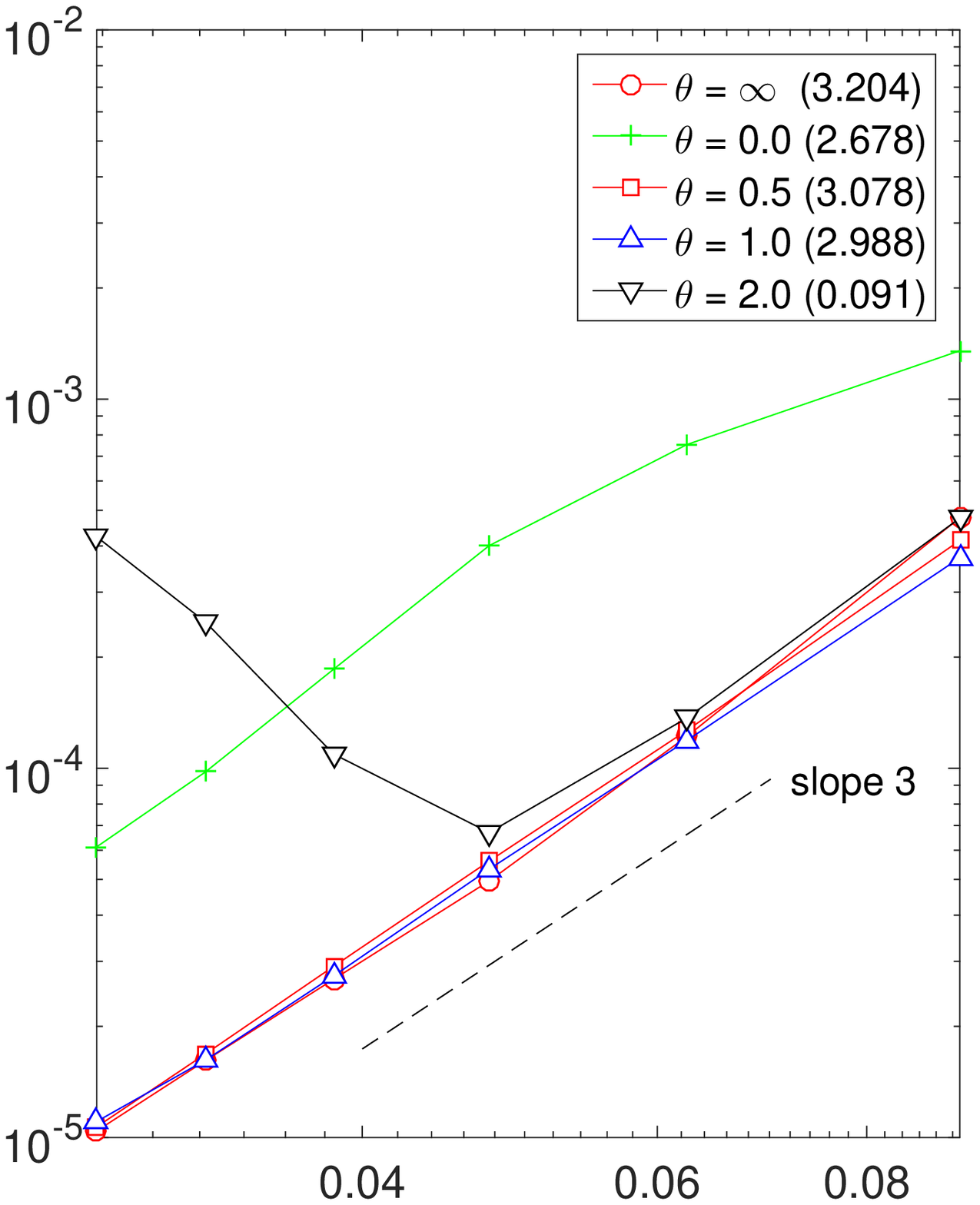} 
  \put (15,85) {$\calu_{Z\cup Y}$ }
  \put (15,95) {$m=5$}
  \end{overpic}
  &
  \begin{overpic}[width=0.4\textwidth,tics=10]{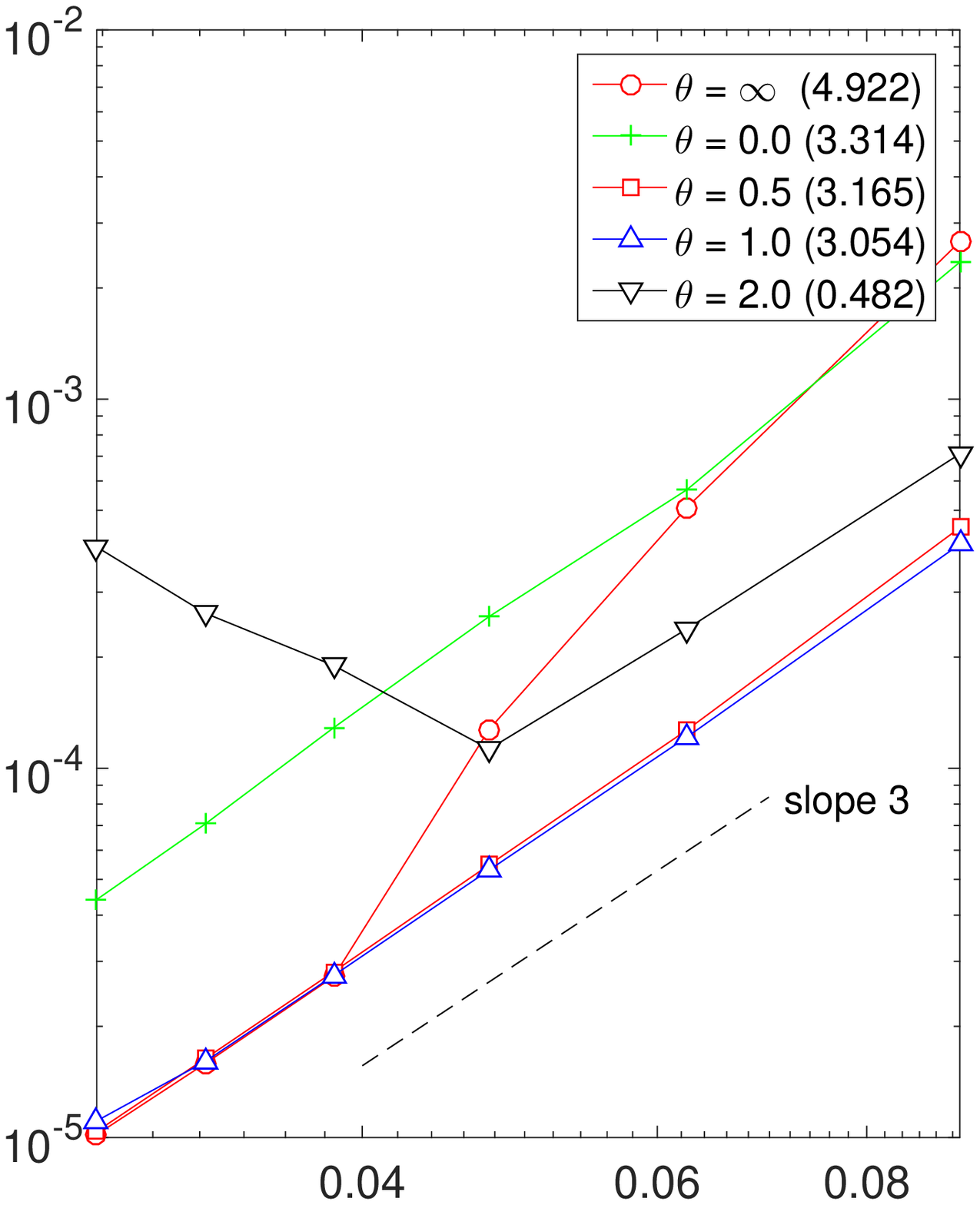} 
  \put (15,85) {$\calu_{Z}$ }
  \end{overpic}
\end{tabular}
  \caption{Example~\ref{examp3}: $H^2(\Omega)$ error profiles for casting the WLS$(\theta)$ formulation in $\calu_{Z\cup Y}$ and $\calu_Z$ with Whittle-Mat\'{e}rn-Sobolev kernels of order $m=3,4,5$ to solve $\Delta u = f$ with exact solution $u^*=\sin(\pi x/2,\pi y/2)$.}\label{fig:3asc}
\end{figure}
%

\end{document}